\documentclass{article}
\usepackage{amssymb}

\input{tcilatex}
\begin{document}

\bigskip \bigskip 

\bigskip \bigskip \bigskip 

\begin{center}
{\LARGE Two Theorems on the Range of Strategy-proof Rules on a Restricted
Domain}\bigskip 

{\Large Donald E. Campbell }$\cdot ${\Large \ Jerry S. Kelly}\bigskip 

\begin{eqnarray*}
&&\text{\textbf{Abstract}} \\
&&\text{Let }g\text{ be a strategy-proof rule on the domain }NP\text{ of} \\
&&\text{profiles where no alternative Pareto-dominates any} \\
&&\text{other and let }g\text{ have range }S\text{ on }NP\text{. \ We
complete} \\
&&\text{the proof of a Gibbard-Satterthwaite result - if S} \\
&&\text{contains more than two elements, then g is dictatorial -} \\
&&\text{by establishing a full range result on two subdomains of }NP\text{.}
\end{eqnarray*}
\end{center}

\bigskip 

\bigskip 

\textbf{1. Introduction.}

\textbf{2. Notation.}

\textbf{3. N}$-$\textbf{Range Theorem: Part 1.}

\textbf{4. N}$-$\textbf{Range Theorem: Part 2.}

\textbf{5. M}$-$\textbf{Range Theorem.\bigskip }

\textbf{1. Introduction.}\medskip

\qquad In Campbell and Kelly (2010), we discussed losses due to manipulation
of social choice rules and gave an example of a non$-$dictatorial rule such
that, for any manipulation, no one else has a loss. Where all individual
preferences are strict, that means that for any manipulation, everyone
gains. We say such rules satisfy universally beneficial manipulation (UBM);
this property is a weakening of Gibbard-Satterthwaite individual stra
tegy-proofness.\medskip

\qquad In Campbell and Kelly (2014a), we present two characterizations of
UBM rules. The proof in that paper uses structural results from the fact
that a UBM rule $g$ must satisfy strategy-proofness on the subdomain $NP$ of
all profiles with the property that no alternative Pareto-dominates any
other.\medskip 

\qquad In Campbell and Kelly (2014b), we derive those structural results: A
rule $g$ that is strategy-proof on $NP$ and has range of at least three
alternatives must be dictatorial. That proof, in turn, required, for two
induction steps, the result that if strategy-proof $g$ is of full range on $%
NP$ it is also of full range on two special subdomains of $NP$. This paper
completes the analysis by proving these two range results.\medskip 

\textbf{2. Notation.}\medskip

\qquad We adopt all terminology and notation from Campbell and Kelly
(2014b). \ In particular, we take as given a finite set $X$ of alternatives
with $|X|$ $=m\geq 3$ and finite set $N=\{1,2,...,n\}$ of individuals with $%
n\geq 3$. A (strong) \textit{ordering} on $X$ is a complete, asymmetric,
transitive relation on $X$ and the set of all such orderings is $L(X)$. For $%
R\in L(X)$ and $Y\subset X$ let $R|Y$ denote the relation $R\cap Y\times Y$
on $y$, the \textit{restriction} of $R$ to $y$. A \textit{profile} $p$ is a
map from $n$ to $L(X)$, where $p=(p(1),p(2),...,p(n))$ and we write $x\succ
_{p(h)}y$ if individual $h$ strongly prefers $x$ to $y$ at profile $p$. The
set of all profiles is $L(X)^{N}$. Let $\wp $ be a nonempty subset of $%
L(X)^{N}$. For each subset $Y$ of $X$ and each profile $p$ in $\wp $, let $%
p|Y$ denote the restriction of profile $p\in \wp $ to $Y$. That is, $p|Y$
represents the function $q\in L(Y)^{N}$ satisfying $q(i)=p(i)|Y$ for all $%
i\in N$. A \textit{social choice rule} on $\wp $ is a function $g:\wp
\rightarrow X$. \ A rule $g$ is \textit{full-range} if Range($g$) $=X$, and $%
g$ is \textit{dictatorial} if there is an individual $i$ such that $g(p)$ is
the highest ranked element in Range($g$) according to ordering $p(i)$%
.\medskip 

\qquad In this paper, we consider social choice rules on the Non-Paretian
domain, $NP$, the set of all profiles $p$ such that for any pair of distinct
alternatives, $x$ and $y$, there exists an individual $i\in N$ such that $%
x\succ _{p(i)}y$ and there exists an individual $j\in N$ such that $y\succ
_{p(j)}x$, so that neither alternative is Pareto-superior to the
other.\medskip

\qquad Two profiles $p$ and $q$ are $h$\textit{-variants}, where $h\in N$,
if $q(i)=p(i)$ for all $i\neq h$. Individual $h$ can \textit{manipulate} the
social choice rule $g:\wp \rightarrow X$ at $p$ via $p^{\prime }$ if $p$ and 
$p^{\prime }$ belong to $\wp $, $p$ and $p^{\prime }$ are $h$-variants, and $%
g(p^{\prime })\succ _{p(h)}g(p)$. And $g$ is \textit{strategy-proof} if no
one can manipulate $g$ at any profile.\medskip

\medskip

\textbf{3. The N}$-$\textbf{Range Theorem: Part 1.}\medskip

As observed earlier, our focus is on the range of $g$ restricted to special
subdomains of $NP$. \ We define $NP^{\ast }$ to be the subdomain of all
profiles $u$ on which individials $n-1$ and $n$ agree: $u(n-1)=u(n)$. \ The
goal is to show \medskip

\textbf{Theorem 3-1}. \ (The N-Range Theorem). \ If. $m=3$ and $g$ on $NP$
is of full range, so is $g|NP^{\ast }$. \medskip

We address the converse: We assume that $g|NP^{\ast }$ is of less than full
range and use that to show $g$ must be of less than full range on $NP$. \ We
will actually do this in two parts by establishing\medskip

\textbf{Theorem 3-2}. \ (The N-Range Theorem, Part 1). \ If $m=3$ and $%
g|NP^{\ast }$ has range of just two alternatives then $g$ must be of less
than full range on $NP$.in this section and then proving \medskip

\textbf{Theorem 4-1}. \ (The N-Range Theorem, Part 2). \ If $m=3$ and $%
g|NP^{\ast }$ has range of just one alternative then $g$ must be of less
than full range on $NP$.in Section 4.\medskip

\qquad The more detailed strategy for proving Theorem 3-2 starts by
determining a very short list $L$ of profiles such that if $x$ is chosen at
some profile in $NP$ it will also have to be chosen at a profile in $L$.
Then for each profile u in $L$, we exploit a decisiveness structure for $%
g|NP^{\ast }$ to show that $x$ is not chosen at $u$ if $x$ is not chosen in $%
NP^{\ast }$. Typically, we will assume $x$ is chosen at $u$ and shoow that
leads to a manipulability at some profile in $NP$, a contradiction of
strategy-proofness.\medskip

\textbf{Section 3-1. The list of profiles}\medskip

\textbf{Lemma 3-3. \ Suppose Range(}$g|NP^{\ast })=\{y,z\}$. \ There exists
a list $L$ of three profiles (along with other profiles that can be
transformed from a member of $L$ by switching $y$ and $z$ everywhere; or
switching preferences of $n-1$ and $n$ or switching orderings within the set
of individuals $i<n-1$) such that, if $x\in $ Range($g$), for strategy-proof 
$g$, then $x=g(u)$ for some $u$ in the list $L$.\medskip

\textbf{Proof of Lemma 3-3:} A profile $u$ in the list $L$ will be described
in terms of where $x$ appears in the orderings for individuals $n-1$ and $n$%
.\medskip

Suppose $x$ is at the bottom of both $u(n-1)$ and $u(n)$. Then it is
possible to switch $y$ and $z$ for one of the two to get a profile $u^{\ast
} $ that is not just in $NP$, but in $NP^{\ast }$, so $g(u^{\ast })\neq x$.
But then the switching individual would manipulate from $u$ to $u^{\ast }$,
a violation of strategy-proofness.\medskip

Similarly, suppose $x$ is at the top of both $u(n-1)$ and $u(n)$. Then it is
possible to switch $y$ and $z$ for one of the two to get a profile $u^{\ast
} $ that is in $NP^{\ast }$, so $g(u^{\ast })\neq x$. But then the switching
individual would manipulate from $u^{\ast }$ to $u$, a violation of
strategy-proofness. There remain four possibilities:\medskip

\qquad\ \ \ \textbf{I.} \qquad $x$ is at the top for one and at the bottom
for the other;\medskip

\qquad\ \ \textbf{II.} \qquad $x$ is at the bottom for one and in the middle
for the other;\medskip

\qquad \textbf{III.} \qquad $x$ is at the top for one and in the middle for
the other;\medskip

\qquad \textbf{IV.} \qquad $x$ is in the middle for both.\medskip

\textbf{Case I.} Suppose $x$ is at the top for $n$ and at the bottom for $%
n-1 $ for profile $u0$ where $g(u0)=x$. Then profile $u1$ is constructed
from $u0 $ by switching $y$ and $z$ for either $n-1$ or $n$ so $y$ and $z$
are ordered by $n-1$ the opposite of their ordering by n. We have $g(u1)=x$
or $g $ is manipulable. Then raise $x$ to the top for each of the $i<n-1$ in
turn. The resulting profile $u2$ will have $g(u_{2})=x$ or else $g$ is
manipulable. Finally change the ordering of $y$ and $z$ for all $i<n-1$ to
agree with the ordering of $y$ and $z$ by $n$. Such a profile looks
like\medskip

\begin{tabular}{|c|c|c|c|c|c|}
\hline
1 & 2 & 3 & $\cdots $ & $n-1$ & n \\ \hline
$%
\begin{array}{c}
x \\ 
y \\ 
z%
\end{array}%
$ & $%
\begin{array}{c}
x \\ 
y \\ 
z%
\end{array}%
$ & $%
\begin{array}{c}
x \\ 
y \\ 
z%
\end{array}%
$ & $\cdots $ & $%
\begin{array}{c}
z \\ 
y \\ 
x%
\end{array}%
$ & $%
\begin{array}{c}
x \\ 
y \\ 
z%
\end{array}%
$ \\ \hline
\end{tabular}%
\medskip \newline
This profile, $L1$, will be in the list $L$ (along with other profiles that
can be transformed from this by switching $y$ and $z$ everywhere; or
switching preferences of $n-1$ and $n$, but all of these will be equivalent
in that if a strategy-proof rule could have $x$ chosen at one but not on $%
NP^{\ast }$, then a rule could be designed that had $x$ chosen at any other
one but not on $NP^{\ast }$.\medskip

\textbf{Case II.} Suppose $x$ is at the bottom for one and in the middle for
the other, say\medskip

\begin{tabular}{|c|c|c|c|c|c|}
\hline
1 & 2 & 3 & $\cdots $ & $n-1$ & n \\ \hline
&  &  & $\cdots $ & $%
\begin{array}{c}
y \\ 
x \\ 
z%
\end{array}%
$ & $%
\begin{array}{c}
z \\ 
y \\ 
x%
\end{array}%
$ \\ \hline
\end{tabular}%
\medskip \newline
Here, if necessary, $y$ and $z$ will be switched for $n$ to be oppositely
ordered from how they are ordered by $n-1$. If that's not possible because
every $i<n-1$ has $y$ above $z$, select one of them, say 1, raise $x$ to the
top of 1's ordering, switch $y$ and $z$ and then continue. But then $x$ can
be raised for $n-1$ and still have $x$ chosen. This puts us back into Case
I, so we do not have to add any profiles from Case II to the list $L$%
.\medskip

\textbf{Case III.} Suppose $x$ is at the top for one and in the middle for
the other, say\medskip

\begin{tabular}{|c|c|c|c|c|c|}
\hline
1 & 2 & 3 & $\cdots $ & $n-1$ & n \\ \hline
&  &  & $\cdots $ & $%
\begin{array}{c}
y \\ 
x \\ 
z%
\end{array}%
$ & $%
\begin{array}{c}
x \\ 
z \\ 
y%
\end{array}%
$ \\ \hline
\end{tabular}%
\medskip \newline
Here, if necessary, $y$ and $z$ will be switched for $n$ to be oppositely
ordered from how they are ordered by $n-1$. If any $i<n-1$ has $y$ preferred
to $x$ we could raise $x$ in $n-1$'s ordering, stay in $NP$, and still have $%
x$ chosen. But this is a case we have already treated. So we may assume
every $i<n-1$ prefers $x$ to $y$. Also to be in $NP$, at least one of them,
say 1, must have $z$ preferred to $x$:\medskip

\begin{tabular}{|c|c|c|c|c|c|}
\hline
1 & 2 & 3 & $\cdots $ & $n-1$ & n \\ \hline
$%
\begin{array}{c}
z \\ 
x \\ 
y%
\end{array}%
$ & $%
\begin{array}{c}
\\ 
x \\ 
y%
\end{array}%
$ & $%
\begin{array}{c}
\\ 
x \\ 
y%
\end{array}%
$ & $\cdots $ & $%
\begin{array}{c}
y \\ 
x \\ 
z%
\end{array}%
$ & $%
\begin{array}{c}
x \\ 
z \\ 
y%
\end{array}%
$ \\ \hline
\end{tabular}%
\medskip \newline
Then for every $i$ such that $2<i<n-1$, we can raise $x$ to the top, and
switch $y$ and $z$ if necessary to agree with $n$, stay in $NP$ and still
have $x$ chosen at the following profile: \medskip

\begin{tabular}{|c|c|c|c|c|c|}
\hline
1 & 2 & 3 & $\cdots $ & $n-1$ & n \\ \hline
$%
\begin{array}{c}
z \\ 
x \\ 
y%
\end{array}%
$ & $%
\begin{array}{c}
x \\ 
z \\ 
y%
\end{array}%
$ & $%
\begin{array}{c}
x \\ 
z \\ 
y%
\end{array}%
$ & $\cdots $ & $%
\begin{array}{c}
y \\ 
x \\ 
z%
\end{array}%
$ & $%
\begin{array}{c}
x \\ 
z \\ 
y%
\end{array}%
$ \\ \hline
\end{tabular}%
\medskip \newline
This profile, $L2$, will also be in list $L$ (along with other profiles that
can be transformed from this by switching $y$ and $z$ everywhere; or
switching preferences of $n-1$ and $n$ or switching orderings within the set
of individuals $i<n-1$).\medskip

\textbf{Case IV}. Suppose $x$ is in the middle for both. Since $x$ is chosen
we are not in $NP^{\ast }$, and the ordering for $n-1$ must be the inverse
of the ordering for n:\medskip

\begin{tabular}{|c|c|c|c|c|c|c|}
\hline
1 & 2 & 3 & $\cdots $ & $n-2$ & $n-1$ & n \\ \hline
&  &  & $\cdots $ &  & $%
\begin{array}{c}
y \\ 
x \\ 
z%
\end{array}%
$ & $%
\begin{array}{c}
z \\ 
x \\ 
y%
\end{array}%
$ \\ \hline
\end{tabular}%
\medskip \newline
If any $i<n-1$ has $y$ preferred to $x$, we could raise $x$ to the top for $%
n-1$, still be in $NP$, and still have $x$ chosen. That would put us in Case
III, and we wouldn't have to add anything to list $L$. Similarly, if any $%
i<n-1$ has $z$ $\succ $ $x$, we could raise $x$ to the top for $n$, still be
in $NP$, and still have $x$ chosen. That would also put us in Case III, and
we wouldn't have to add anything to list $L$. So we may assume that $x$ is
at the top for every $i<n-1$. Then $y$ and $z$ could be switched for each $%
i<n-1$ if necessary to be ordered the same as for n:\medskip

\begin{tabular}{|c|c|c|c|c|c|c|}
\hline
1 & 2 & 3 & $\cdots $ & $n-2$ & $n-1$ & n \\ \hline
$%
\begin{array}{c}
x \\ 
z \\ 
y%
\end{array}%
$ & $%
\begin{array}{c}
x \\ 
z \\ 
y%
\end{array}%
$ & $%
\begin{array}{c}
x \\ 
z \\ 
y%
\end{array}%
$ & $\cdots $ & $%
\begin{array}{c}
x \\ 
z \\ 
y%
\end{array}%
$ & $%
\begin{array}{c}
y \\ 
x \\ 
z%
\end{array}%
$ & $%
\begin{array}{c}
z \\ 
x \\ 
y%
\end{array}%
$ \\ \hline
\end{tabular}%
\medskip \newline
This profile, $L3$, will be the third and last in the list $L$ (along with
other profiles that can be transformed from this by switching $y$ and $z$
everywhere; or switching preferences of $n-1$ and $n$). \ \ \ \ $\square $%
\medskip

\textbf{Section 3-2. Decisiveness structures}\medskip \textbf{Proof of
Theorem 3-2}. \ For each of the three profile types in $L$, we exploit
decisiveness structures to show that $x$ is not chosen at that profile.

\bigskip .

Given strategy-proof rule $g$ on the $NP$ domain for $m$ alternatives and $n$
individuals with Range($g$) $=\{y,z\}$, we define a rule $g^{\ast }$ on the $%
NP$ domain for $m$ alternatives and $n-1$ individuals. \ At profile $%
u=(u_{1},u_{2},...,u_{n-1})$, we set$g^{\ast
}(u)=g(u_{1},u_{2},...,u_{n-1},u_{n})$, where $u_{n}=u_{n-1}$. so $g$ is
operating on a profile in $NP^{\ast }$. \ We have observed earlier (Campbell
and Kelly, 2014b) that $g^{\ast }$ is strategy-proof; it clearly has range $%
\{y,z\}$. \ A result of Barber\'{a} et al (2010) can be modified to show
that for $g^{\ast }$ there is a collection of coalitions decisive for $y$
against $z$ and a related collection of coalitions decisive for $z$ against $%
y$, each collection satisfying a monotonicity condition: supersets of
members are also members.\medskip 

\qquad Correspondingly, then, we know a decisiveness structure for $%
g|NP^{\ast }$. \ Let $C$ be a coalition in $\{1,2,...,n-2\}$.\medskip

\qquad \qquad 1. \ If $C$ is decisive for $g^{\ast }$ for $y$ against $z$ ($%
z $ against $y$), then $C$ is decisive for $g|NP^{\ast }$ for $y$ against $z$
($z$ against $y$).

\qquad \qquad 2. \ If $C$ is minimally decisive for $g^{\ast }$ for $y$
against $z$ ($z$ against $y$), then $C$ is minimally decisive for $%
g|NP^{\ast }$ for $y$ against $z$ ($z$ against $y$).

\qquad \qquad 3. \ If $C\cup \{n-1\}$ is decisive for $g^{\ast }$ for $y$
against $z$ ($z$ against $y$), then $C\cup \{n-1,n\}$ is decisive for $%
g|NP^{\ast }$ for $y$ against $z$ ($z$ against $y$).

\qquad \qquad 4. \ If $C\cup \{n-1\}$ is minimally decisive for $g^{\ast }$
for $y$ against $z$ ($z$ against $y$), there is no proper subset $C^{\ast }$
of $C$ such that $C^{\ast }\cup \{n-1,n\}$ is decisive for $g|NP^{\ast }$
for $y$ against $z$ ($z$ against $y$). With some license, we will say $C\cup
\{n-1,n\}$ is minimally decisive for $g|NP^{\ast }$.\medskip

\qquad There are two possible categories of decisiveness structures for $%
g^{\ast }$. In the first category, some minimal decisive coalition for $y$
against $z$ for $g^{\ast }$ is contained in $\{1,...,n-2\}$. (Or,
alternatively, some minimum decisive coalition for $z$ against $y$ for $g$
is contained in $\{1,...,n-2\}$.) Within this first category, there are two
possibilities to consider:\medskip

\textbf{Case A}. Some subset $S$ of $\{1,...,n-2\}$, containing at least two
elements, is a minimal decisive coalition for $y$ against $z$
(alternatively, some subset of $\{1,...,n-2\}$, containing at least two
elements, is a minimal decisive coalition for $z$ against y). Since $N$
contains more than two individuals, we then automatically have useful
coalitions decisive (though possibly not minimally) for $z$ against $y$,
namely $X\backslash S$ together with a non$-$empty proper subset of $S$%
.\medskip

\textbf{Case B}. Some singleton subset $S$ of $\{1,...,n-2\}$ is a minimal
decisive coalition for $y$ against $z$. When $S$ is a singleton, we
sometimes must think carefully about coalitions decisive for $z$ against $y$
(of course, since $Range(g^{\ast })=\{y,z\}$, some such coalitions
exist).\medskip

\qquad In the second category, every minimal decisive coalition for $y$
against $z$ for $g^{\ast }$ includes $n-1$. The possibilities are that
either $\{n-1\}$ itself is a minimal decisive coalition for $y$ against $z$
or that every minimal decisive coalition for $y$ against $z$ includes both $%
n-1$ and some individual in $\{1,...,n-2\}$. But in the latter case, some
subset of $\{1,...,n-2\}$ is a minimal decisive coalition for $z$ against $y$
and we are back in Case A or Case B. So we only need to treat the following
possibility:\medskip

\textbf{Case C}. $\{n-1\}$ is a minimal decisive coalition for $y$ against $%
z $ for $g^{\ast }$ and $\{n-1\}$ is also a minimal decisive coalition for $%
z $ against $y$ for $g^{\ast }$. So $\{n-1,n\}$ is a decisive coalition for $%
y$ against $z$ and for $z$ against $y$ for $g|NP^{\ast }$.\medskip

Accordingly, with three profiles in $L$ to treat and, for each of those
profiles, three kinds of decisiveness structures to consider, we complete
the proof by carrying nine tasks.\medskip

\textbf{Section 3-3.} $L1$.\medskip

We want to show that for any collection of decisive coalitions, $x\notin
Range(g|NP^{\ast })$ and strategy-proofness of $g$ will exclude the
possibility that $x$ is chosen at profile $L1$: \medskip

\begin{tabular}{c|c|c|c|c|c|c|c|}
\cline{2-8}
$L1$ & $1$ & $2$ & $3$ & $\cdots $ & $n-2$ & $n-1$ & $n$ \\ \cline{2-8}
& $%
\begin{array}{c}
x \\ 
y \\ 
z%
\end{array}%
$ & $%
\begin{array}{c}
x \\ 
y \\ 
z%
\end{array}%
$ & $%
\begin{array}{c}
x \\ 
y \\ 
z%
\end{array}%
$ & $\cdots $ & $%
\begin{array}{c}
x \\ 
y \\ 
z%
\end{array}%
$ & $%
\begin{array}{c}
z \\ 
y \\ 
x%
\end{array}%
$ & $%
\begin{array}{c}
x \\ 
y \\ 
z%
\end{array}%
$ \\ \cline{2-8}
\end{tabular}%
\medskip \newline
(All profiles displayed in this paper are elements of $NP$; this is easily
checked and will not be further remarked upon.)\medskip

\qquad \textbf{Case }$L1$.A. Some subset $S$ of $\{1,...,n-2\}$, containing
at least two elements, is a minimal decisive coalition for $y$ against $z$
(alternatively, some subset of $\{1,...,n-2\}$, containing at least two
elements, is a minimal decisive coalition for $z$ against y). We treat here
the case where $\{1,2,..,k\}$ is a \textit{minimal} decisive set for $y$
against $z$ (the case for a subset of $\{1,...,n-2\}$ being decisive for $z$
against $y$ can be dealt with similarly). Note that $k=n-2$ is
allowed.\medskip

We want to show that for minimal decisive coalition $\{1,2,..,k\}$, a
violation of strategy-proofness of $g$ will follow from an assumption that $%
x $ is chosen at profile $L1$:\medskip

\begin{tabular}{c|c|c|c|c|c|c|c|c|c|c|c|}
\cline{2-12}
$L1$ & $1$ & $2$ & $3$ & $\cdots $ & $k-1$ & $k$ & $k+1$ & $\cdots $ & $n-2$
& $n-1$ & $n$ \\ \cline{2-12}
& $%
\begin{array}{c}
x \\ 
y \\ 
z%
\end{array}%
$ & $%
\begin{array}{c}
x \\ 
y \\ 
z%
\end{array}%
$ & $%
\begin{array}{c}
x \\ 
y \\ 
z%
\end{array}%
$ & $\cdots $ & $%
\begin{array}{c}
x \\ 
y \\ 
z%
\end{array}%
$ & $%
\begin{array}{c}
x \\ 
y \\ 
z%
\end{array}%
$ & $%
\begin{array}{c}
x \\ 
y \\ 
z%
\end{array}%
$ & $\cdots $ & $%
\begin{array}{c}
x \\ 
y \\ 
z%
\end{array}%
$ & $%
\begin{array}{c}
z \\ 
y \\ 
x%
\end{array}%
$ & $%
\begin{array}{c}
x \\ 
y \\ 
z%
\end{array}%
$ \\ \cline{2-12}
\end{tabular}%
\medskip \newline
Alternative $x$ is chosen at $L1$ if and only if $x$ is chosen at $L1^{\ast
} $:\medskip

\begin{tabular}{c|c|c|c|c|c|c|c|c|c|c|c|}
\cline{2-12}
$L1^{\ast }$ & $1$ & $2$ & $3$ & $\cdots $ & $k-1$ & $k$ & $k+1$ & $\cdots $
& $n-2$ & $n-1$ & $n$ \\ \cline{2-12}
& $%
\begin{array}{c}
x \\ 
y \\ 
z%
\end{array}%
$ & $%
\begin{array}{c}
x \\ 
y \\ 
z%
\end{array}%
$ & $%
\begin{array}{c}
x \\ 
y \\ 
z%
\end{array}%
$ & $\cdots $ & $%
\begin{array}{c}
x \\ 
y \\ 
z%
\end{array}%
$ & $%
\begin{array}{c}
x \\ 
y \\ 
z%
\end{array}%
$ & $%
\begin{array}{c}
x \\ 
z \\ 
y%
\end{array}%
$ & $\cdots $ & $%
\begin{array}{c}
x \\ 
z \\ 
y%
\end{array}%
$ & $%
\begin{array}{c}
z \\ 
y \\ 
x%
\end{array}%
$ & $%
\begin{array}{c}
x \\ 
y \\ 
z%
\end{array}%
$ \\ \cline{2-12}
\end{tabular}%
\medskip \newline
with $y$ and $z$ reversed below $x$ for $\{k+1,...,n-2\}$.\medskip

At $n$-variant $u1$ in $NP^{\ast }$:\medskip

\begin{tabular}{c|c|c|c|c|c|c|c|c|c|c|c|}
\cline{2-12}
$u1$ & $1$ & $2$ & $3$ & $\cdots $ & $k-1$ & $k$ & $k+1$ & $\cdots $ & $n-2$
& $n-1$ & $n$ \\ \cline{2-12}
& $%
\begin{array}{c}
x \\ 
y \\ 
z%
\end{array}%
$ & $%
\begin{array}{c}
x \\ 
y \\ 
z%
\end{array}%
$ & $%
\begin{array}{c}
x \\ 
y \\ 
z%
\end{array}%
$ & $\cdots $ & $%
\begin{array}{c}
x \\ 
y \\ 
z%
\end{array}%
$ & $%
\begin{array}{c}
x \\ 
y \\ 
z%
\end{array}%
$ & $%
\begin{array}{c}
x \\ 
z \\ 
y%
\end{array}%
$ & $\cdots $ & $%
\begin{array}{c}
x \\ 
z \\ 
y%
\end{array}%
$ & $%
\begin{array}{c}
z \\ 
y \\ 
x%
\end{array}%
$ & $%
\begin{array}{c}
z \\ 
y \\ 
x%
\end{array}%
$ \\ \cline{2-12}
\end{tabular}%
\medskip \newline
$g(u1)=y$ by decisiveness of $\{1,...,k\}$ for $y$ against $z$ and then at
2-variant profile $u2$, also in $NP^{\ast }$:\medskip

\begin{tabular}{c|c|c|c|c|c|c|c|c|c|c|c|}
\cline{2-12}
$u2$ & $1$ & $2$ & $3$ & $\cdots $ & $k-1$ & $k$ & $k+1$ & $\cdots $ & $n-2$
& $n-1$ & $n$ \\ \cline{2-12}
& $%
\begin{array}{c}
x \\ 
y \\ 
z%
\end{array}%
$ & $%
\begin{array}{c}
x \\ 
z \\ 
y%
\end{array}%
$ & $%
\begin{array}{c}
x \\ 
y \\ 
z%
\end{array}%
$ & $\cdots $ & $%
\begin{array}{c}
x \\ 
y \\ 
z%
\end{array}%
$ & $%
\begin{array}{c}
x \\ 
y \\ 
z%
\end{array}%
$ & $%
\begin{array}{c}
x \\ 
z \\ 
y%
\end{array}%
$ & $\cdots $ & $%
\begin{array}{c}
x \\ 
z \\ 
y%
\end{array}%
$ & $%
\begin{array}{c}
z \\ 
y \\ 
x%
\end{array}%
$ & $%
\begin{array}{c}
z \\ 
y \\ 
x%
\end{array}%
$ \\ \cline{2-12}
\end{tabular}%
\medskip \newline
$g(u2)=z$ since $\{1,...,k\}$ is a minimal decisive set for $y$ against $z$
on $NP^{\ast }$. Next, at $n$-variant profile $u3$:\medskip

\begin{tabular}{c|c|c|c|c|c|c|c|c|c|c|c|}
\cline{2-12}
$u3$ & $1$ & $2$ & $3$ & $\cdots $ & $k-1$ & $k$ & $k+1$ & $\cdots $ & $n-2$
& $n-1$ & $n$ \\ \cline{2-12}
& $%
\begin{array}{c}
x \\ 
y \\ 
z%
\end{array}%
$ & $%
\begin{array}{c}
x \\ 
z \\ 
y%
\end{array}%
$ & $%
\begin{array}{c}
x \\ 
y \\ 
z%
\end{array}%
$ & $\cdots $ & $%
\begin{array}{c}
x \\ 
y \\ 
z%
\end{array}%
$ & $%
\begin{array}{c}
x \\ 
y \\ 
z%
\end{array}%
$ & $%
\begin{array}{c}
x \\ 
z \\ 
y%
\end{array}%
$ & $\cdots $ & $%
\begin{array}{c}
x \\ 
z \\ 
y%
\end{array}%
$ & $%
\begin{array}{c}
z \\ 
y \\ 
x%
\end{array}%
$ & $%
\begin{array}{c}
z \\ 
x \\ 
y%
\end{array}%
$ \\ \cline{2-12}
\end{tabular}%
\medskip \newline
$g(u3$) = $z$ or $n$ manipulates from $u3$ to $u2$. Next consider 2-variant
profile $u4$:\medskip

\begin{tabular}{c|c|c|c|c|c|c|c|c|c|c|c|}
\cline{2-12}
$u4$ & $1$ & $2$ & $3$ & $\cdots $ & $k-1$ & $k$ & $k+1$ & $\cdots $ & $n-2$
& $n-1$ & $n$ \\ \cline{2-12}
& $%
\begin{array}{c}
x \\ 
y \\ 
z%
\end{array}%
$ & $%
\begin{array}{c}
x \\ 
y \\ 
z%
\end{array}%
$ & $%
\begin{array}{c}
x \\ 
y \\ 
z%
\end{array}%
$ & $\cdots $ & $%
\begin{array}{c}
x \\ 
y \\ 
z%
\end{array}%
$ & $%
\begin{array}{c}
x \\ 
y \\ 
z%
\end{array}%
$ & $%
\begin{array}{c}
x \\ 
z \\ 
y%
\end{array}%
$ & $\cdots $ & $%
\begin{array}{c}
x \\ 
z \\ 
y%
\end{array}%
$ & $%
\begin{array}{c}
z \\ 
y \\ 
x%
\end{array}%
$ & $%
\begin{array}{c}
z \\ 
x \\ 
y%
\end{array}%
$ \\ \cline{2-12}
\end{tabular}%
\medskip \newline
$g(u4)\neq z$ or $n$ manipulates from $u1$ to $u4$. And $g(u4)\neq x$ or 2
manipulates from $u3$ to $u4$. So $g(u4)=y$. But then, if $g(L1^{\ast })=x$,
voter $n$ would manipulate from $u4$ to $L1^{\ast }$:\medskip

\begin{tabular}{c|c|c|c|c|c|c|c|c|c|c|c|}
\cline{2-12}
$L1^{\ast }$ & $1$ & $2$ & $3$ & $\cdots $ & $k-1$ & $k$ & $k+1$ & $\cdots $
& $n-2$ & $n-1$ & $n$ \\ \cline{2-12}
& $%
\begin{array}{c}
x \\ 
y \\ 
z%
\end{array}%
$ & $%
\begin{array}{c}
x \\ 
y \\ 
z%
\end{array}%
$ & $%
\begin{array}{c}
x \\ 
y \\ 
z%
\end{array}%
$ & $\cdots $ & $%
\begin{array}{c}
x \\ 
y \\ 
z%
\end{array}%
$ & $%
\begin{array}{c}
x \\ 
y \\ 
z%
\end{array}%
$ & $%
\begin{array}{c}
x \\ 
z \\ 
y%
\end{array}%
$ & $\cdots $ & $%
\begin{array}{c}
x \\ 
z \\ 
y%
\end{array}%
$ & $%
\begin{array}{c}
z \\ 
y \\ 
x%
\end{array}%
$ & $%
\begin{array}{c}
x \\ 
y \\ 
z%
\end{array}%
$ \\ \cline{2-12}
\end{tabular}%
\medskip \newline
\qquad \textbf{Subcase }$L1$.B. Without loss of generality, $\{1\}$ is a
(minimal) decisive coalition for $y$ against $z$ for $g|NP^{\ast }$. There
are many possible coalitions $C$ decisive for $z$ against $y$. These must
include 1, but won't have to be minimal, so we choose them as large as
possible (though they still have to exclude at least one individual) since a
small coalitions being decisive implies supersets also decisive. We
distinguish between two possibilities, regarding which kind of individual is
excluded:\medskip

\qquad 1. The individual excluded for is in $\{1,...,n-2\}$, say $n-2$, so $%
\{1,...,n-3,n-1,n\}$ is decisive for $z$ against $y$ for $g|NP^{\ast }$%
.\medskip

\qquad 2. The individual excluded for is $n-1$, so $\{1,...,n-2\}$ is
decisive for $z$ against $y$ for $g|NP^{\ast }$. If the minimal decisive
coalition for $z$ against $y$ in $\{1,...,n-2\}$ contains an individual in
addition to 1, we are in the situation already covered in Case A. So we may
assume that $\{1\}$ is decisive for $z$ against $y$ as well as for $y$
against $z$ in $NP^{\ast }$. But then also $\{1,...,n-3,n-1,n\}$ is decisive
for $z$ against $y$ for $g|NP^{\ast }$, and we are back to the first
possibility.\medskip

\qquad So we only need to treat the case where coalition $\{1\}$ is decisive
for $y$ against $z$ on $NP^{\ast }$ and $\{1,...,n-3,n-1,n\}$ is decisive
for $z$ against $y$ on $NP^{\ast }$. At $L1$, $x$ is chosen by
assumption:\medskip

\begin{tabular}{c|c|c|c|c|c|c|c|c|}
\cline{2-9}
$L1$ & $1$ & $2$ & $3$ & $\cdots $ & $n-3$ & $n-2$ & $n-1$ & $n$ \\ 
\cline{2-9}
& $%
\begin{array}{c}
x \\ 
y \\ 
z%
\end{array}%
$ & $%
\begin{array}{c}
x \\ 
y \\ 
z%
\end{array}%
$ & $%
\begin{array}{c}
x \\ 
y \\ 
z%
\end{array}%
$ & $\cdots $ & $%
\begin{array}{c}
x \\ 
y \\ 
z%
\end{array}%
$ & $%
\begin{array}{c}
x \\ 
y \\ 
z%
\end{array}%
$ & $%
\begin{array}{c}
z \\ 
y \\ 
x%
\end{array}%
$ & $%
\begin{array}{c}
x \\ 
y \\ 
z%
\end{array}%
$ \\ \cline{2-9}
\end{tabular}%
\medskip \newline
Then $x$ must also be chosen at $L1^{\ast }$:\medskip

\begin{tabular}{c|c|c|c|c|c|c|c|c|}
\cline{2-9}
$L1^{\ast }$ & $1$ & $2$ & $3$ & $\cdots $ & $n-3$ & $n-2$ & $n-1$ & $n$ \\ 
\cline{2-9}
& $%
\begin{array}{c}
x \\ 
y \\ 
z%
\end{array}%
$ & $%
\begin{array}{c}
x \\ 
z \\ 
y%
\end{array}%
$ & $%
\begin{array}{c}
x \\ 
z \\ 
y%
\end{array}%
$ & $\cdots $ & $%
\begin{array}{c}
x \\ 
z \\ 
y%
\end{array}%
$ & $%
\begin{array}{c}
x \\ 
y \\ 
z%
\end{array}%
$ & $%
\begin{array}{c}
z \\ 
y \\ 
x%
\end{array}%
$ & $%
\begin{array}{c}
x \\ 
y \\ 
z%
\end{array}%
$ \\ \cline{2-9}
\end{tabular}%
\medskip \newline
obtained by switching $y$ and $z$ for individuals $2,...,n-3$.\medskip

At $n$-variant $u1$:\medskip

\begin{tabular}{c|c|c|c|c|c|c|c|c|}
\cline{2-9}
$u1$ & $1$ & $2$ & $3$ & $\cdots $ & $n-3$ & $n-2$ & $n-1$ & $n$ \\ 
\cline{2-9}
& $%
\begin{array}{c}
x \\ 
y \\ 
z%
\end{array}%
$ & $%
\begin{array}{c}
x \\ 
z \\ 
y%
\end{array}%
$ & $%
\begin{array}{c}
x \\ 
z \\ 
y%
\end{array}%
$ & $\cdots $ & $%
\begin{array}{c}
x \\ 
z \\ 
y%
\end{array}%
$ & $%
\begin{array}{c}
x \\ 
y \\ 
z%
\end{array}%
$ & $%
\begin{array}{c}
z \\ 
y \\ 
x%
\end{array}%
$ & $%
\begin{array}{c}
z \\ 
x \\ 
y%
\end{array}%
$ \\ \cline{2-9}
\end{tabular}%
\medskip \newline
$g(u1)\neq y$ or $n$ will manipulate to $L1^{\ast }$. Next consider $u2$,
another $n$-variant of $u1$, but one that is in $NP^{\ast }$:\medskip

\begin{tabular}{c|c|c|c|c|c|c|c|c|}
\cline{2-9}
$u2$ & $1$ & $2$ & $3$ & $\cdots $ & $n-3$ & $n-2$ & $n-1$ & $n$ \\ 
\cline{2-9}
& $%
\begin{array}{c}
x \\ 
y \\ 
z%
\end{array}%
$ & $%
\begin{array}{c}
x \\ 
z \\ 
y%
\end{array}%
$ & $%
\begin{array}{c}
x \\ 
z \\ 
y%
\end{array}%
$ & $\cdots $ & $%
\begin{array}{c}
x \\ 
z \\ 
y%
\end{array}%
$ & $%
\begin{array}{c}
x \\ 
y \\ 
z%
\end{array}%
$ & $%
\begin{array}{c}
z \\ 
y \\ 
x%
\end{array}%
$ & $%
\begin{array}{c}
z \\ 
y \\ 
x%
\end{array}%
$ \\ \cline{2-9}
\end{tabular}%
\medskip \newline
$g(u2)=y$, since $\{1\}$ is decisive for $y$ against $z$ on $NP^{\ast }$.
Then $g(u1)\neq z$ or $n$ will manipulate from $u2$ to $u1$. Combining, $%
g(u1)=x$.\medskip

\qquad Now look at profile $u3$, also in $NP^{\ast }$:\medskip

\begin{tabular}{c|c|c|c|c|c|c|c|c|}
\cline{2-9}
$u3$ & $1$ & $2$ & $3$ & $\cdots $ & $n-3$ & $n-2$ & $n-1$ & $n$ \\ 
\cline{2-9}
& $%
\begin{array}{c}
x \\ 
z \\ 
y%
\end{array}%
$ & $%
\begin{array}{c}
x \\ 
z \\ 
y%
\end{array}%
$ & $%
\begin{array}{c}
x \\ 
z \\ 
y%
\end{array}%
$ & $\cdots $ & $%
\begin{array}{c}
x \\ 
z \\ 
y%
\end{array}%
$ & $%
\begin{array}{c}
x \\ 
y \\ 
z%
\end{array}%
$ & $%
\begin{array}{c}
z \\ 
y \\ 
x%
\end{array}%
$ & $%
\begin{array}{c}
z \\ 
y \\ 
x%
\end{array}%
$ \\ \cline{2-9}
\end{tabular}%
\medskip \newline
$g(u3)=z$ by the decisiveness of $\{1,2,...,n-3,n-1,n\}$ for $z$ against $y$
on $NP^{\ast }$. Then at $n$-variant profile $u4$:\medskip

\begin{tabular}{c|c|c|c|c|c|c|c|c|}
\cline{2-9}
$u4$ & $1$ & $2$ & $3$ & $\cdots $ & $n-3$ & $n-2$ & $n-1$ & $n$ \\ 
\cline{2-9}
& $%
\begin{array}{c}
x \\ 
z \\ 
y%
\end{array}%
$ & $%
\begin{array}{c}
x \\ 
z \\ 
y%
\end{array}%
$ & $%
\begin{array}{c}
x \\ 
z \\ 
y%
\end{array}%
$ & $\cdots $ & $%
\begin{array}{c}
x \\ 
z \\ 
y%
\end{array}%
$ & $%
\begin{array}{c}
x \\ 
y \\ 
z%
\end{array}%
$ & $%
\begin{array}{c}
z \\ 
y \\ 
x%
\end{array}%
$ & $%
\begin{array}{c}
z \\ 
x \\ 
y%
\end{array}%
$ \\ \cline{2-9}
\end{tabular}%
\medskip \newline
$g(u4)=z$ or $n$ will manipulate from $u4$ to $u3$. But $u4$ is a 1-variant
of $u1$ and 1 will manipulate from $u4$ to $u1$.\medskip

\textbf{Case }$L1$.C: In this case $\{n-1\}$ is decisive for $y$ against $z$
and also for $z$ against $y$ for rule $g^{\ast }$, i.e., $\{n-1,n\}$ is
decisive both ways for rule $g|NP^{\ast }$. We first trace out two results
of decisiveness and strategy-proofness. \medskip

\textbf{Result \#1}\medskip

We assume $x$ is chosen at $L1$:\medskip

\begin{tabular}{c|c|c|c|c|c|c|c|}
\cline{2-8}
$L1$ & $1$ & $2$ & $3$ & $\cdots $ & $n-2$ & $n-1$ & $n$ \\ \cline{2-8}
& $%
\begin{array}{c}
x \\ 
y \\ 
z%
\end{array}%
$ & $%
\begin{array}{c}
x \\ 
y \\ 
z%
\end{array}%
$ & $%
\begin{array}{c}
x \\ 
y \\ 
z%
\end{array}%
$ & $\cdots $ & $%
\begin{array}{c}
x \\ 
y \\ 
z%
\end{array}%
$ & $%
\begin{array}{c}
z \\ 
y \\ 
x%
\end{array}%
$ & $%
\begin{array}{c}
x \\ 
y \\ 
z%
\end{array}%
$ \\ \cline{2-8}
\end{tabular}%
\medskip \newline
and seek a contradiction. We must have $g(u1)=x$ at 1-variant $u1$:\medskip

\begin{tabular}{c|c|c|c|c|c|c|c|}
\cline{2-8}
$u1$ & $1$ & $2$ & $3$ & $\cdots $ & $n-2$ & $n-1$ & $n$ \\ \cline{2-8}
& $%
\begin{array}{c}
x \\ 
z \\ 
y%
\end{array}%
$ & $%
\begin{array}{c}
x \\ 
y \\ 
z%
\end{array}%
$ & $%
\begin{array}{c}
x \\ 
y \\ 
z%
\end{array}%
$ & $\cdots $ & $%
\begin{array}{c}
x \\ 
y \\ 
z%
\end{array}%
$ & $%
\begin{array}{c}
z \\ 
y \\ 
x%
\end{array}%
$ & $%
\begin{array}{c}
x \\ 
y \\ 
z%
\end{array}%
$ \\ \cline{2-8}
\end{tabular}%
\medskip \newline
or 1 manipulates from $u1$ to $L1$. Then $g(u2)=x$ at $(n-1)$-variant $u2$%
:\medskip

\begin{tabular}{c|c|c|c|c|c|c|c|}
\cline{2-8}
$u2$ & $1$ & $2$ & $3$ & $\cdots $ & $n-2$ & $n-1$ & $n$ \\ \cline{2-8}
& $%
\begin{array}{c}
x \\ 
z \\ 
y%
\end{array}%
$ & $%
\begin{array}{c}
x \\ 
y \\ 
z%
\end{array}%
$ & $%
\begin{array}{c}
x \\ 
y \\ 
z%
\end{array}%
$ & $\cdots $ & $%
\begin{array}{c}
x \\ 
y \\ 
z%
\end{array}%
$ & $%
\begin{array}{c}
y \\ 
z \\ 
x%
\end{array}%
$ & $%
\begin{array}{c}
x \\ 
y \\ 
z%
\end{array}%
$ \\ \cline{2-8}
\end{tabular}%
\medskip \newline
or $n-1$ manipulates from $u1$ to $u2$. Now $g(u3)\neq y$ at 1-variant $u3$%
:\medskip

\begin{tabular}{c|c|c|c|c|c|c|c|}
\cline{2-8}
$u3$ & $1$ & $2$ & $3$ & $\cdots $ & $n-2$ & $n-1$ & $n$ \\ \cline{2-8}
& $%
\begin{array}{c}
z \\ 
x \\ 
y%
\end{array}%
$ & $%
\begin{array}{c}
x \\ 
y \\ 
z%
\end{array}%
$ & $%
\begin{array}{c}
x \\ 
y \\ 
z%
\end{array}%
$ & $\cdots $ & $%
\begin{array}{c}
x \\ 
y \\ 
z%
\end{array}%
$ & $%
\begin{array}{c}
y \\ 
z \\ 
x%
\end{array}%
$ & $%
\begin{array}{c}
x \\ 
y \\ 
z%
\end{array}%
$ \\ \cline{2-8}
\end{tabular}%
\medskip \newline
or 1 manipulates from $u3$ to $u2$. We next show $g(u3)\neq z$.

\qquad At $u4$ in $NP^{\ast }$:\medskip

\begin{tabular}{c|c|c|c|c|c|c|c|}
\cline{2-8}
$u4$ & $1$ & $2$ & $3$ & $\cdots $ & $n-2$ & $n-1$ & $n$ \\ \cline{2-8}
& $%
\begin{array}{c}
z \\ 
y \\ 
x%
\end{array}%
$ & $%
\begin{array}{c}
x \\ 
y \\ 
z%
\end{array}%
$ & $%
\begin{array}{c}
x \\ 
y \\ 
z%
\end{array}%
$ & $\cdots $ & $%
\begin{array}{c}
x \\ 
y \\ 
z%
\end{array}%
$ & $%
\begin{array}{c}
x \\ 
y \\ 
z%
\end{array}%
$ & $%
\begin{array}{c}
x \\ 
y \\ 
z%
\end{array}%
$ \\ \cline{2-8}
\end{tabular}%
\medskip \newline
$g(u4)=y$ by the decisiveness of $\{n-1,n\}$ on $NP^{\ast }$. That implies $%
g(u5)=y$ at $(n-1)$-variant profile $u5$:\medskip

\begin{tabular}{c|c|c|c|c|c|c|c|}
\cline{2-8}
$u5$ & $1$ & $2$ & $3$ & $\cdots $ & $n-2$ & $n-1$ & $n$ \\ \cline{2-8}
& $%
\begin{array}{c}
z \\ 
y \\ 
x%
\end{array}%
$ & $%
\begin{array}{c}
x \\ 
y \\ 
z%
\end{array}%
$ & $%
\begin{array}{c}
x \\ 
y \\ 
z%
\end{array}%
$ & $\cdots $ & $%
\begin{array}{c}
x \\ 
y \\ 
z%
\end{array}%
$ & $%
\begin{array}{c}
y \\ 
z \\ 
x%
\end{array}%
$ & $%
\begin{array}{c}
x \\ 
y \\ 
z%
\end{array}%
$ \\ \cline{2-8}
\end{tabular}%
\medskip \newline
or $n-1$ manipulates from $u5$ to $u4$. But $u5$ is a 1-variant of $u3$%
:\medskip

\begin{tabular}{c|c|c|c|c|c|c|c|}
\cline{2-8}
$u3$ & $1$ & $2$ & $3$ & $\cdots $ & $n-2$ & $n-1$ & $n$ \\ \cline{2-8}
& $%
\begin{array}{c}
z \\ 
x \\ 
y%
\end{array}%
$ & $%
\begin{array}{c}
x \\ 
y \\ 
z%
\end{array}%
$ & $%
\begin{array}{c}
x \\ 
y \\ 
z%
\end{array}%
$ & $\cdots $ & $%
\begin{array}{c}
x \\ 
y \\ 
z%
\end{array}%
$ & $%
\begin{array}{c}
y \\ 
z \\ 
x%
\end{array}%
$ & $%
\begin{array}{c}
x \\ 
y \\ 
z%
\end{array}%
$ \\ \cline{2-8}
\end{tabular}%
\medskip \newline
So, if $g(u3)=z$, 1 would manipulate from $u5$ to $u3$. Therefore $g(u3)\neq
z$. Combining, $g(u3)=x$.\medskip

\textbf{Result \#2}\medskip

Earlier we saw that if $x$ is selected at $L1$ then $g(u2)=x$ at $u2$%
:\medskip

\begin{tabular}{c|c|c|c|c|c|c|c|}
\cline{2-8}
$u2$ & $1$ & $2$ & $3$ & $\cdots $ & $n-2$ & $n-1$ & $n$ \\ \cline{2-8}
& $%
\begin{array}{c}
x \\ 
z \\ 
y%
\end{array}%
$ & $%
\begin{array}{c}
x \\ 
y \\ 
z%
\end{array}%
$ & $%
\begin{array}{c}
x \\ 
y \\ 
z%
\end{array}%
$ & $\cdots $ & $%
\begin{array}{c}
x \\ 
y \\ 
z%
\end{array}%
$ & $%
\begin{array}{c}
y \\ 
z \\ 
x%
\end{array}%
$ & $%
\begin{array}{c}
x \\ 
y \\ 
z%
\end{array}%
$ \\ \cline{2-8}
\end{tabular}%
\medskip \newline
Then $g(u6)=x$ at $n$-variant $u6$:\medskip

\begin{tabular}{c|c|c|c|c|c|c|c|}
\cline{2-8}
$u6$ & $1$ & $2$ & $3$ & $\cdots $ & $n-2$ & $n-1$ & $n$ \\ \cline{2-8}
& $%
\begin{array}{c}
x \\ 
z \\ 
y%
\end{array}%
$ & $%
\begin{array}{c}
x \\ 
y \\ 
z%
\end{array}%
$ & $%
\begin{array}{c}
x \\ 
y \\ 
z%
\end{array}%
$ & $\cdots $ & $%
\begin{array}{c}
x \\ 
y \\ 
z%
\end{array}%
$ & $%
\begin{array}{c}
y \\ 
z \\ 
x%
\end{array}%
$ & $%
\begin{array}{c}
x \\ 
z \\ 
y%
\end{array}%
$ \\ \cline{2-8}
\end{tabular}%
\medskip \newline
or $n$ manipulates from $u6$ to $u2$. Therefore $g(u7)\neq y$ at 1-variant $%
u7$:\medskip

\begin{tabular}{c|c|c|c|c|c|c|c|}
\cline{2-8}
$u7$ & $1$ & $2$ & $3$ & $\cdots $ & $n-2$ & $n-1$ & $n$ \\ \cline{2-8}
& $%
\begin{array}{c}
z \\ 
x \\ 
y%
\end{array}%
$ & $%
\begin{array}{c}
x \\ 
y \\ 
z%
\end{array}%
$ & $%
\begin{array}{c}
x \\ 
y \\ 
z%
\end{array}%
$ & $\cdots $ & $%
\begin{array}{c}
x \\ 
y \\ 
z%
\end{array}%
$ & $%
\begin{array}{c}
y \\ 
z \\ 
x%
\end{array}%
$ & $%
\begin{array}{c}
x \\ 
z \\ 
y%
\end{array}%
$ \\ \cline{2-8}
\end{tabular}%
\medskip \newline
or 1 manipulates from $u7$ to $u6$. We next show $g(u7)\neq x$.\medskip

\qquad At $u8$ in $NP^{\ast }$\medskip

\begin{tabular}{c|c|c|c|c|c|c|c|}
\cline{2-8}
$u8$ & $1$ & $2$ & $3$ & $\cdots $ & $n-2$ & $n-1$ & $n$ \\ \cline{2-8}
& $%
\begin{array}{c}
z \\ 
x \\ 
y%
\end{array}%
$ & $%
\begin{array}{c}
y \\ 
x \\ 
z%
\end{array}%
$ & $%
\begin{array}{c}
x \\ 
y \\ 
z%
\end{array}%
$ & $\cdots $ & $%
\begin{array}{c}
x \\ 
y \\ 
z%
\end{array}%
$ & $%
\begin{array}{c}
x \\ 
z \\ 
y%
\end{array}%
$ & $%
\begin{array}{c}
x \\ 
z \\ 
y%
\end{array}%
$ \\ \cline{2-8}
\end{tabular}%
\medskip \newline
$g(u8)=z$ by the decisiveness of $\{n-1,n\}$ on $NP^{\ast }$. Then $g(u9)=z$
at $(n-1)$-variant $u9$:\medskip

\begin{tabular}{c|c|c|c|c|c|c|c|}
\cline{2-8}
$u9$ & $1$ & $2$ & $3$ & $\cdots $ & $n-2$ & $n-1$ & $n$ \\ \cline{2-8}
& $%
\begin{array}{c}
z \\ 
x \\ 
y%
\end{array}%
$ & $%
\begin{array}{c}
y \\ 
x \\ 
z%
\end{array}%
$ & $%
\begin{array}{c}
x \\ 
y \\ 
z%
\end{array}%
$ & $\cdots $ & $%
\begin{array}{c}
x \\ 
y \\ 
z%
\end{array}%
$ & $%
\begin{array}{c}
z \\ 
y \\ 
x%
\end{array}%
$ & $%
\begin{array}{c}
x \\ 
z \\ 
y%
\end{array}%
$ \\ \cline{2-8}
\end{tabular}%
\medskip \newline
or $n-1$ manipulates from $u9$ to $u8$. Then $g(u10)=z$ at 2-variant $u10$%
:\medskip

\begin{tabular}{c|c|c|c|c|c|c|c|}
\cline{2-8}
$u10$ & $1$ & $2$ & $3$ & $\cdots $ & $n-2$ & $n-1$ & $n$ \\ \cline{2-8}
& $%
\begin{array}{c}
z \\ 
x \\ 
y%
\end{array}%
$ & $%
\begin{array}{c}
x \\ 
y \\ 
z%
\end{array}%
$ & $%
\begin{array}{c}
x \\ 
y \\ 
z%
\end{array}%
$ & $\cdots $ & $%
\begin{array}{c}
x \\ 
y \\ 
z%
\end{array}%
$ & $%
\begin{array}{c}
z \\ 
y \\ 
x%
\end{array}%
$ & $%
\begin{array}{c}
x \\ 
z \\ 
y%
\end{array}%
$ \\ \cline{2-8}
\end{tabular}%
\medskip \newline
or 2 manipulates from $u9$ to $u10$. But $u10$ is an $(n-1)$-variant of $u7$%
:\medskip

\begin{tabular}{c|c|c|c|c|c|c|c|}
\cline{2-8}
$u7$ & $1$ & $2$ & $3$ & $\cdots $ & $n-2$ & $n-1$ & $n$ \\ \cline{2-8}
& $%
\begin{array}{c}
z \\ 
x \\ 
y%
\end{array}%
$ & $%
\begin{array}{c}
x \\ 
y \\ 
z%
\end{array}%
$ & $%
\begin{array}{c}
x \\ 
y \\ 
z%
\end{array}%
$ & $\cdots $ & $%
\begin{array}{c}
x \\ 
y \\ 
z%
\end{array}%
$ & $%
\begin{array}{c}
y \\ 
z \\ 
x%
\end{array}%
$ & $%
\begin{array}{c}
x \\ 
z \\ 
y%
\end{array}%
$ \\ \cline{2-8}
\end{tabular}%
\medskip \newline
If $g(u7)=x$, then $n-1$ manipulates from $u7$ to $u10$. Therefore, $%
g(u7)\neq x$. Combining, $g(u7)=z$.\medskip

\textbf{Main Thread}\medskip

\qquad By Result \#1, $x$ is chosen at $u3$:\medskip

\begin{tabular}{c|c|c|c|c|c|c|c|}
\cline{2-8}
$u3$ & $1$ & $2$ & $3$ & $\cdots $ & $n-2$ & $n-1$ & $n$ \\ \cline{2-8}
& $%
\begin{array}{c}
z \\ 
x \\ 
y%
\end{array}%
$ & $%
\begin{array}{c}
x \\ 
y \\ 
z%
\end{array}%
$ & $%
\begin{array}{c}
x \\ 
y \\ 
z%
\end{array}%
$ & $\cdots $ & $%
\begin{array}{c}
x \\ 
y \\ 
z%
\end{array}%
$ & $%
\begin{array}{c}
y \\ 
z \\ 
x%
\end{array}%
$ & $%
\begin{array}{c}
x \\ 
y \\ 
z%
\end{array}%
$ \\ \cline{2-8}
\end{tabular}%
\medskip \newline
and by Result \#2, $z$ is chosen at $u7$:\medskip

\begin{tabular}{c|c|c|c|c|c|c|c|}
\cline{2-8}
$u7$ & $1$ & $2$ & $3$ & $\cdots $ & $n-2$ & $n-1$ & $n$ \\ \cline{2-8}
& $%
\begin{array}{c}
z \\ 
x \\ 
y%
\end{array}%
$ & $%
\begin{array}{c}
x \\ 
y \\ 
z%
\end{array}%
$ & $%
\begin{array}{c}
x \\ 
y \\ 
z%
\end{array}%
$ & $\cdots $ & $%
\begin{array}{c}
x \\ 
y \\ 
z%
\end{array}%
$ & $%
\begin{array}{c}
y \\ 
z \\ 
x%
\end{array}%
$ & $%
\begin{array}{c}
x \\ 
z \\ 
y%
\end{array}%
$ \\ \cline{2-8}
\end{tabular}%
\medskip \newline
Then $n$ manipulates from $u7$ to $u3$ and $g$ violates strategy-proofness.
This establishes that $x$ is not selected at $L1$ if $x$ is not in the range
of $g|NP^{\ast }$.\medskip

\textbf{Section 3-4}. $L2$.

We assume that $x$ is chosen at $L2$:\medskip

\begin{tabular}{c|c|c|c|c|c|c|c|}
\cline{2-8}
$L2$ & $1$ & $2$ & $3$ & $\cdots $ & $n-2$ & $n-1$ & $n$ \\ \cline{2-8}
& $%
\begin{array}{c}
z \\ 
x \\ 
y%
\end{array}%
$ & $%
\begin{array}{c}
x \\ 
z \\ 
y%
\end{array}%
$ & $%
\begin{array}{c}
x \\ 
z \\ 
y%
\end{array}%
$ & $\cdots $ & $%
\begin{array}{c}
x \\ 
z \\ 
y%
\end{array}%
$ & $%
\begin{array}{c}
y \\ 
x \\ 
z%
\end{array}%
$ & $%
\begin{array}{c}
x \\ 
z \\ 
y%
\end{array}%
$ \\ \cline{2-8}
\end{tabular}%
\newline

and seek a contradiction to strategy-proofness of $g$. This is more
complicated than the analysis for $L1$ or $L3$ since $L2(1)$ is different
from $L2(2)$, ..., and $L2(n-2)$ whereas for $L1$ and $L3$ all of the first $%
n-2$ individuals have the same ordering.\medskip

We first explain the organization of the proof, laying out a set of cases to
be considered, and then later fill in an analysis of each case and
subcase.\medskip

\textbf{Case }$L2$.A: A subset $S$ of $\{1,2,...,n-2\}$ is a minimal
decisive set for $y$ against $z$ for $g^{\ast }$ and so also for $g|NP^{\ast
}$ and the smallest such minimal subset has at least two individuals.
Because L2(1) is different, we must now distinguish between two subcases
according to whether or not $S$ contains 1.\medskip\ 

\textbf{Subcase }$L2$.A.1: A subset $S$ of $\{2,...,n-2\}$, say $%
S=\{2,...,k\}$ for $3\leq k\leq n-2$ of at least two individuals is a
minimal decisive set for $y$ against $z$ for $g|NP^{\ast }$, i.e., $S$ does
not contain 1. This requires $n\geq 5$.\medskip

\qquad We also need to consider possible coalitions $T$ decisive for $z$
against $y$ for $g|NP^{\ast }$. Since $\{1\}$ is not decisive for $y$
against $z$ for $g^{\ast }$, we see $\{2,3,...,n-1,n\}$ is decisive for $z$
against $y$ for $g|NP^{\ast }$.\medskip

\textbf{Subcase }$L2$.A.2: A subset $S=\{1,...,k\}$ for $2\leq k\leq n-2$ is
a minimal decisive set for $y$ against $z$ for $g^{\ast }$ and so also for $%
g|NP^{\ast }$, i.e., $S$ does contain 1 and 2). Within this subcase we must
again consider possibilities regarding coalitions $T$ decisive for $z$
against $y$ for $g$. But we know some coalitions decisive for $z$ against $y$
due to the minimality of S: For example, $\{2,k+1,...,n-1\}$ is decisive for 
$z$ against $y$ for $g^{\ast }$ and so $\{2,3,...,n\}$ is decisive for $z$
against $y$ for $g|NP^{\ast }$.\medskip

\textbf{Case }$L2$.B: A singleton from $\{1,2,...,n-2\}$ is a minimal
decisive set for $g|NP^{\ast }$. We again distinguish two subcases.\medskip

\qquad \textbf{Subcase }$L2$B.1: $\{1\}$ is a minimal decisive set for $y$
against $z$ for $g|NP^{\ast }$.\medskip

\qquad We also need to consider possible coalitions $T$ decisive for $z$
against $y$ for $g|NP^{\ast }$.\medskip

\qquad \textbf{Possibility a}. 2 is excluded and \{1,3,..., n-1\} is
decisive for $z$ against

$\qquad \qquad \qquad \qquad \qquad y$ for $g^{\ast }$ and so $\{1,3,...,n\}$
is decisive for $z$ against

$\qquad \qquad \qquad \qquad \qquad y$ for $g|NP^{\ast }$.\medskip

\qquad \textbf{Possibility b}. $n-1$ is excluded and $\{1,2,...,n-2\}$ is
decisive for

$\qquad \qquad \qquad \qquad \qquad z$ against $y$ for for $g^{\ast }$ and
so also for $g|NP^{\ast }$.\medskip

\qquad \textbf{Subcase }$L2$B.2: A singleton from $\{2,...,n-2\}$, say $%
\{2\} $, is a minimal decisive set for $y$ against $z$ for $g|NP^{\ast }$%
.\medskip

\qquad We also need to consider possible coalitions $T$ decisive for $z$
against $y$ for $g|NP^{\ast }$. Alternative 2 must also be in $T$; and some
element can be excluded from $T$ since we can't have Pareto
domination.\medskip

\qquad \textbf{Possibility a}. 1 is excluded and $\{2,...,n\}$ is decisive
for $z$ against

$\qquad \qquad \qquad \qquad \qquad y$ for $g|NP^{\ast }$.\medskip

\qquad \textbf{Possibility b}. An individual in $\{3,...,n-2\}$, say 3, is
excluded

\qquad \qquad \qquad \qquad \qquad and $\{1,3,...,n\}$ is decisive for $z$
against $y$ for $g$

\qquad \qquad \qquad \qquad \qquad on $g|NP^{\ast }$.\medskip

\qquad \textbf{Possibility c}. $n-1$ is excluded and $\{1,2,...,n-2\}$ is
decisive

\qquad \qquad \qquad \qquad \qquad for $z$ against $y$ for $g|NP^{\ast }$%
.\medskip

\textbf{Case }$L2$.C: The coalition $\{n-1,n\}$ is a (minimal) decisive set
for $z$ against $y$ for $g|NP^{\ast }$.\medskip

Now we analyze each subcase in turn.\medskip

\textbf{Subcase }$L2$.A.1: $\{2,...,k\}$ is decisive for $y$ against $z$ and 
$\{2,...,n\}$ is decisive for $z$ against $y$ for $g|NP^{\ast }$.\medskip

\textbf{Result 1}\medskip

At profile $u1$:\medskip

\begin{tabular}{c|c|c|c|c|c|c|c|c|c|c|}
\cline{2-11}
$u1$ & $1$ & $2$ & $3$ & $\cdots $ & $k$ & $k+1$ & $\cdots $ & $n-2$ & $n-1$
& $n$ \\ \cline{2-11}
& $%
\begin{array}{c}
z \\ 
x \\ 
y%
\end{array}%
$ & $%
\begin{array}{c}
x \\ 
y \\ 
z%
\end{array}%
$ & $%
\begin{array}{c}
x \\ 
y \\ 
z%
\end{array}%
$ & $\cdots $ & $%
\begin{array}{c}
x \\ 
y \\ 
z%
\end{array}%
$ & $%
\begin{array}{c}
x \\ 
z \\ 
y%
\end{array}%
$ & $\cdots $ & $%
\begin{array}{c}
x \\ 
z \\ 
y%
\end{array}%
$ & $%
\begin{array}{c}
z \\ 
y \\ 
x%
\end{array}%
$ & $%
\begin{array}{c}
z \\ 
y \\ 
x%
\end{array}%
$ \\ \cline{2-11}
\end{tabular}%
\medskip \newline
we have $g(u1)=y$ since $\{2,...,k\}$ is decisive for $y$ against $z$ for $%
g|NP^{\ast }$. Then at $(n-1)$-variant profile $u2$:\medskip

\begin{tabular}{c|c|c|c|c|c|c|c|c|c|c|}
\cline{2-11}
$u2$ & $1$ & $2$ & $3$ & $\cdots $ & $k$ & $k+1$ & $\cdots $ & $n-2$ & $n-1$
& $n$ \\ \cline{2-11}
& $%
\begin{array}{c}
z \\ 
x \\ 
y%
\end{array}%
$ & $%
\begin{array}{c}
x \\ 
y \\ 
z%
\end{array}%
$ & $%
\begin{array}{c}
x \\ 
y \\ 
z%
\end{array}%
$ & $\cdots $ & $%
\begin{array}{c}
x \\ 
y \\ 
z%
\end{array}%
$ & $%
\begin{array}{c}
x \\ 
z \\ 
y%
\end{array}%
$ & $\cdots $ & $%
\begin{array}{c}
x \\ 
z \\ 
y%
\end{array}%
$ & $%
\begin{array}{c}
y \\ 
z \\ 
x%
\end{array}%
$ & $%
\begin{array}{c}
z \\ 
y \\ 
x%
\end{array}%
$ \\ \cline{2-11}
\end{tabular}%
\medskip \newline
$y$ is also chosen or else $n-1$ will manipulate from $u2$ to $u1$. Next, at 
$n$-variant profile $u3$,\medskip

\begin{tabular}{c|c|c|c|c|c|c|c|c|c|c|}
\cline{2-11}
$u3$ & $1$ & $2$ & $3$ & $\cdots $ & $k$ & $k+1$ & $\cdots $ & $n-2$ & $n-1$
& $n$ \\ \cline{2-11}
& $%
\begin{array}{c}
z \\ 
x \\ 
y%
\end{array}%
$ & $%
\begin{array}{c}
x \\ 
y \\ 
z%
\end{array}%
$ & $%
\begin{array}{c}
x \\ 
y \\ 
z%
\end{array}%
$ & $\cdots $ & $%
\begin{array}{c}
x \\ 
y \\ 
z%
\end{array}%
$ & $%
\begin{array}{c}
x \\ 
z \\ 
y%
\end{array}%
$ & $\cdots $ & $%
\begin{array}{c}
x \\ 
z \\ 
y%
\end{array}%
$ & $%
\begin{array}{c}
y \\ 
z \\ 
x%
\end{array}%
$ & $%
\begin{array}{c}
x \\ 
z \\ 
y%
\end{array}%
$ \\ \cline{2-11}
\end{tabular}%
\medskip \newline
we have $g(u3$) $\neq $ $z$ or $n$ would manipulate from $u2$ to $u3$.
Similarly, at $n$-variant $u4$,\medskip

\begin{tabular}{c|c|c|c|c|c|c|c|c|c|c|}
\cline{2-11}
$u4$ & $1$ & $2$ & $3$ & $\cdots $ & $k$ & $k+1$ & $\cdots $ & $n-2$ & $n-1$
& $n$ \\ \cline{2-11}
& $%
\begin{array}{c}
z \\ 
x \\ 
y%
\end{array}%
$ & $%
\begin{array}{c}
x \\ 
y \\ 
z%
\end{array}%
$ & $%
\begin{array}{c}
x \\ 
y \\ 
z%
\end{array}%
$ & $\cdots $ & $%
\begin{array}{c}
x \\ 
y \\ 
z%
\end{array}%
$ & $%
\begin{array}{c}
x \\ 
z \\ 
y%
\end{array}%
$ & $\cdots $ & $%
\begin{array}{c}
x \\ 
z \\ 
y%
\end{array}%
$ & $%
\begin{array}{c}
y \\ 
z \\ 
x%
\end{array}%
$ & $%
\begin{array}{c}
x \\ 
y \\ 
z%
\end{array}%
$ \\ \cline{2-11}
\end{tabular}%
\medskip \newline
We have $g(u4)\neq z$ or $n$ would manipulate from $u2$ to $u4$.\medskip

\textbf{Result 2}\medskip

We are assuming $x$ is chosen at $L2$:\medskip

\begin{tabular}{c|c|c|c|c|c|c|c|c|c|c|}
\cline{2-11}
$L2$ & $1$ & $2$ & $3$ & $\cdots $ & $k$ & $k+1$ & $\cdots $ & $n-2$ & $n-1$
& $n$ \\ \cline{2-11}
& $%
\begin{array}{c}
z \\ 
x \\ 
y%
\end{array}%
$ & $%
\begin{array}{c}
x \\ 
z \\ 
y%
\end{array}%
$ & $%
\begin{array}{c}
x \\ 
z \\ 
y%
\end{array}%
$ & $\cdots $ & $%
\begin{array}{c}
x \\ 
z \\ 
y%
\end{array}%
$ & $%
\begin{array}{c}
x \\ 
z \\ 
y%
\end{array}%
$ & $\cdots $ & $%
\begin{array}{c}
x \\ 
z \\ 
y%
\end{array}%
$ & $%
\begin{array}{c}
y \\ 
x \\ 
z%
\end{array}%
$ & $%
\begin{array}{c}
x \\ 
z \\ 
y%
\end{array}%
$ \\ \cline{2-11}
\end{tabular}%
\medskip \newline
Then $x$ is also chosen at $L2^{\ast }$:\medskip

\begin{tabular}{c|c|c|c|c|c|c|c|c|c|c|}
\cline{2-11}
$L2^{\ast }$ & $1$ & $2$ & $3$ & $\cdots $ & $k$ & $k+1$ & $\cdots $ & $n-2$
& $n-1$ & $n$ \\ \cline{2-11}
& $%
\begin{array}{c}
z \\ 
x \\ 
y%
\end{array}%
$ & $%
\begin{array}{c}
x \\ 
y \\ 
z%
\end{array}%
$ & $%
\begin{array}{c}
x \\ 
y \\ 
z%
\end{array}%
$ & $\cdots $ & $%
\begin{array}{c}
x \\ 
y \\ 
z%
\end{array}%
$ & $%
\begin{array}{c}
x \\ 
z \\ 
y%
\end{array}%
$ & $\cdots $ & $%
\begin{array}{c}
x \\ 
z \\ 
y%
\end{array}%
$ & $%
\begin{array}{c}
y \\ 
x \\ 
z%
\end{array}%
$ & $%
\begin{array}{c}
x \\ 
z \\ 
y%
\end{array}%
$ \\ \cline{2-11}
\end{tabular}%
\medskip \newline
obtained by interchanging $y$ and $z$ below $x$ for individuals $2,...,k$.
In turn, this implies that at $(n-1)$-variant profile $u3$, we get $%
g(u3)\neq y$ or $n-1$ manipulates from $L2^{\ast }$ to $u3$. Then $g(u3)=x$%
.\medskip

\textbf{Result 3}\medskip

From earlier analysis, we know $x$ is not chosen at $L1$:\medskip

\begin{tabular}{c|c|c|c|c|c|c|c|c|c|c|}
\cline{2-11}
$L1$ & $1$ & $2$ & $3$ & $\cdots $ & $k$ & $k+1$ & $\cdots $ & $n-2$ & $n-1$
& $n$ \\ \cline{2-11}
& $%
\begin{array}{c}
x \\ 
y \\ 
z%
\end{array}%
$ & $%
\begin{array}{c}
x \\ 
y \\ 
z%
\end{array}%
$ & $%
\begin{array}{c}
x \\ 
y \\ 
z%
\end{array}%
$ & $\cdots $ & $%
\begin{array}{c}
x \\ 
y \\ 
z%
\end{array}%
$ & $%
\begin{array}{c}
x \\ 
y \\ 
z%
\end{array}%
$ & $\cdots $ & $%
\begin{array}{c}
x \\ 
y \\ 
z%
\end{array}%
$ & $%
\begin{array}{c}
z \\ 
y \\ 
x%
\end{array}%
$ & $%
\begin{array}{c}
x \\ 
y \\ 
z%
\end{array}%
$ \\ \cline{2-11}
\end{tabular}%
\medskip \newline
Then $x$ is also not chosen at profile $L1^{\ast }$, obtained from $L1$ by
interchanging $y$ and $z$ below $x$ for individuals $t+1,...,n-2$:\medskip

\begin{tabular}{c|c|c|c|c|c|c|c|c|c|c|}
\cline{2-11}
$L1^{\ast }$ & $1$ & $2$ & $3$ & $\cdots $ & $k$ & $k+1$ & $\cdots $ & $n-2$
& $n-1$ & $n$ \\ \cline{2-11}
& $%
\begin{array}{c}
x \\ 
y \\ 
z%
\end{array}%
$ & $%
\begin{array}{c}
x \\ 
y \\ 
z%
\end{array}%
$ & $%
\begin{array}{c}
x \\ 
y \\ 
z%
\end{array}%
$ & $\cdots $ & $%
\begin{array}{c}
x \\ 
y \\ 
z%
\end{array}%
$ & $%
\begin{array}{c}
x \\ 
z \\ 
y%
\end{array}%
$ & $\cdots $ & $%
\begin{array}{c}
x \\ 
z \\ 
y%
\end{array}%
$ & $%
\begin{array}{c}
z \\ 
y \\ 
x%
\end{array}%
$ & $%
\begin{array}{c}
x \\ 
y \\ 
z%
\end{array}%
$ \\ \cline{2-11}
\end{tabular}%
\medskip \newline
So also $x$ is not chosen at 1-variant profile $u5$:\medskip

\begin{tabular}{c|c|c|c|c|c|c|c|c|c|c|}
\cline{2-11}
$u5$ & $1$ & $2$ & $3$ & $\cdots $ & $k$ & $k+1$ & $\cdots $ & $n-2$ & $n-1$
& $n$ \\ \cline{2-11}
& $%
\begin{array}{c}
z \\ 
x \\ 
y%
\end{array}%
$ & $%
\begin{array}{c}
x \\ 
y \\ 
z%
\end{array}%
$ & $%
\begin{array}{c}
x \\ 
y \\ 
z%
\end{array}%
$ & $\cdots $ & $%
\begin{array}{c}
x \\ 
y \\ 
z%
\end{array}%
$ & $%
\begin{array}{c}
x \\ 
z \\ 
y%
\end{array}%
$ & $\cdots $ & $%
\begin{array}{c}
x \\ 
z \\ 
y%
\end{array}%
$ & $%
\begin{array}{c}
z \\ 
y \\ 
x%
\end{array}%
$ & $%
\begin{array}{c}
x \\ 
y \\ 
z%
\end{array}%
$ \\ \cline{2-11}
\end{tabular}%
\medskip \newline
or 1 would manipulate from $L1^{\ast }$ to $u5$. Then $g(u4)\neq x$ or $n-1$
manipulates from $u4$ to $u5$.\medskip

\textbf{Main Thread}\medskip

By Result 1 and Result 2, $x$ is chosen at $u3$:\medskip

\begin{tabular}{c|c|c|c|c|c|c|c|c|c|c|}
\cline{2-11}
$u3$ & $1$ & $2$ & $3$ & $\cdots $ & $k$ & $k+1$ & $\cdots $ & $n-2$ & $n-1$
& $n$ \\ \cline{2-11}
& $%
\begin{array}{c}
z \\ 
x \\ 
y%
\end{array}%
$ & $%
\begin{array}{c}
x \\ 
y \\ 
z%
\end{array}%
$ & $%
\begin{array}{c}
x \\ 
y \\ 
z%
\end{array}%
$ & $\cdots $ & $%
\begin{array}{c}
x \\ 
y \\ 
z%
\end{array}%
$ & $%
\begin{array}{c}
x \\ 
z \\ 
y%
\end{array}%
$ & $\cdots $ & $%
\begin{array}{c}
x \\ 
z \\ 
y%
\end{array}%
$ & $%
\begin{array}{c}
y \\ 
z \\ 
x%
\end{array}%
$ & $%
\begin{array}{c}
x \\ 
z \\ 
y%
\end{array}%
$ \\ \cline{2-11}
\end{tabular}%
\medskip \newline
By Result 1 and Result 3, $y$ is chosen at $u4$:\medskip

\begin{tabular}{c|c|c|c|c|c|c|c|c|c|c|}
\cline{2-11}
$u4$ & $1$ & $2$ & $3$ & $\cdots $ & $k$ & $k+1$ & $\cdots $ & $n-2$ & $n-1$
& $n$ \\ \cline{2-11}
& $%
\begin{array}{c}
z \\ 
x \\ 
y%
\end{array}%
$ & $%
\begin{array}{c}
x \\ 
y \\ 
z%
\end{array}%
$ & $%
\begin{array}{c}
x \\ 
y \\ 
z%
\end{array}%
$ & $\cdots $ & $%
\begin{array}{c}
x \\ 
y \\ 
z%
\end{array}%
$ & $%
\begin{array}{c}
x \\ 
z \\ 
y%
\end{array}%
$ & $\cdots $ & $%
\begin{array}{c}
x \\ 
z \\ 
y%
\end{array}%
$ & $%
\begin{array}{c}
y \\ 
z \\ 
x%
\end{array}%
$ & $%
\begin{array}{c}
x \\ 
y \\ 
z%
\end{array}%
$ \\ \cline{2-11}
\end{tabular}%
\medskip \newline
But then $n$ manipulates from $u4$ to $u3$, and $g$ violates
strategy-proofness.\medskip

\textbf{Subcase }$L2$A.2: A subset $S$ = \{1, ..., k\} for 2 $\leq $ k $\leq 
$ $n-2$ is a minimal decisive set for $y$ against $z$ for $g|NP^{\ast }$,
i.e., $S$ does contain 1 and 2). Within this subcase we must again consider
possibilities regarding coalitions $T$ decisive for $z$ against $y$ for $%
g|NP^{\ast }$. But we know some coalitions decisive for $z$ against $y$ due
to the minimality of S: For example, $\{2,k+1,...,n-1\}$ is decisive for $z$
against $y$ for $g|NP^{\ast }$.\medskip

\textbf{Result 1}.\medskip

At $L1$:\medskip

\begin{tabular}{c|c|c|c|c|c|c|c|c|c|c|}
\cline{2-11}
$L1$ & $1$ & $2$ & $3$ & $\cdots $ & $k$ & $k+1$ & $\cdots $ & $n-2$ & $n-1$
& $n$ \\ \cline{2-11}
& $%
\begin{array}{c}
x \\ 
y \\ 
z%
\end{array}%
$ & $%
\begin{array}{c}
x \\ 
y \\ 
z%
\end{array}%
$ & $%
\begin{array}{c}
x \\ 
y \\ 
z%
\end{array}%
$ & $\cdots $ & $%
\begin{array}{c}
x \\ 
y \\ 
z%
\end{array}%
$ & $%
\begin{array}{c}
x \\ 
y \\ 
z%
\end{array}%
$ & $\cdots $ & $%
\begin{array}{c}
x \\ 
y \\ 
z%
\end{array}%
$ & $%
\begin{array}{c}
z \\ 
y \\ 
x%
\end{array}%
$ & $%
\begin{array}{c}
x \\ 
y \\ 
z%
\end{array}%
$ \\ \cline{2-11}
\end{tabular}%
\medskip \newline
we have $g(L1)\neq x$ by our previous analysis. \ Then $x$ is also not
chosen at $L1^{\ast }$, obtained by reversing $y$ and $z$ below $x$ for $%
k+1,...,n-2$ or else a standard sequence argument will show a violation of
strategy-proofness.\medskip

\begin{tabular}{c|c|c|c|c|c|c|c|c|c|c|}
\cline{2-11}
$L1^{\ast }$ & $1$ & $2$ & $3$ & $\cdots $ & $k$ & $k+1$ & $\cdots $ & $n-2$
& $n-1$ & $n$ \\ \cline{2-11}
& $%
\begin{array}{c}
x \\ 
y \\ 
z%
\end{array}%
$ & $%
\begin{array}{c}
x \\ 
y \\ 
z%
\end{array}%
$ & $%
\begin{array}{c}
x \\ 
y \\ 
z%
\end{array}%
$ & $\cdots $ & $%
\begin{array}{c}
x \\ 
y \\ 
z%
\end{array}%
$ & $%
\begin{array}{c}
x \\ 
z \\ 
y%
\end{array}%
$ & $\cdots $ & $%
\begin{array}{c}
x \\ 
z \\ 
y%
\end{array}%
$ & $%
\begin{array}{c}
z \\ 
y \\ 
x%
\end{array}%
$ & $%
\begin{array}{c}
x \\ 
y \\ 
z%
\end{array}%
$ \\ \cline{2-11}
\end{tabular}%
\medskip \newline
Then at 1-variant profile $u1$:\medskip

\begin{tabular}{c|c|c|c|c|c|c|c|c|c|c|}
\cline{2-11}
$u1$ & $1$ & $2$ & $3$ & $\cdots $ & $k$ & $k+1$ & $\cdots $ & $n-2$ & $n-1$
& $n$ \\ \cline{2-11}
& $%
\begin{array}{c}
z \\ 
x \\ 
y%
\end{array}%
$ & $%
\begin{array}{c}
x \\ 
y \\ 
z%
\end{array}%
$ & $%
\begin{array}{c}
x \\ 
y \\ 
z%
\end{array}%
$ & $\cdots $ & $%
\begin{array}{c}
x \\ 
y \\ 
z%
\end{array}%
$ & $%
\begin{array}{c}
x \\ 
z \\ 
y%
\end{array}%
$ & $\cdots $ & $%
\begin{array}{c}
x \\ 
z \\ 
y%
\end{array}%
$ & $%
\begin{array}{c}
z \\ 
y \\ 
x%
\end{array}%
$ & $%
\begin{array}{c}
x \\ 
y \\ 
z%
\end{array}%
$ \\ \cline{2-11}
\end{tabular}%
\medskip \newline
we have $g(u1)\neq x$ or 1 manipulates from $L1^{\ast }$ to $u1$. Then at $%
(n-1)$-variant $u2$:\medskip

\begin{tabular}{c|c|c|c|c|c|c|c|c|c|c|}
\cline{2-11}
$u2$ & $1$ & $2$ & $3$ & $\cdots $ & $k$ & $k+1$ & $\cdots $ & $n-2$ & $n-1$
& $n$ \\ \cline{2-11}
& $%
\begin{array}{c}
z \\ 
x \\ 
y%
\end{array}%
$ & $%
\begin{array}{c}
x \\ 
y \\ 
z%
\end{array}%
$ & $%
\begin{array}{c}
x \\ 
y \\ 
z%
\end{array}%
$ & $\cdots $ & $%
\begin{array}{c}
x \\ 
y \\ 
z%
\end{array}%
$ & $%
\begin{array}{c}
x \\ 
z \\ 
y%
\end{array}%
$ & $\cdots $ & $%
\begin{array}{c}
x \\ 
z \\ 
y%
\end{array}%
$ & $%
\begin{array}{c}
y \\ 
z \\ 
x%
\end{array}%
$ & $%
\begin{array}{c}
x \\ 
y \\ 
z%
\end{array}%
$ \\ \cline{2-11}
\end{tabular}%
\medskip \newline
we have $g(u2)\neq x$ or $n-1$ manipulates from $u2$ to $u1$.\medskip

Next, at $n$-variant $u3$:\medskip

\begin{tabular}{c|c|c|c|c|c|c|c|c|c|c|}
\cline{2-11}
$u3$ & $1$ & $2$ & $3$ & $\cdots $ & $k$ & $k+1$ & $\cdots $ & $n-2$ & $n-1$
& $n$ \\ \cline{2-11}
& $%
\begin{array}{c}
z \\ 
x \\ 
y%
\end{array}%
$ & $%
\begin{array}{c}
x \\ 
y \\ 
z%
\end{array}%
$ & $%
\begin{array}{c}
x \\ 
y \\ 
z%
\end{array}%
$ & $\cdots $ & $%
\begin{array}{c}
x \\ 
y \\ 
z%
\end{array}%
$ & $%
\begin{array}{c}
x \\ 
z \\ 
y%
\end{array}%
$ & $\cdots $ & $%
\begin{array}{c}
x \\ 
z \\ 
y%
\end{array}%
$ & $%
\begin{array}{c}
y \\ 
z \\ 
x%
\end{array}%
$ & $%
\begin{array}{c}
x \\ 
z \\ 
y%
\end{array}%
$ \\ \cline{2-11}
\end{tabular}%
\medskip \newline
we get $g(u3)\neq x$ or $n$ manipulates from $u2$ to $u3$.\medskip

\textbf{Result 2}.\medskip

\qquad We are assuming $x$ is chosen at profile $L2$:\medskip

\medskip 
\begin{tabular}{c|c|c|c|c|c|c|c|c|c|c|}
\cline{2-11}
$L2$ & $1$ & $2$ & $3$ & $\cdots $ & $k$ & $k+1$ & $\cdots $ & $n-2$ & $n-1$
& $n$ \\ \cline{2-11}
& $%
\begin{array}{c}
z \\ 
x \\ 
y%
\end{array}%
$ & $%
\begin{array}{c}
x \\ 
z \\ 
y%
\end{array}%
$ & $%
\begin{array}{c}
x \\ 
z \\ 
y%
\end{array}%
$ & $\cdots $ & $%
\begin{array}{c}
x \\ 
z \\ 
y%
\end{array}%
$ & $%
\begin{array}{c}
x \\ 
z \\ 
y%
\end{array}%
$ & $\cdots $ & $%
\begin{array}{c}
x \\ 
z \\ 
y%
\end{array}%
$ & $%
\begin{array}{c}
y \\ 
x \\ 
z%
\end{array}%
$ & $%
\begin{array}{c}
x \\ 
z \\ 
y%
\end{array}%
$ \\ \cline{2-11}
\end{tabular}%
\medskip \newline
Then at profile $u4$, obtained by interchanging $y$ and $z$ below $x$ for 2
to $k$:\medskip

\begin{tabular}{c|c|c|c|c|c|c|c|c|c|c|}
\cline{2-11}
$u4$ & $1$ & $2$ & $3$ & $\cdots $ & $k$ & $k+1$ & $\cdots $ & $n-2$ & $n-1$
& $n$ \\ \cline{2-11}
& $%
\begin{array}{c}
z \\ 
x \\ 
y%
\end{array}%
$ & $%
\begin{array}{c}
x \\ 
y \\ 
z%
\end{array}%
$ & $%
\begin{array}{c}
x \\ 
y \\ 
z%
\end{array}%
$ & $\cdots $ & $%
\begin{array}{c}
x \\ 
y \\ 
z%
\end{array}%
$ & $%
\begin{array}{c}
x \\ 
z \\ 
y%
\end{array}%
$ & $\cdots $ & $%
\begin{array}{c}
x \\ 
z \\ 
y%
\end{array}%
$ & $%
\begin{array}{c}
y \\ 
x \\ 
z%
\end{array}%
$ & $%
\begin{array}{c}
x \\ 
z \\ 
y%
\end{array}%
$ \\ \cline{2-11}
\end{tabular}%
\medskip \newline
$g(u4)=x$ by strategy-proofness.\medskip

\qquad Profile $u4$ is also an $(n-1)$-variant of $u3$, and $g(u4)=x$
implies $g(u3)\neq y$ or $n-1$ manipulates from $u4$ to $u3$. Hence $g(u3)=z$%
.\medskip

\textbf{Result 3}.\medskip

\qquad At profile $u5$:\medskip

\begin{tabular}{c|c|c|c|c|c|c|c|c|c|c|}
\cline{2-11}
$u5$ & $1$ & $2$ & $3$ & $\cdots $ & $k$ & $k+1$ & $\cdots $ & $n-2$ & $n-1$
& $n$ \\ \cline{2-11}
& $%
\begin{array}{c}
z \\ 
y \\ 
x%
\end{array}%
$ & $%
\begin{array}{c}
x \\ 
y \\ 
z%
\end{array}%
$ & $%
\begin{array}{c}
x \\ 
y \\ 
z%
\end{array}%
$ & $\cdots $ & $%
\begin{array}{c}
x \\ 
y \\ 
z%
\end{array}%
$ & $%
\begin{array}{c}
x \\ 
z \\ 
y%
\end{array}%
$ & $\cdots $ & $%
\begin{array}{c}
x \\ 
z \\ 
y%
\end{array}%
$ & $%
\begin{array}{c}
x \\ 
z \\ 
y%
\end{array}%
$ & $%
\begin{array}{c}
x \\ 
z \\ 
y%
\end{array}%
$ \\ \cline{2-11}
\end{tabular}%
\medskip \newline
$g(u5)=z$ since $\{1,...,k\}$ is minimal decisive for $y$ against $z$ on $%
NP^{\ast }$. Then at $(n-1)$-variant profile $u6$:\medskip

\begin{tabular}{c|c|c|c|c|c|c|c|c|c|c|}
\cline{2-11}
$u6$ & $1$ & $2$ & $3$ & $\cdots $ & $k$ & $k+1$ & $\cdots $ & $n-2$ & $n-1$
& $n$ \\ \cline{2-11}
& $%
\begin{array}{c}
z \\ 
y \\ 
x%
\end{array}%
$ & $%
\begin{array}{c}
x \\ 
y \\ 
z%
\end{array}%
$ & $%
\begin{array}{c}
x \\ 
y \\ 
z%
\end{array}%
$ & $\cdots $ & $%
\begin{array}{c}
x \\ 
y \\ 
z%
\end{array}%
$ & $%
\begin{array}{c}
x \\ 
z \\ 
y%
\end{array}%
$ & $\cdots $ & $%
\begin{array}{c}
x \\ 
z \\ 
y%
\end{array}%
$ & $%
\begin{array}{c}
y \\ 
x \\ 
z%
\end{array}%
$ & $%
\begin{array}{c}
x \\ 
z \\ 
y%
\end{array}%
$ \\ \cline{2-11}
\end{tabular}%
\medskip \newline
we have $g(u6)\neq x$ or $n-1$ manipulates from $u5$ to $u6$. Profile $u6$
is also a 1-variant of $u4$ and so $g(u6)\neq z$ or 1 manipulates from $u4$
to $u6$. Summarizing, $g(u6)=y$.\medskip

\textbf{Main Thread}\medskip

At profile $u7$, an $(n-1)$-variant of $u6$:\medskip

\begin{tabular}{c|c|c|c|c|c|c|c|c|c|c|}
\cline{2-11}
$u7$ & $1$ & $2$ & $3$ & $\cdots $ & $k$ & $k+1$ & $\cdots $ & $n-2$ & $n-1$
& $n$ \\ \cline{2-11}
& $%
\begin{array}{c}
z \\ 
y \\ 
x%
\end{array}%
$ & $%
\begin{array}{c}
x \\ 
y \\ 
z%
\end{array}%
$ & $%
\begin{array}{c}
x \\ 
y \\ 
z%
\end{array}%
$ & $\cdots $ & $%
\begin{array}{c}
x \\ 
y \\ 
z%
\end{array}%
$ & $%
\begin{array}{c}
x \\ 
z \\ 
y%
\end{array}%
$ & $\cdots $ & $%
\begin{array}{c}
x \\ 
z \\ 
y%
\end{array}%
$ & $%
\begin{array}{c}
y \\ 
z \\ 
x%
\end{array}%
$ & $%
\begin{array}{c}
x \\ 
z \\ 
y%
\end{array}%
$ \\ \cline{2-11}
\end{tabular}%
\medskip \newline
$g(u7)=y$ or $n-1$ manipulates from $u7$ to $u6$. But $u7$ is also a 1
variant of $u3$:\medskip

\begin{tabular}{c|c|c|c|c|c|c|c|c|c|c|}
\cline{2-11}
$u3$ & $1$ & $2$ & $3$ & $\cdots $ & $k$ & $k+1$ & $\cdots $ & $n-2$ & $n-1$
& $n$ \\ \cline{2-11}
& $%
\begin{array}{c}
z \\ 
x \\ 
y%
\end{array}%
$ & $%
\begin{array}{c}
x \\ 
y \\ 
z%
\end{array}%
$ & $%
\begin{array}{c}
x \\ 
y \\ 
z%
\end{array}%
$ & $\cdots $ & $%
\begin{array}{c}
x \\ 
y \\ 
z%
\end{array}%
$ & $%
\begin{array}{c}
x \\ 
z \\ 
y%
\end{array}%
$ & $\cdots $ & $%
\begin{array}{c}
x \\ 
z \\ 
y%
\end{array}%
$ & $%
\begin{array}{c}
y \\ 
z \\ 
x%
\end{array}%
$ & $%
\begin{array}{c}
x \\ 
z \\ 
y%
\end{array}%
$ \\ \cline{2-11}
\end{tabular}%
\medskip \newline
and $g(u3)=z$ from Result 1 combined with Result 2. Therefore 1 manipulates
from $u7$ to $u3$, and $g$ violates strategy-proofness.\medskip

\textbf{Case }$L2$.B: A singleton from $\{1,2,...,n-2\}$ is a minimal
decisive set for $y$ against $z$ for $g|NP^{\ast }$. We again distinguish
two subcases.\medskip

\qquad \textbf{Subcase }$L2$B.1: $\{1\}$ is a minimal decisive set for $y$
against $z$ for $g|NP^{\ast }$.\medskip

\qquad We also need to consider possible coalitions $T$ decisive for $z$
against $y$ for $g|NP^{\ast }$.\medskip

\qquad Possibility a. 2 is excluded and $\{1,3,...,n\}$ is decisive for $z$
against $y$

\qquad \qquad \qquad \qquad \qquad for $g|NP^{\ast }$.\medskip

\textbf{Result 1}.\medskip

\qquad At profile $L2$,\medskip

\begin{tabular}{c|c|c|c|c|c|c|c|}
\cline{2-8}
$L2$ & $1$ & $2$ & $3$ & $\cdots $ & $n-2$ & $n-1$ & $n$ \\ \cline{2-8}
& $%
\begin{array}{c}
z \\ 
x \\ 
y%
\end{array}%
$ & $%
\begin{array}{c}
x \\ 
z \\ 
y%
\end{array}%
$ & $%
\begin{array}{c}
x \\ 
z \\ 
y%
\end{array}%
$ & $\cdots $ & $%
\begin{array}{c}
x \\ 
z \\ 
y%
\end{array}%
$ & $%
\begin{array}{c}
y \\ 
x \\ 
z%
\end{array}%
$ & $%
\begin{array}{c}
x \\ 
z \\ 
y%
\end{array}%
$ \\ \cline{2-8}
\end{tabular}%
\medskip \newline
we get $g(L2)=x$ by assumption. Then at 2-variant $u1$,\medskip

\begin{tabular}{c|c|c|c|c|c|c|c|}
\cline{2-8}
$u1$ & $1$ & $2$ & $3$ & $\cdots $ & $n-2$ & $n-1$ & $n$ \\ \cline{2-8}
& $%
\begin{array}{c}
z \\ 
x \\ 
y%
\end{array}%
$ & $%
\begin{array}{c}
x \\ 
y \\ 
z%
\end{array}%
$ & $%
\begin{array}{c}
x \\ 
z \\ 
y%
\end{array}%
$ & $\cdots $ & $%
\begin{array}{c}
x \\ 
z \\ 
y%
\end{array}%
$ & $%
\begin{array}{c}
y \\ 
x \\ 
z%
\end{array}%
$ & $%
\begin{array}{c}
x \\ 
z \\ 
y%
\end{array}%
$ \\ \cline{2-8}
\end{tabular}%
\medskip \newline
we have $g(u1)=x$ or 2 manipulates from $u1$ to $L2$. Then at $(n-1)$%
-variant profile $u2$:\medskip

\begin{tabular}{c|c|c|c|c|c|c|c|}
\cline{2-8}
$u2$ & $1$ & $2$ & $3$ & $\cdots $ & $n-2$ & $n-1$ & $n$ \\ \cline{2-8}
& $%
\begin{array}{c}
z \\ 
x \\ 
y%
\end{array}%
$ & $%
\begin{array}{c}
x \\ 
y \\ 
z%
\end{array}%
$ & $%
\begin{array}{c}
x \\ 
z \\ 
y%
\end{array}%
$ & $\cdots $ & $%
\begin{array}{c}
x \\ 
z \\ 
y%
\end{array}%
$ & $%
\begin{array}{c}
y \\ 
z \\ 
x%
\end{array}%
$ & $%
\begin{array}{c}
x \\ 
z \\ 
y%
\end{array}%
$ \\ \cline{2-8}
\end{tabular}%
\medskip \newline
we get $g(u2)\neq y$ or $n-1$ manipulates from $u1$ to $u2$. At $n$-variant
profile $u3$,\medskip

\begin{tabular}{c|c|c|c|c|c|c|c|}
\cline{2-8}
$u3$ & $1$ & $2$ & $3$ & $\cdots $ & $n-2$ & $n-1$ & $n$ \\ \cline{2-8}
& $%
\begin{array}{c}
z \\ 
x \\ 
y%
\end{array}%
$ & $%
\begin{array}{c}
x \\ 
y \\ 
z%
\end{array}%
$ & $%
\begin{array}{c}
x \\ 
z \\ 
y%
\end{array}%
$ & $\cdots $ & $%
\begin{array}{c}
x \\ 
z \\ 
y%
\end{array}%
$ & $%
\begin{array}{c}
y \\ 
z \\ 
x%
\end{array}%
$ & $%
\begin{array}{c}
z \\ 
x \\ 
y%
\end{array}%
$ \\ \cline{2-8}
\end{tabular}%
\medskip \newline
we have $g(u3)\neq y$ or $n$ manipulates from $u3$ to $u2$.\medskip

\textbf{Result 2}.\medskip

\qquad At profile $L1$,\medskip

\begin{tabular}{c|c|c|c|c|c|c|c|}
\cline{2-8}
$L1$ & $1$ & $2$ & $3$ & $\cdots $ & $n-2$ & $n-1$ & $n$ \\ \cline{2-8}
& $%
\begin{array}{c}
x \\ 
y \\ 
z%
\end{array}%
$ & $%
\begin{array}{c}
x \\ 
y \\ 
z%
\end{array}%
$ & $%
\begin{array}{c}
x \\ 
y \\ 
z%
\end{array}%
$ & $\cdots $ & $%
\begin{array}{c}
x \\ 
y \\ 
z%
\end{array}%
$ & $%
\begin{array}{c}
z \\ 
y \\ 
x%
\end{array}%
$ & $%
\begin{array}{c}
x \\ 
y \\ 
z%
\end{array}%
$ \\ \cline{2-8}
\end{tabular}%
\medskip \newline
we get $g(L1)\neq x$ by an earlier analysis. Then $x$ is also not chosen at $%
L1^{\ast }$, obtained by switching $y$ and $z$ below $x$ for individuals 3,
..., $n-2$:\medskip

\begin{tabular}{c|c|c|c|c|c|c|c|}
\cline{2-8}
$L1^{\ast }$ & $1$ & $2$ & $3$ & $\cdots $ & $n-2$ & $n-1$ & $n$ \\ 
\cline{2-8}
& $%
\begin{array}{c}
x \\ 
y \\ 
z%
\end{array}%
$ & $%
\begin{array}{c}
x \\ 
y \\ 
z%
\end{array}%
$ & $%
\begin{array}{c}
x \\ 
z \\ 
y%
\end{array}%
$ & $\cdots $ & $%
\begin{array}{c}
x \\ 
z \\ 
y%
\end{array}%
$ & $%
\begin{array}{c}
z \\ 
y \\ 
x%
\end{array}%
$ & $%
\begin{array}{c}
x \\ 
y \\ 
z%
\end{array}%
$ \\ \cline{2-8}
\end{tabular}%
\medskip \newline
Then at 1-variant profile $u4$:\medskip

\begin{tabular}{c|c|c|c|c|c|c|c|}
\cline{2-8}
$u4$ & $1$ & $2$ & $3$ & $\cdots $ & $n-2$ & $n-1$ & $n$ \\ \cline{2-8}
& $%
\begin{array}{c}
z \\ 
x \\ 
y%
\end{array}%
$ & $%
\begin{array}{c}
x \\ 
y \\ 
z%
\end{array}%
$ & $%
\begin{array}{c}
x \\ 
z \\ 
y%
\end{array}%
$ & $\cdots $ & $%
\begin{array}{c}
x \\ 
z \\ 
y%
\end{array}%
$ & $%
\begin{array}{c}
z \\ 
y \\ 
x%
\end{array}%
$ & $%
\begin{array}{c}
x \\ 
y \\ 
z%
\end{array}%
$ \\ \cline{2-8}
\end{tabular}%
\medskip \newline
we also have $x$ not chosen or 1 manipulates from $L1^{\ast }$ to $u4$.
Next, at $(n-1)$-variant profile $u5$:\medskip

\begin{tabular}{c|c|c|c|c|c|c|c|}
\cline{2-8}
$u5$ & $1$ & $2$ & $3$ & $\cdots $ & $n-2$ & $n-1$ & $n$ \\ \cline{2-8}
& $%
\begin{array}{c}
z \\ 
x \\ 
y%
\end{array}%
$ & $%
\begin{array}{c}
x \\ 
y \\ 
z%
\end{array}%
$ & $%
\begin{array}{c}
x \\ 
z \\ 
y%
\end{array}%
$ & $\cdots $ & $%
\begin{array}{c}
x \\ 
z \\ 
y%
\end{array}%
$ & $%
\begin{array}{c}
y \\ 
z \\ 
x%
\end{array}%
$ & $%
\begin{array}{c}
x \\ 
y \\ 
z%
\end{array}%
$ \\ \cline{2-8}
\end{tabular}%
\medskip \newline
we have $g(u5)\neq x$ or $n-1$ manipulates from $u5$ to $u4$. Then $%
g(u3)\neq x$ or $n$ manipulates from $u5$ to $u3$. Hence $g(u3)=z$\medskip

\textbf{Result 3}.\medskip

At profile $L2$:\medskip

\begin{tabular}{c|c|c|c|c|c|c|c|}
\cline{2-8}
$L2$ & $1$ & $2$ & $3$ & $\cdots $ & $n-2$ & $n-1$ & $n$ \\ \cline{2-8}
& $%
\begin{array}{c}
z \\ 
x \\ 
y%
\end{array}%
$ & $%
\begin{array}{c}
x \\ 
z \\ 
y%
\end{array}%
$ & $%
\begin{array}{c}
x \\ 
z \\ 
y%
\end{array}%
$ & $\cdots $ & $%
\begin{array}{c}
x \\ 
z \\ 
y%
\end{array}%
$ & $%
\begin{array}{c}
y \\ 
x \\ 
z%
\end{array}%
$ & $%
\begin{array}{c}
x \\ 
z \\ 
y%
\end{array}%
$ \\ \cline{2-8}
\end{tabular}%
\medskip \newline
we have $g(L2)=x$ by assumption. Then at 2-variant profile $u6$:\medskip

\begin{tabular}{c|c|c|c|c|c|c|c|}
\cline{2-8}
$u6$ & $1$ & $2$ & $3$ & $\cdots $ & $n-2$ & $n-1$ & $n$ \\ \cline{2-8}
& $%
\begin{array}{c}
z \\ 
x \\ 
y%
\end{array}%
$ & $%
\begin{array}{c}
y \\ 
x \\ 
z%
\end{array}%
$ & $%
\begin{array}{c}
x \\ 
z \\ 
y%
\end{array}%
$ & $\cdots $ & $%
\begin{array}{c}
x \\ 
z \\ 
y%
\end{array}%
$ & $%
\begin{array}{c}
y \\ 
x \\ 
z%
\end{array}%
$ & $%
\begin{array}{c}
x \\ 
z \\ 
y%
\end{array}%
$ \\ \cline{2-8}
\end{tabular}%
\medskip \newline
$g(u6)\neq z$ or 2 manipulates from $u6$ to $L2$.\medskip

\textbf{Result 4}.\medskip

\begin{tabular}{c|c|c|c|c|c|c|c|}
\cline{2-8}
$u7$ & $1$ & $2$ & $3$ & $\cdots $ & $n-2$ & $n-1$ & $n$ \\ \cline{2-8}
& $%
\begin{array}{c}
z \\ 
x \\ 
y%
\end{array}%
$ & $%
\begin{array}{c}
y \\ 
x \\ 
z%
\end{array}%
$ & $%
\begin{array}{c}
x \\ 
z \\ 
y%
\end{array}%
$ & $\cdots $ & $%
\begin{array}{c}
x \\ 
z \\ 
y%
\end{array}%
$ & $%
\begin{array}{c}
x \\ 
z \\ 
y%
\end{array}%
$ & $%
\begin{array}{c}
x \\ 
z \\ 
y%
\end{array}%
$ \\ \cline{2-8}
\end{tabular}%
\medskip \newline
$g(u7)=z$ since $\{1,3,...,n\}$ is decisive for $z$ against $y$ on $NP^{\ast
}$. If $g(u6)=x$, then $n-1$ would manipulate from $u7$ to $u6$. Therefore, $%
g(u6)\neq x$. Hence $g(u6)=y$.\medskip

\textbf{Main Thread}\medskip

By Result 3 and Result 4, $y$ is chosen at $u6$:\medskip \newline
\begin{tabular}{c|c|c|c|c|c|c|c|}
\cline{2-8}
$u6$ & $1$ & $2$ & $3$ & $\cdots $ & $n-2$ & $n-1$ & $n$ \\ \cline{2-8}
& $%
\begin{array}{c}
z \\ 
x \\ 
y%
\end{array}%
$ & $%
\begin{array}{c}
y \\ 
x \\ 
z%
\end{array}%
$ & $%
\begin{array}{c}
x \\ 
z \\ 
y%
\end{array}%
$ & $\cdots $ & $%
\begin{array}{c}
x \\ 
z \\ 
y%
\end{array}%
$ & $%
\begin{array}{c}
y \\ 
x \\ 
z%
\end{array}%
$ & $%
\begin{array}{c}
x \\ 
z \\ 
y%
\end{array}%
$ \\ \cline{2-8}
\end{tabular}%
\medskip \newline
So $y$ is also chosen at $(n-1)$-variant $u8$:\medskip

\begin{tabular}{c|c|c|c|c|c|c|c|}
\cline{2-8}
$u8$ & $1$ & $2$ & $3$ & $\cdots $ & $n-2$ & $n-1$ & $n$ \\ \cline{2-8}
& $%
\begin{array}{c}
z \\ 
x \\ 
y%
\end{array}%
$ & $%
\begin{array}{c}
y \\ 
x \\ 
z%
\end{array}%
$ & $%
\begin{array}{c}
x \\ 
z \\ 
y%
\end{array}%
$ & $\cdots $ & $%
\begin{array}{c}
x \\ 
z \\ 
y%
\end{array}%
$ & $%
\begin{array}{c}
y \\ 
z \\ 
x%
\end{array}%
$ & $%
\begin{array}{c}
x \\ 
z \\ 
y%
\end{array}%
$ \\ \cline{2-8}
\end{tabular}%
\medskip \newline
or $n-1$ will manipulate from $u8$ to $u6$. Then again, $y$ is chosen at $n$%
-variant profile $u9$:\medskip

\begin{tabular}{c|c|c|c|c|c|c|c|}
\cline{2-8}
$u9$ & $1$ & $2$ & $3$ & $\cdots $ & $n-2$ & $n-1$ & $n$ \\ \cline{2-8}
& $%
\begin{array}{c}
z \\ 
x \\ 
y%
\end{array}%
$ & $%
\begin{array}{c}
y \\ 
x \\ 
z%
\end{array}%
$ & $%
\begin{array}{c}
x \\ 
z \\ 
y%
\end{array}%
$ & $\cdots $ & $%
\begin{array}{c}
x \\ 
z \\ 
y%
\end{array}%
$ & $%
\begin{array}{c}
y \\ 
z \\ 
x%
\end{array}%
$ & $%
\begin{array}{c}
z \\ 
x \\ 
y%
\end{array}%
$ \\ \cline{2-8}
\end{tabular}%
\medskip \newline
or $n$ manipulates from $u8$ to $u9$.\medskip

\qquad But by Result 1 and Result 2, $z$ is chosen at $u3$:\medskip

\begin{tabular}{c|c|c|c|c|c|c|c|}
\cline{2-8}
$u3$ & $1$ & $2$ & $3$ & $\cdots $ & $n-2$ & $n-1$ & $n$ \\ \cline{2-8}
& $%
\begin{array}{c}
z \\ 
x \\ 
y%
\end{array}%
$ & $%
\begin{array}{c}
x \\ 
y \\ 
z%
\end{array}%
$ & $%
\begin{array}{c}
x \\ 
z \\ 
y%
\end{array}%
$ & $\cdots $ & $%
\begin{array}{c}
x \\ 
z \\ 
y%
\end{array}%
$ & $%
\begin{array}{c}
y \\ 
z \\ 
x%
\end{array}%
$ & $%
\begin{array}{c}
z \\ 
x \\ 
y%
\end{array}%
$ \\ \cline{2-8}
\end{tabular}%
\medskip \newline
So 2 would manipulate from $u3$ to $u9$, showing that $g$ violates
strategy-proofness.\medskip

\qquad Possibility b. $n-1$ is excluded and $\{1,2,...,n-2\}$ is decisive
for $z$ against $y$ for $g$ and so $\{1,2,...,n-2\}$ is decisive for $z$
against $y$ for $g$ on $NP^{\ast }$.\medskip

\qquad This is covered by $L2$.A, except when $T$ is the singleton $\{1\}$,
so here we assume $\{1\}$ is decisive for $y$ against $z$ and also $z$
against $y$ on $NP^{\ast }$.\medskip

We are assuming $x$ is chosen at $L2$:\medskip

\begin{tabular}{c|c|c|c|c|c|c|c|}
\cline{2-8}
$L2$ & $1$ & $2$ & $3$ & $\cdots $ & $n-2$ & $n-1$ & $n$ \\ \cline{2-8}
& $%
\begin{array}{c}
z \\ 
x \\ 
y%
\end{array}%
$ & $%
\begin{array}{c}
x \\ 
z \\ 
y%
\end{array}%
$ & $%
\begin{array}{c}
x \\ 
z \\ 
y%
\end{array}%
$ & $\cdots $ & $%
\begin{array}{c}
x \\ 
z \\ 
y%
\end{array}%
$ & $%
\begin{array}{c}
y \\ 
x \\ 
z%
\end{array}%
$ & $%
\begin{array}{c}
x \\ 
z \\ 
y%
\end{array}%
$ \\ \cline{2-8}
\end{tabular}%
\medskip \newline
So at profile $u1$, where $y$ and $z$ are interchanged below $x$ for $%
2,...,n-2$, strategy-proofness implies $x$ is still chosen:\medskip

\begin{tabular}{c|c|c|c|c|c|c|c|}
\cline{2-8}
$u1$ & $1$ & $2$ & $3$ & $\cdots $ & $n-2$ & $n-1$ & $n$ \\ \cline{2-8}
& $%
\begin{array}{c}
z \\ 
x \\ 
y%
\end{array}%
$ & $%
\begin{array}{c}
x \\ 
y \\ 
z%
\end{array}%
$ & $%
\begin{array}{c}
x \\ 
y \\ 
z%
\end{array}%
$ & $\cdots $ & $%
\begin{array}{c}
x \\ 
y \\ 
z%
\end{array}%
$ & $%
\begin{array}{c}
y \\ 
x \\ 
z%
\end{array}%
$ & $%
\begin{array}{c}
x \\ 
z \\ 
y%
\end{array}%
$ \\ \cline{2-8}
\end{tabular}%
\medskip \newline
This implies that at 1-variant profile $u2$:\medskip

\begin{tabular}{c|c|c|c|c|c|c|c|}
\cline{2-8}
$u2$ & $1$ & $2$ & $3$ & $\cdots $ & $n-2$ & $n-1$ & $n$ \\ \cline{2-8}
& $%
\begin{array}{c}
z \\ 
y \\ 
x%
\end{array}%
$ & $%
\begin{array}{c}
x \\ 
y \\ 
z%
\end{array}%
$ & $%
\begin{array}{c}
x \\ 
y \\ 
z%
\end{array}%
$ & $\cdots $ & $%
\begin{array}{c}
x \\ 
y \\ 
z%
\end{array}%
$ & $%
\begin{array}{c}
y \\ 
x \\ 
z%
\end{array}%
$ & $%
\begin{array}{c}
x \\ 
z \\ 
y%
\end{array}%
$ \\ \cline{2-8}
\end{tabular}%
\medskip \newline
$g(u2)\neq z$ or 1 will manipulate from $u1$ to $u2$.\medskip

\qquad Now $u2$ is also an $(n-1)$-variant of $u3$:\medskip

\begin{tabular}{c|c|c|c|c|c|c|c|}
\cline{2-8}
$u3$ & $1$ & $2$ & $3$ & $\cdots $ & $n-2$ & $n-1$ & $n$ \\ \cline{2-8}
& $%
\begin{array}{c}
z \\ 
y \\ 
x%
\end{array}%
$ & $%
\begin{array}{c}
x \\ 
y \\ 
z%
\end{array}%
$ & $%
\begin{array}{c}
x \\ 
y \\ 
z%
\end{array}%
$ & $\cdots $ & $%
\begin{array}{c}
x \\ 
y \\ 
z%
\end{array}%
$ & $%
\begin{array}{c}
x \\ 
y \\ 
z%
\end{array}%
$ & $%
\begin{array}{c}
x \\ 
z \\ 
y%
\end{array}%
$ \\ \cline{2-8}
\end{tabular}%
\medskip \newline
If $g(u3)=z$, $n-1$ will manipulate from $u3$ to $u2$. So $g(u3)\neq z$%
.\medskip

\qquad But then consider $n$ variant profile $u4$:\medskip

\begin{tabular}{c|c|c|c|c|c|c|c|}
\cline{2-8}
$u4$ & $1$ & $2$ & $3$ & $\cdots $ & $n-2$ & $n-1$ & $n$ \\ \cline{2-8}
& $%
\begin{array}{c}
z \\ 
y \\ 
x%
\end{array}%
$ & $%
\begin{array}{c}
x \\ 
y \\ 
z%
\end{array}%
$ & $%
\begin{array}{c}
x \\ 
y \\ 
z%
\end{array}%
$ & $\cdots $ & $%
\begin{array}{c}
x \\ 
y \\ 
z%
\end{array}%
$ & $%
\begin{array}{c}
x \\ 
y \\ 
z%
\end{array}%
$ & $%
\begin{array}{c}
x \\ 
y \\ 
z%
\end{array}%
$ \\ \cline{2-8}
\end{tabular}%
\medskip \newline
where $g(u4)=z$ since $\{1\}$ is decisive for $z$ against $y$ on $NP^{\ast }$%
. But then $n$ manipulates from $u4$ to $u3$ showing that $g$ violates
strategy-proofness.\medskip

\qquad \textbf{Subcase }$L2$B.2: A singleton from $\{2,...,n-2\}$, say $%
\{2\} $, is a minimal decisive set for $y$ against $z$ for $g|NP^{\ast }$.
This is covered by Case $L2$.A since we did not have to assume that T had
more than two members. Since that proof made no use of coalitions decisive
for $z$ against $y$, we also don't need to treat such issues here.\medskip

\qquad \textbf{Case }$L2$.C: This time we assume $\{n-1,n\}$ is a decisive
coalition both ways for $g|NP^{\ast }$. We first establish two intermediate
results.\medskip

\textbf{Result 1}\medskip

\qquad At profile $u1$ in $NP^{\ast }$:\medskip

\begin{tabular}{c|c|c|c|c|c|c|c|}
\cline{2-8}
$u1$ & $1$ & $2$ & $3$ & $\cdots $ & $n-2$ & $n-1$ & $n$ \\ \cline{2-8}
& $%
\begin{array}{c}
z \\ 
y \\ 
x%
\end{array}%
$ & $%
\begin{array}{c}
x \\ 
y \\ 
z%
\end{array}%
$ & $%
\begin{array}{c}
x \\ 
y \\ 
z%
\end{array}%
$ & $\cdots $ & $%
\begin{array}{c}
x \\ 
y \\ 
z%
\end{array}%
$ & $%
\begin{array}{c}
x \\ 
z \\ 
y%
\end{array}%
$ & $%
\begin{array}{c}
x \\ 
z \\ 
y%
\end{array}%
$ \\ \cline{2-8}
\end{tabular}%
\medskip \newline
we have $g(u1)=z$ by the decisiveness of $\{n-1,n\}$ on $NP^{\ast }$.
Therefore, $g(u2)\neq x$ at $(n-1)$-variant profile $u2$:\medskip

\begin{tabular}{c|c|c|c|c|c|c|c|}
\cline{2-8}
$u2$ & $1$ & $2$ & $3$ & $\cdots $ & $n-2$ & $n-1$ & $n$ \\ \cline{2-8}
& $%
\begin{array}{c}
z \\ 
y \\ 
x%
\end{array}%
$ & $%
\begin{array}{c}
x \\ 
y \\ 
z%
\end{array}%
$ & $%
\begin{array}{c}
x \\ 
y \\ 
z%
\end{array}%
$ & $\cdots $ & $%
\begin{array}{c}
x \\ 
y \\ 
z%
\end{array}%
$ & $%
\begin{array}{c}
y \\ 
x \\ 
z%
\end{array}%
$ & $%
\begin{array}{c}
x \\ 
z \\ 
y%
\end{array}%
$ \\ \cline{2-8}
\end{tabular}%
\medskip \newline
or $n-1$ manipulates from $u1$ to $u2$. We next show $g(u2)\neq z$.\medskip

\qquad We are assuming $x$ is chosen at $L2$:\medskip

\begin{tabular}{c|c|c|c|c|c|c|c|}
\cline{2-8}
$L2$ & $1$ & $2$ & $3$ & $\cdots $ & $n-2$ & $n-1$ & $n$ \\ \cline{2-8}
& $%
\begin{array}{c}
z \\ 
x \\ 
y%
\end{array}%
$ & $%
\begin{array}{c}
x \\ 
z \\ 
y%
\end{array}%
$ & $%
\begin{array}{c}
x \\ 
z \\ 
y%
\end{array}%
$ & $\cdots $ & $%
\begin{array}{c}
x \\ 
z \\ 
y%
\end{array}%
$ & $%
\begin{array}{c}
y \\ 
x \\ 
z%
\end{array}%
$ & $%
\begin{array}{c}
x \\ 
z \\ 
y%
\end{array}%
$ \\ \cline{2-8}
\end{tabular}%
\medskip \newline
Clearly a sequence of switches of $y$ and $z$ for individuals $2,...,n-2$
will leave $x$ still chosen at profile $u3$:\medskip

\begin{tabular}{c|c|c|c|c|c|c|c|}
\cline{2-8}
$u3$ & $1$ & $2$ & $3$ & $\cdots $ & $n-2$ & $n-1$ & $n$ \\ \cline{2-8}
& $%
\begin{array}{c}
z \\ 
x \\ 
y%
\end{array}%
$ & $%
\begin{array}{c}
x \\ 
y \\ 
z%
\end{array}%
$ & $%
\begin{array}{c}
x \\ 
y \\ 
z%
\end{array}%
$ & $\cdots $ & $%
\begin{array}{c}
x \\ 
y \\ 
z%
\end{array}%
$ & $%
\begin{array}{c}
y \\ 
x \\ 
z%
\end{array}%
$ & $%
\begin{array}{c}
x \\ 
z \\ 
y%
\end{array}%
$ \\ \cline{2-8}
\end{tabular}%
\medskip \newline
But then, if $g(u2)=z$, 1 would manipulate from $u3$ to $u2$. So $g(u2)\neq
z $. Combining, $g(u2)=y$. This, in turn, implies $g(u4)=y$ at $u4$:\medskip

\begin{tabular}{c|c|c|c|c|c|c|c|}
\cline{2-8}
$u4$ & $1$ & $2$ & $3$ & $\cdots $ & $n-2$ & $n-1$ & $n$ \\ \cline{2-8}
& $%
\begin{array}{c}
z \\ 
y \\ 
x%
\end{array}%
$ & $%
\begin{array}{c}
x \\ 
y \\ 
z%
\end{array}%
$ & $%
\begin{array}{c}
x \\ 
y \\ 
z%
\end{array}%
$ & $\cdots $ & $%
\begin{array}{c}
x \\ 
y \\ 
z%
\end{array}%
$ & $%
\begin{array}{c}
y \\ 
z \\ 
x%
\end{array}%
$ & $%
\begin{array}{c}
x \\ 
z \\ 
y%
\end{array}%
$ \\ \cline{2-8}
\end{tabular}%
\medskip \newline
or $n-1$ manipulates from $u4$ to $u2$.\medskip

\textbf{Result 2}\medskip

\qquad From earlier analysis, we know $x$ is not chosen at $L1$:\medskip

\begin{tabular}{c|c|c|c|c|c|c|c|}
\cline{2-8}
$L1$ & $1$ & $2$ & $3$ & $\cdots $ & $n-2$ & $n-1$ & $n$ \\ \cline{2-8}
& $%
\begin{array}{c}
x \\ 
y \\ 
z%
\end{array}%
$ & $%
\begin{array}{c}
x \\ 
y \\ 
z%
\end{array}%
$ & $%
\begin{array}{c}
x \\ 
y \\ 
z%
\end{array}%
$ & $\cdots $ & $%
\begin{array}{c}
x \\ 
y \\ 
z%
\end{array}%
$ & $%
\begin{array}{c}
z \\ 
y \\ 
x%
\end{array}%
$ & $%
\begin{array}{c}
x \\ 
y \\ 
z%
\end{array}%
$ \\ \cline{2-8}
\end{tabular}%
\medskip \newline
But then at $u5$:\medskip

\begin{tabular}{c|c|c|c|c|c|c|c|}
\cline{2-8}
$u5$ & $1$ & $2$ & $3$ & $\cdots $ & $n-2$ & $n-1$ & $n$ \\ \cline{2-8}
& $%
\begin{array}{c}
z \\ 
x \\ 
y%
\end{array}%
$ & $%
\begin{array}{c}
x \\ 
y \\ 
z%
\end{array}%
$ & $%
\begin{array}{c}
x \\ 
y \\ 
z%
\end{array}%
$ & $\cdots $ & $%
\begin{array}{c}
x \\ 
y \\ 
z%
\end{array}%
$ & $%
\begin{array}{c}
z \\ 
y \\ 
x%
\end{array}%
$ & $%
\begin{array}{c}
x \\ 
y \\ 
z%
\end{array}%
$ \\ \cline{2-8}
\end{tabular}%
\medskip \newline
$g(u5)\neq x$ or 1 manipulates from $L1$ to $u5$. But then at $(n-1)$%
-variant profile $u6$:\medskip

\begin{tabular}{c|c|c|c|c|c|c|c|}
\cline{2-8}
$u6$ & $1$ & $2$ & $3$ & $\cdots $ & $n-2$ & $n-1$ & $n$ \\ \cline{2-8}
& $%
\begin{array}{c}
z \\ 
x \\ 
y%
\end{array}%
$ & $%
\begin{array}{c}
x \\ 
y \\ 
z%
\end{array}%
$ & $%
\begin{array}{c}
x \\ 
y \\ 
z%
\end{array}%
$ & $\cdots $ & $%
\begin{array}{c}
x \\ 
y \\ 
z%
\end{array}%
$ & $%
\begin{array}{c}
y \\ 
z \\ 
x%
\end{array}%
$ & $%
\begin{array}{c}
x \\ 
y \\ 
z%
\end{array}%
$ \\ \cline{2-8}
\end{tabular}%
\medskip \newline
$g(u6)\neq x$ or $n-1$ would manipulate from $u6$ to $u5$. But then consider 
$n$-variant $u7$:\medskip

\begin{tabular}{c|c|c|c|c|c|c|c|}
\cline{2-8}
$u7$ & $1$ & $2$ & $3$ & $\cdots $ & $n-2$ & $n-1$ & $n$ \\ \cline{2-8}
& $%
\begin{array}{c}
z \\ 
x \\ 
y%
\end{array}%
$ & $%
\begin{array}{c}
x \\ 
y \\ 
z%
\end{array}%
$ & $%
\begin{array}{c}
x \\ 
y \\ 
z%
\end{array}%
$ & $\cdots $ & $%
\begin{array}{c}
x \\ 
y \\ 
z%
\end{array}%
$ & $%
\begin{array}{c}
y \\ 
z \\ 
x%
\end{array}%
$ & $%
\begin{array}{c}
x \\ 
z \\ 
y%
\end{array}%
$ \\ \cline{2-8}
\end{tabular}%
\medskip \newline
$g(u7)\neq x$ or $n$ would manipulate from $u6$ to $u7$. Note that $u7$ is
also an $(n-1)$-variant of $u3$ where $g(u3)=x$. If $g(u7)=y$, then $n-1$
would manipulate from $u3$ to $u7$. Therefore, $g(u7$) $\neq $ $y$.
Combining, $g(u7$) = $z$.\medskip

\textbf{Main Thread}\medskip

\qquad From Result 1, $y$ is chosen at $u4$:\medskip

\begin{tabular}{c|c|c|c|c|c|c|c|}
\cline{2-8}
$u4$ & $1$ & $2$ & $3$ & $\cdots $ & $n-2$ & $n-1$ & $n$ \\ \cline{2-8}
& $%
\begin{array}{c}
z \\ 
y \\ 
x%
\end{array}%
$ & $%
\begin{array}{c}
x \\ 
y \\ 
z%
\end{array}%
$ & $%
\begin{array}{c}
x \\ 
y \\ 
z%
\end{array}%
$ & $\cdots $ & $%
\begin{array}{c}
x \\ 
y \\ 
z%
\end{array}%
$ & $%
\begin{array}{c}
y \\ 
z \\ 
x%
\end{array}%
$ & $%
\begin{array}{c}
x \\ 
z \\ 
y%
\end{array}%
$ \\ \cline{2-8}
\end{tabular}%
\medskip \newline
and from Result 2, $z$ is chosen at $u7$:\medskip

\begin{tabular}{c|c|c|c|c|c|c|c|}
\cline{2-8}
$u7$ & $1$ & $2$ & $3$ & $\cdots $ & $n-2$ & $n-1$ & $n$ \\ \cline{2-8}
& $%
\begin{array}{c}
z \\ 
x \\ 
y%
\end{array}%
$ & $%
\begin{array}{c}
x \\ 
y \\ 
z%
\end{array}%
$ & $%
\begin{array}{c}
x \\ 
y \\ 
z%
\end{array}%
$ & $\cdots $ & $%
\begin{array}{c}
x \\ 
y \\ 
z%
\end{array}%
$ & $%
\begin{array}{c}
y \\ 
z \\ 
x%
\end{array}%
$ & $%
\begin{array}{c}
x \\ 
z \\ 
y%
\end{array}%
$ \\ \cline{2-8}
\end{tabular}%
\medskip \newline
Then 1 manipulates from $u4$ to $u7$, and $g$ violates
strategy-proofness.\medskip

Summarizing, we have established that $x$ is not selected at $L2$ or $L1$ if
the range of $g|NP^{\ast }$ is $\{y,z\}$.\medskip

\textbf{Section 5}. $L3$.\medskip

\textbf{Case }$L3$.A:\medskip

In this Case, we assume $\{1,...,k\}$ with $1<k\leq n-2$, is a minimal set
decisive for $y$ against $z$. We want to show $x$ is not chosen at $L3$. \
Assume to the contrary that $g(L3)=x$.\medskip

\begin{tabular}{c|c|c|c|c|c|c|c|c|c|c|c|}
\cline{2-12}
$L3$ & $1$ & $2$ & $3$ & $\cdots $ & $k-1$ & $k$ & $k+1$ & $\cdots $ & $n-2$
& $n-1$ & $n$ \\ \cline{2-12}
& $%
\begin{array}{c}
x \\ 
z \\ 
y%
\end{array}%
$ & $%
\begin{array}{c}
x \\ 
z \\ 
y%
\end{array}%
$ & $%
\begin{array}{c}
x \\ 
z \\ 
y%
\end{array}%
$ & $\cdots $ & $%
\begin{array}{c}
x \\ 
z \\ 
y%
\end{array}%
$ & $%
\begin{array}{c}
x \\ 
z \\ 
y%
\end{array}%
$ & $%
\begin{array}{c}
x \\ 
z \\ 
y%
\end{array}%
$ & $\cdots $ & $%
\begin{array}{c}
x \\ 
z \\ 
y%
\end{array}%
$ & $%
\begin{array}{c}
y \\ 
x \\ 
z%
\end{array}%
$ & $%
\begin{array}{c}
z \\ 
x \\ 
y%
\end{array}%
$ \\ \cline{2-12}
\end{tabular}%
\medskip \newline
But $x$ is chosen at $L3$ if and only if $x$ is chosen at $L3$*:\medskip

\begin{tabular}{c|c|c|c|c|c|c|c|c|c|c|c|}
\cline{2-12}
$L3$* & $1$ & $2$ & $3$ & $\cdots $ & $k-1$ & $k$ & $k+1$ & $\cdots $ & $n-2$
& $n-1$ & $n$ \\ \cline{2-12}
& $%
\begin{array}{c}
x \\ 
z \\ 
y%
\end{array}%
$ & $%
\begin{array}{c}
x \\ 
z \\ 
y%
\end{array}%
$ & $%
\begin{array}{c}
x \\ 
y \\ 
z%
\end{array}%
$ & $\cdots $ & $%
\begin{array}{c}
x \\ 
y \\ 
z%
\end{array}%
$ & $%
\begin{array}{c}
x \\ 
y \\ 
z%
\end{array}%
$ & $%
\begin{array}{c}
x \\ 
z \\ 
y%
\end{array}%
$ & $\cdots $ & $%
\begin{array}{c}
x \\ 
z \\ 
y%
\end{array}%
$ & $%
\begin{array}{c}
y \\ 
x \\ 
z%
\end{array}%
$ & $%
\begin{array}{c}
z \\ 
x \\ 
y%
\end{array}%
$ \\ \cline{2-12}
\end{tabular}%
\medskip \newline
with $y$ and $z$ switched for individuals 3 to $k$. Next observe that at
profile $u1$ in $NP^{\ast }$:\medskip

\begin{tabular}{c|c|c|c|c|c|c|c|c|c|c|c|}
\cline{2-12}
$u1$ & $1$ & $2$ & $3$ & $\cdots $ & $k-1$ & $k$ & $k+1$ & $\cdots $ & $n-2$
& $n-1$ & $n$ \\ \cline{2-12}
& $%
\begin{array}{c}
x \\ 
y \\ 
z%
\end{array}%
$ & $%
\begin{array}{c}
x \\ 
z \\ 
y%
\end{array}%
$ & $%
\begin{array}{c}
x \\ 
y \\ 
z%
\end{array}%
$ & $\cdots $ & $%
\begin{array}{c}
x \\ 
y \\ 
z%
\end{array}%
$ & $%
\begin{array}{c}
x \\ 
y \\ 
z%
\end{array}%
$ & $%
\begin{array}{c}
x \\ 
z \\ 
y%
\end{array}%
$ & $\cdots $ & $%
\begin{array}{c}
x \\ 
z \\ 
y%
\end{array}%
$ & $%
\begin{array}{c}
z \\ 
y \\ 
x%
\end{array}%
$ & $%
\begin{array}{c}
z \\ 
y \\ 
x%
\end{array}%
$ \\ \cline{2-12}
\end{tabular}%
\medskip \newline
we have $g(u1)=z$ since $\{1,...,k\}$ is \textit{minimal }decisive for $y$
against $z$ on $NP^{\ast }$. Therefore at $n$-variant $u2$:\medskip

\begin{tabular}{c|c|c|c|c|c|c|c|c|c|c|c|}
\cline{2-12}
$u2$ & $1$ & $2$ & $3$ & $\cdots $ & $k-1$ & $k$ & $k+1$ & $\cdots $ & $n-2$
& $n-1$ & $n$ \\ \cline{2-12}
& $%
\begin{array}{c}
x \\ 
y \\ 
z%
\end{array}%
$ & $%
\begin{array}{c}
x \\ 
z \\ 
y%
\end{array}%
$ & $%
\begin{array}{c}
x \\ 
y \\ 
z%
\end{array}%
$ & $\cdots $ & $%
\begin{array}{c}
x \\ 
y \\ 
z%
\end{array}%
$ & $%
\begin{array}{c}
x \\ 
y \\ 
z%
\end{array}%
$ & $%
\begin{array}{c}
x \\ 
z \\ 
y%
\end{array}%
$ & $\cdots $ & $%
\begin{array}{c}
x \\ 
z \\ 
y%
\end{array}%
$ & $%
\begin{array}{c}
z \\ 
y \\ 
x%
\end{array}%
$ & $%
\begin{array}{c}
z \\ 
x \\ 
y%
\end{array}%
$ \\ \cline{2-12}
\end{tabular}%
\medskip \newline
we also have $g(u2)=z$ or $n$ would manipulate from $u2$ to $u1$.\medskip

\qquad At profile $u3$ in $NP^{\ast }$:\medskip

\begin{tabular}{c|c|c|c|c|c|c|c|c|c|c|c|}
\cline{2-12}
$u3$ & $1$ & $2$ & $3$ & $\cdots $ & $k-1$ & $k$ & $k+1$ & $\cdots $ & $n-2$
& $n-1$ & $n$ \\ \cline{2-12}
& $%
\begin{array}{c}
x \\ 
y \\ 
z%
\end{array}%
$ & $%
\begin{array}{c}
y \\ 
x \\ 
z%
\end{array}%
$ & $%
\begin{array}{c}
x \\ 
y \\ 
z%
\end{array}%
$ & $\cdots $ & $%
\begin{array}{c}
x \\ 
y \\ 
z%
\end{array}%
$ & $%
\begin{array}{c}
x \\ 
y \\ 
z%
\end{array}%
$ & $%
\begin{array}{c}
x \\ 
z \\ 
y%
\end{array}%
$ & $\cdots $ & $%
\begin{array}{c}
x \\ 
z \\ 
y%
\end{array}%
$ & $%
\begin{array}{c}
z \\ 
x \\ 
y%
\end{array}%
$ & $%
\begin{array}{c}
z \\ 
x \\ 
y%
\end{array}%
$ \\ \cline{2-12}
\end{tabular}%
\medskip \newline
we get $g(u3)=y$ since $\{1,...,k\}$ is decisive for $y$ against $z$ on $%
NP^{\ast }$. Therefore at $(n-1)$-variant $u4$:\medskip

\begin{tabular}{c|c|c|c|c|c|c|c|c|c|c|c|}
\cline{2-12}
$u4$ & $1$ & $2$ & $3$ & $\cdots $ & $k-1$ & $k$ & $k+1$ & $\cdots $ & $n-2$
& $n-1$ & $n$ \\ \cline{2-12}
& $%
\begin{array}{c}
x \\ 
y \\ 
z%
\end{array}%
$ & $%
\begin{array}{c}
y \\ 
x \\ 
z%
\end{array}%
$ & $%
\begin{array}{c}
x \\ 
y \\ 
z%
\end{array}%
$ & $\cdots $ & $%
\begin{array}{c}
x \\ 
y \\ 
z%
\end{array}%
$ & $%
\begin{array}{c}
x \\ 
y \\ 
z%
\end{array}%
$ & $%
\begin{array}{c}
x \\ 
z \\ 
y%
\end{array}%
$ & $\cdots $ & $%
\begin{array}{c}
x \\ 
z \\ 
y%
\end{array}%
$ & $%
\begin{array}{c}
z \\ 
y \\ 
x%
\end{array}%
$ & $%
\begin{array}{c}
z \\ 
x \\ 
y%
\end{array}%
$ \\ \cline{2-12}
\end{tabular}%
\medskip \newline
we also get $g(u4)=y$ or $n-1$ manipulates from $u3$ to $u4$.\medskip

\qquad Proceeding, at $u5$, a 2-variant of both $u2$ and $u4$:\medskip

\begin{tabular}{c|c|c|c|c|c|c|c|c|c|c|c|}
\cline{2-12}
$u5$ & $1$ & $2$ & $3$ & $\cdots $ & $k-1$ & $k$ & $k+1$ & $\cdots $ & $n-2$
& $n-1$ & $n$ \\ \cline{2-12}
& $%
\begin{array}{c}
x \\ 
y \\ 
z%
\end{array}%
$ & $%
\begin{array}{c}
x \\ 
y \\ 
z%
\end{array}%
$ & $%
\begin{array}{c}
x \\ 
y \\ 
z%
\end{array}%
$ & $\cdots $ & $%
\begin{array}{c}
x \\ 
y \\ 
z%
\end{array}%
$ & $%
\begin{array}{c}
x \\ 
y \\ 
z%
\end{array}%
$ & $%
\begin{array}{c}
x \\ 
z \\ 
y%
\end{array}%
$ & $\cdots $ & $%
\begin{array}{c}
x \\ 
z \\ 
y%
\end{array}%
$ & $%
\begin{array}{c}
z \\ 
y \\ 
x%
\end{array}%
$ & $%
\begin{array}{c}
z \\ 
x \\ 
y%
\end{array}%
$ \\ \cline{2-12}
\end{tabular}%
\medskip \newline
we have $g(u5)\neq x$ or 2 manipulates from $u2$ to $u5$ and $g(u5)\neq z$
or 2 manipulates from $u5$ to $u4$. Therefore $g(u5)=y$. But then at $(n-1)$%
-variant $u6$:\medskip

\begin{tabular}{c|c|c|c|c|c|c|c|c|c|c|c|}
\cline{2-12}
$u6$ & $1$ & $2$ & $3$ & $\cdots $ & $k-1$ & $k$ & $k+1$ & $\cdots $ & $n-2$
& $n-1$ & $n$ \\ \cline{2-12}
& $%
\begin{array}{c}
x \\ 
y \\ 
z%
\end{array}%
$ & $%
\begin{array}{c}
x \\ 
y \\ 
z%
\end{array}%
$ & $%
\begin{array}{c}
x \\ 
y \\ 
z%
\end{array}%
$ & $\cdots $ & $%
\begin{array}{c}
x \\ 
y \\ 
z%
\end{array}%
$ & $%
\begin{array}{c}
x \\ 
y \\ 
z%
\end{array}%
$ & $%
\begin{array}{c}
x \\ 
z \\ 
y%
\end{array}%
$ & $\cdots $ & $%
\begin{array}{c}
x \\ 
z \\ 
y%
\end{array}%
$ & $%
\begin{array}{c}
y \\ 
x \\ 
z%
\end{array}%
$ & $%
\begin{array}{c}
z \\ 
x \\ 
y%
\end{array}%
$ \\ \cline{2-12}
\end{tabular}%
\medskip \newline
we get $g(u6)=y$ or $n-1$ manipulates from $u6$ to $u5$. That in turn
implies $g(u7)\neq x$ at 2-variant $u7$:\medskip

\begin{tabular}{c|c|c|c|c|c|c|c|c|c|c|c|}
\cline{2-12}
$u7$ & $1$ & $2$ & $3$ & $\cdots $ & $k-1$ & $k$ & $k+1$ & $\cdots $ & $n-2$
& $n-1$ & $n$ \\ \cline{2-12}
& $%
\begin{array}{c}
x \\ 
y \\ 
z%
\end{array}%
$ & $%
\begin{array}{c}
x \\ 
z \\ 
y%
\end{array}%
$ & $%
\begin{array}{c}
x \\ 
y \\ 
z%
\end{array}%
$ & $\cdots $ & $%
\begin{array}{c}
x \\ 
y \\ 
z%
\end{array}%
$ & $%
\begin{array}{c}
x \\ 
y \\ 
z%
\end{array}%
$ & $%
\begin{array}{c}
x \\ 
z \\ 
y%
\end{array}%
$ & $\cdots $ & $%
\begin{array}{c}
x \\ 
z \\ 
y%
\end{array}%
$ & $%
\begin{array}{c}
y \\ 
x \\ 
z%
\end{array}%
$ & $%
\begin{array}{c}
z \\ 
x \\ 
y%
\end{array}%
$ \\ \cline{2-12}
\end{tabular}%
\medskip \newline
or 2 manipulates from $u6$ to $u7$. And then $g(L3^{\ast })\neq x$ or 1
manipulates from $u7$ to $L3$*:\medskip

\begin{tabular}{c|c|c|c|c|c|c|c|c|c|c|c|}
\cline{2-12}
$L3^{\ast }$ & $1$ & $2$ & $3$ & $\cdots $ & $k-1$ & $k$ & $k+1$ & $\cdots $
& $n-2$ & $n-1$ & $n$ \\ \cline{2-12}
& $%
\begin{array}{c}
x \\ 
z \\ 
y%
\end{array}%
$ & $%
\begin{array}{c}
x \\ 
z \\ 
y%
\end{array}%
$ & $%
\begin{array}{c}
x \\ 
y \\ 
z%
\end{array}%
$ & $\cdots $ & $%
\begin{array}{c}
x \\ 
y \\ 
z%
\end{array}%
$ & $%
\begin{array}{c}
x \\ 
y \\ 
z%
\end{array}%
$ & $%
\begin{array}{c}
x \\ 
z \\ 
y%
\end{array}%
$ & $\cdots $ & $%
\begin{array}{c}
x \\ 
z \\ 
y%
\end{array}%
$ & $%
\begin{array}{c}
y \\ 
x \\ 
z%
\end{array}%
$ & $%
\begin{array}{c}
z \\ 
x \\ 
y%
\end{array}%
$ \\ \cline{2-12}
\end{tabular}%
\medskip \newline
\textbf{Case }$L3$.B. The analysis here parallels that for $L1$.B: $\{1\}$
is a (minimal) decisive coalition for $y$ against $z$ for $g|NP^{\ast }$.
There are many possible coalitions $C$ decisive for $z$ against $y$. These
must include 1, but won't have to be minimal, so we choose them as large as
possible (though they still have to exclude at least one individual) since a
small coalition being decisive implies supersets also decisive. We
distinguish between two possibilities:\medskip

\qquad 1. The individual excluded for $g^{\ast }$ is in $\{1,...,n-2\}$, say 
$n-2$, so $\{1,...,n-3,n-1,n\}$ is decisive for $z$ against $y$ for $%
g|NP^{\ast }$.\medskip

\qquad 2. The individual excluded for $g^{\ast }$ is $n-1$, so $%
\{1,...,n-2\} $ is decisive for $z$ against $y$ for $g|NP^{\ast }$. If the
minimal decisive coalition for $z$ against $y$ in $\{1,...,n-2\}$ contains
an alternative other than 1, we are in the situation already covered in Case
A. So we may assume that $\{1\}$ is decisive for $z$ against $y$ as well as
for $y$ against $z$ for $g|NP^{\ast }$. But then also $\{1,...,n-3,n-1,n\}$
is decisive for $z$ against $y$ for $g|NP^{\ast }$, and we are back to the
first possibility.\medskip

\qquad So we only need to treat the case where coalition $\{1\}$ is decisive
for $y$ against $z$ for $g|NP^{\ast }$ and $\{1,...,n-3,n-1,n\}$ is decisive
for $z$ against $y$ for $g|NP^{\ast }$.\medskip

\textbf{Result 1}.\medskip

At profile $L3$,\medskip

\begin{tabular}{c|c|c|c|c|c|c|c|c|}
\cline{2-9}
$L3$ & $1$ & $2$ & $3$ & $\cdots $ & $n-3$ & $n-2$ & $n-1$ & $n$ \\ 
\cline{2-9}
& $%
\begin{array}{c}
x \\ 
z \\ 
y%
\end{array}%
$ & $%
\begin{array}{c}
x \\ 
z \\ 
y%
\end{array}%
$ & $%
\begin{array}{c}
x \\ 
z \\ 
y%
\end{array}%
$ & $\cdots $ & $%
\begin{array}{c}
x \\ 
z \\ 
y%
\end{array}%
$ & $%
\begin{array}{c}
x \\ 
z \\ 
y%
\end{array}%
$ & $%
\begin{array}{c}
y \\ 
x \\ 
z%
\end{array}%
$ & $%
\begin{array}{c}
z \\ 
x \\ 
y%
\end{array}%
$ \\ \cline{2-9}
\end{tabular}%
\medskip \newline
we have by assumption, $g(L3)=x$. Then at 1-variant profile $u1$:\medskip

\begin{tabular}{c|c|c|c|c|c|c|c|c|}
\cline{2-9}
$u1$ & $1$ & $2$ & $3$ & $\cdots $ & $n-3$ & $n-2$ & $n-1$ & $n$ \\ 
\cline{2-9}
& $%
\begin{array}{c}
x \\ 
y \\ 
z%
\end{array}%
$ & $%
\begin{array}{c}
x \\ 
z \\ 
y%
\end{array}%
$ & $%
\begin{array}{c}
x \\ 
z \\ 
y%
\end{array}%
$ & $\cdots $ & $%
\begin{array}{c}
x \\ 
z \\ 
y%
\end{array}%
$ & $%
\begin{array}{c}
x \\ 
z \\ 
y%
\end{array}%
$ & $%
\begin{array}{c}
y \\ 
x \\ 
z%
\end{array}%
$ & $%
\begin{array}{c}
z \\ 
x \\ 
y%
\end{array}%
$ \\ \cline{2-9}
\end{tabular}%
\medskip \newline
we also have $g(u1)=x$ or 1 manipulates from $u1$ to $L3$. Then at 2-variant
profile $u2$:\medskip

\begin{tabular}{c|c|c|c|c|c|c|c|c|}
\cline{2-9}
$u2$ & $1$ & $2$ & $3$ & $\cdots $ & $n-3$ & $n-2$ & $n-1$ & $n$ \\ 
\cline{2-9}
& $%
\begin{array}{c}
x \\ 
y \\ 
z%
\end{array}%
$ & $%
\begin{array}{c}
x \\ 
y \\ 
z%
\end{array}%
$ & $%
\begin{array}{c}
x \\ 
z \\ 
y%
\end{array}%
$ & $\cdots $ & $%
\begin{array}{c}
x \\ 
z \\ 
y%
\end{array}%
$ & $%
\begin{array}{c}
x \\ 
z \\ 
y%
\end{array}%
$ & $%
\begin{array}{c}
y \\ 
x \\ 
z%
\end{array}%
$ & $%
\begin{array}{c}
z \\ 
x \\ 
y%
\end{array}%
$ \\ \cline{2-9}
\end{tabular}%
\medskip \newline
we also have $g(u2)=x$ or 2 manipulates from $u2$ to $u1$. Then consider $%
(n-1)$-variant profile $u3$:\medskip

\begin{tabular}{c|c|c|c|c|c|c|c|c|}
\cline{2-9}
$u3$ & $1$ & $2$ & $3$ & $\cdots $ & $n-3$ & $n-2$ & $n-1$ & $n$ \\ 
\cline{2-9}
& $%
\begin{array}{c}
x \\ 
y \\ 
z%
\end{array}%
$ & $%
\begin{array}{c}
x \\ 
y \\ 
z%
\end{array}%
$ & $%
\begin{array}{c}
x \\ 
z \\ 
y%
\end{array}%
$ & $\cdots $ & $%
\begin{array}{c}
x \\ 
z \\ 
y%
\end{array}%
$ & $%
\begin{array}{c}
x \\ 
z \\ 
y%
\end{array}%
$ & $%
\begin{array}{c}
z \\ 
y \\ 
x%
\end{array}%
$ & $%
\begin{array}{c}
z \\ 
x \\ 
y%
\end{array}%
$ \\ \cline{2-9}
\end{tabular}%
\medskip \newline
We must have $g(u3)\neq y$ or $n-1$ would manipulate from $u2$ to $u3$%
.\medskip

\textbf{Result 2}.\medskip

At profile $u4$, an $n$-variant of $u3$:\medskip

\begin{tabular}{c|c|c|c|c|c|c|c|c|}
\cline{2-9}
$u4$ & $1$ & $2$ & $3$ & $\cdots $ & $n-3$ & $n-2$ & $n-1$ & $n$ \\ 
\cline{2-9}
& $%
\begin{array}{c}
x \\ 
y \\ 
z%
\end{array}%
$ & $%
\begin{array}{c}
x \\ 
y \\ 
z%
\end{array}%
$ & $%
\begin{array}{c}
x \\ 
z \\ 
y%
\end{array}%
$ & $\cdots $ & $%
\begin{array}{c}
x \\ 
z \\ 
y%
\end{array}%
$ & $%
\begin{array}{c}
x \\ 
z \\ 
y%
\end{array}%
$ & $%
\begin{array}{c}
z \\ 
y \\ 
x%
\end{array}%
$ & $%
\begin{array}{c}
z \\ 
y \\ 
x%
\end{array}%
$ \\ \cline{2-9}
\end{tabular}%
\medskip \newline
$g(u4)=y$ since $\{1\}$ is decisive for $y$ against $z$ on $NP^{\ast }$. But
then $g(u3)\neq z$, or $n$ manipulates from $u4$ to $u3$. Hence $g(u3)=x$%
.\medskip

\textbf{Result 3}.\medskip

At profile $u5$:\medskip

\begin{tabular}{c|c|c|c|c|c|c|c|c|}
\cline{2-9}
$u5$ & $1$ & $2$ & $3$ & $\cdots $ & $n-3$ & $n-2$ & $n-1$ & $n$ \\ 
\cline{2-9}
& $%
\begin{array}{c}
x \\ 
z \\ 
y%
\end{array}%
$ & $%
\begin{array}{c}
x \\ 
z \\ 
y%
\end{array}%
$ & $%
\begin{array}{c}
x \\ 
z \\ 
y%
\end{array}%
$ & $\cdots $ & $%
\begin{array}{c}
x \\ 
z \\ 
y%
\end{array}%
$ & $%
\begin{array}{c}
x \\ 
y \\ 
z%
\end{array}%
$ & $%
\begin{array}{c}
z \\ 
y \\ 
x%
\end{array}%
$ & $%
\begin{array}{c}
z \\ 
y \\ 
x%
\end{array}%
$ \\ \cline{2-9}
\end{tabular}%
\medskip \newline
we have $g(u5)=z$ since $\{1,...,n-3,n-1,n\}$ is decisive for $z$ against $y$
on $NP^{\ast }$. Then at $n$-variant profile $u6$:\medskip

\begin{tabular}{c|c|c|c|c|c|c|c|c|}
\cline{2-9}
$u6$ & $1$ & $2$ & $3$ & $\cdots $ & $n-3$ & $n-2$ & $n-1$ & $n$ \\ 
\cline{2-9}
& $%
\begin{array}{c}
x \\ 
z \\ 
y%
\end{array}%
$ & $%
\begin{array}{c}
x \\ 
z \\ 
y%
\end{array}%
$ & $%
\begin{array}{c}
x \\ 
z \\ 
y%
\end{array}%
$ & $\cdots $ & $%
\begin{array}{c}
x \\ 
z \\ 
y%
\end{array}%
$ & $%
\begin{array}{c}
x \\ 
y \\ 
z%
\end{array}%
$ & $%
\begin{array}{c}
z \\ 
y \\ 
x%
\end{array}%
$ & $%
\begin{array}{c}
z \\ 
x \\ 
y%
\end{array}%
$ \\ \cline{2-9}
\end{tabular}%
\medskip \newline
we get $g(u6)=z$ or $n$ manipulates from $u6$ to $u5$.\medskip

\textbf{Main Thread}\medskip

Combining Result 1 and Result 2, $g(u3)=x$ at $u3$:\medskip

\begin{tabular}{c|c|c|c|c|c|c|c|c|}
\cline{2-9}
$u3$ & $1$ & $2$ & $3$ & $\cdots $ & $n-3$ & $n-2$ & $n-1$ & $n$ \\ 
\cline{2-9}
& $%
\begin{array}{c}
x \\ 
y \\ 
z%
\end{array}%
$ & $%
\begin{array}{c}
x \\ 
y \\ 
z%
\end{array}%
$ & $%
\begin{array}{c}
x \\ 
z \\ 
y%
\end{array}%
$ & $\cdots $ & $%
\begin{array}{c}
x \\ 
z \\ 
y%
\end{array}%
$ & $%
\begin{array}{c}
x \\ 
z \\ 
y%
\end{array}%
$ & $%
\begin{array}{c}
z \\ 
y \\ 
x%
\end{array}%
$ & $%
\begin{array}{c}
z \\ 
x \\ 
y%
\end{array}%
$ \\ \cline{2-9}
\end{tabular}%
\medskip \newline
and also $x$ is chosen at new profile $u7$:\medskip

\begin{tabular}{c|c|c|c|c|c|c|c|c|}
\cline{2-9}
$u7$ & $1$ & $2$ & $3$ & $\cdots $ & $n-3$ & $n-2$ & $n-1$ & $n$ \\ 
\cline{2-9}
& $%
\begin{array}{c}
x \\ 
y \\ 
z%
\end{array}%
$ & $%
\begin{array}{c}
x \\ 
z \\ 
y%
\end{array}%
$ & $%
\begin{array}{c}
x \\ 
z \\ 
y%
\end{array}%
$ & $\cdots $ & $%
\begin{array}{c}
x \\ 
z \\ 
y%
\end{array}%
$ & $%
\begin{array}{c}
x \\ 
y \\ 
z%
\end{array}%
$ & $%
\begin{array}{c}
z \\ 
y \\ 
x%
\end{array}%
$ & $%
\begin{array}{c}
z \\ 
x \\ 
y%
\end{array}%
$ \\ \cline{2-9}
\end{tabular}%
\medskip \newline
where $y$ and $z$ are switched for individuals 2 and $n-2$. But then 1
manipulates from $u6$ to $u7$.\medskip

\textbf{Case }$L3$.C: For our final case, we assume $\{n-1,n\}$ is a minimal
set decisive for $y$ against $z$ and for $z$ against $y$ for $g|NP^{\ast }$.
We assume $x$ is chosen at $L3$:\medskip

\begin{tabular}{c|c|c|c|c|c|c|c|}
\cline{2-8}
$L3$ & $1$ & $2$ & $3$ & $\cdots $ & $n-2$ & $n-1$ & $n$ \\ \cline{2-8}
& $%
\begin{array}{c}
x \\ 
z \\ 
y%
\end{array}%
$ & $%
\begin{array}{c}
x \\ 
z \\ 
y%
\end{array}%
$ & $%
\begin{array}{c}
x \\ 
z \\ 
y%
\end{array}%
$ & $\cdots $ & $%
\begin{array}{c}
x \\ 
z \\ 
y%
\end{array}%
$ & $%
\begin{array}{c}
y \\ 
x \\ 
z%
\end{array}%
$ & $%
\begin{array}{c}
z \\ 
x \\ 
y%
\end{array}%
$ \\ \cline{2-8}
\end{tabular}%
\medskip \newline
and seek a contradiction.\medskip

\textbf{Result 1}.\medskip

\qquad We must also have $x$ chosen at $u1$:\medskip

\begin{tabular}{c|c|c|c|c|c|c|c|}
\cline{2-8}
$u1$ & $1$ & $2$ & $3$ & $\cdots $ & $n-2$ & $n-1$ & $n$ \\ \cline{2-8}
& $%
\begin{array}{c}
x \\ 
z \\ 
y%
\end{array}%
$ & $%
\begin{array}{c}
x \\ 
y \\ 
z%
\end{array}%
$ & $%
\begin{array}{c}
x \\ 
z \\ 
y%
\end{array}%
$ & $\cdots $ & $%
\begin{array}{c}
x \\ 
z \\ 
y%
\end{array}%
$ & $%
\begin{array}{c}
y \\ 
x \\ 
z%
\end{array}%
$ & $%
\begin{array}{c}
z \\ 
x \\ 
y%
\end{array}%
$ \\ \cline{2-8}
\end{tabular}%
\medskip \newline
or 2 manipulates from $u1$ to $L3$. That implies $g(u2)\neq y$ at 1-variant $%
u2$:\medskip

\begin{tabular}{c|c|c|c|c|c|c|c|}
\cline{2-8}
$u2$ & $1$ & $2$ & $3$ & $\cdots $ & $n-2$ & $n-1$ & $n$ \\ \cline{2-8}
& $%
\begin{array}{c}
z \\ 
x \\ 
y%
\end{array}%
$ & $%
\begin{array}{c}
x \\ 
y \\ 
z%
\end{array}%
$ & $%
\begin{array}{c}
x \\ 
z \\ 
y%
\end{array}%
$ & $\cdots $ & $%
\begin{array}{c}
x \\ 
z \\ 
y%
\end{array}%
$ & $%
\begin{array}{c}
y \\ 
x \\ 
z%
\end{array}%
$ & $%
\begin{array}{c}
z \\ 
x \\ 
y%
\end{array}%
$ \\ \cline{2-8}
\end{tabular}%
\medskip \newline
or 1 manipulates from $u2$ to $u1$. Then $g(u3)\neq y$ at $n$-variant $u3$%
:\medskip

\begin{tabular}{c|c|c|c|c|c|c|c|}
\cline{2-8}
$u3$ & $1$ & $2$ & $3$ & $\cdots $ & $n-2$ & $n-1$ & $n$ \\ \cline{2-8}
& $%
\begin{array}{c}
z \\ 
x \\ 
y%
\end{array}%
$ & $%
\begin{array}{c}
x \\ 
y \\ 
z%
\end{array}%
$ & $%
\begin{array}{c}
x \\ 
z \\ 
y%
\end{array}%
$ & $\cdots $ & $%
\begin{array}{c}
x \\ 
z \\ 
y%
\end{array}%
$ & $%
\begin{array}{c}
y \\ 
x \\ 
z%
\end{array}%
$ & $%
\begin{array}{c}
x \\ 
z \\ 
y%
\end{array}%
$ \\ \cline{2-8}
\end{tabular}%
\medskip \newline
or $n$ manipulates from $u3$ to $u2$.\medskip

\qquad From earlier analysis, $g(L2)\neq x$ at $L2$:\medskip

\begin{tabular}{c|c|c|c|c|c|c|c|}
\cline{2-8}
$L2$ & $1$ & $2$ & $3$ & $\cdots $ & $n-2$ & $n-1$ & $n$ \\ \cline{2-8}
& $%
\begin{array}{c}
z \\ 
x \\ 
y%
\end{array}%
$ & $%
\begin{array}{c}
x \\ 
z \\ 
y%
\end{array}%
$ & $%
\begin{array}{c}
x \\ 
z \\ 
y%
\end{array}%
$ & $\cdots $ & $%
\begin{array}{c}
x \\ 
z \\ 
y%
\end{array}%
$ & $%
\begin{array}{c}
y \\ 
x \\ 
z%
\end{array}%
$ & $%
\begin{array}{c}
x \\ 
z \\ 
y%
\end{array}%
$ \\ \cline{2-8}
\end{tabular}%
\medskip \newline
Therefore $g(u3)\neq x$ or 2 manipulates from $L2$ to $u3$. Combining, $%
g(u3)=z$. Then $g(u4)=z$ at $n$-variant $u4$:\medskip

\begin{tabular}{c|c|c|c|c|c|c|c|}
\cline{2-8}
$u4$ & $1$ & $2$ & $3$ & $\cdots $ & $n-2$ & $n-1$ & $n$ \\ \cline{2-8}
& $%
\begin{array}{c}
z \\ 
x \\ 
y%
\end{array}%
$ & $%
\begin{array}{c}
x \\ 
y \\ 
z%
\end{array}%
$ & $%
\begin{array}{c}
x \\ 
z \\ 
y%
\end{array}%
$ & $\cdots $ & $%
\begin{array}{c}
x \\ 
z \\ 
y%
\end{array}%
$ & $%
\begin{array}{c}
y \\ 
x \\ 
z%
\end{array}%
$ & $%
\begin{array}{c}
z \\ 
y \\ 
x%
\end{array}%
$ \\ \cline{2-8}
\end{tabular}%
\medskip \newline
or $n$ manipulates from $u4$ to $u3$.\medskip

\textbf{Result 2}.\medskip

\qquad From $g(u1$) = $x$, we have $g(u5)\neq z$ at $n$-variant $u5$:\medskip

\begin{tabular}{c|c|c|c|c|c|c|c|}
\cline{2-8}
$u5$ & $1$ & $2$ & $3$ & $\cdots $ & $n-2$ & $n-1$ & $n$ \\ \cline{2-8}
& $%
\begin{array}{c}
x \\ 
z \\ 
y%
\end{array}%
$ & $%
\begin{array}{c}
x \\ 
y \\ 
z%
\end{array}%
$ & $%
\begin{array}{c}
x \\ 
z \\ 
y%
\end{array}%
$ & $\cdots $ & $%
\begin{array}{c}
x \\ 
z \\ 
y%
\end{array}%
$ & $%
\begin{array}{c}
y \\ 
x \\ 
z%
\end{array}%
$ & $%
\begin{array}{c}
z \\ 
y \\ 
x%
\end{array}%
$ \\ \cline{2-8}
\end{tabular}%
\medskip \newline
or $n$ manipulates from $u1$ to $u5$. We get more detail about $g(u5)$ by
observing that $y$ is chosen at $u6$ in $NP^{\ast }$:\medskip

\begin{tabular}{c|c|c|c|c|c|c|c|}
\cline{2-8}
$u6$ & $1$ & $2$ & $3$ & $\cdots $ & $n-2$ & $n-1$ & $n$ \\ \cline{2-8}
& $%
\begin{array}{c}
x \\ 
z \\ 
y%
\end{array}%
$ & $%
\begin{array}{c}
x \\ 
y \\ 
z%
\end{array}%
$ & $%
\begin{array}{c}
x \\ 
z \\ 
y%
\end{array}%
$ & $\cdots $ & $%
\begin{array}{c}
x \\ 
z \\ 
y%
\end{array}%
$ & $%
\begin{array}{c}
y \\ 
z \\ 
x%
\end{array}%
$ & $%
\begin{array}{c}
y \\ 
z \\ 
x%
\end{array}%
$ \\ \cline{2-8}
\end{tabular}%
\medskip \newline
because $\{n-1,n\}$ is decisive for $y$ against $z$ on $NP^{\ast }$.
Therefore at $(n-1)$-variant $u7$\medskip

\begin{tabular}{c|c|c|c|c|c|c|c|}
\cline{2-8}
$u7$ & $1$ & $2$ & $3$ & $\cdots $ & $n-2$ & $n-1$ & $n$ \\ \cline{2-8}
& $%
\begin{array}{c}
x \\ 
z \\ 
y%
\end{array}%
$ & $%
\begin{array}{c}
x \\ 
y \\ 
z%
\end{array}%
$ & $%
\begin{array}{c}
x \\ 
z \\ 
y%
\end{array}%
$ & $\cdots $ & $%
\begin{array}{c}
x \\ 
z \\ 
y%
\end{array}%
$ & $%
\begin{array}{c}
y \\ 
x \\ 
z%
\end{array}%
$ & $%
\begin{array}{c}
y \\ 
z \\ 
x%
\end{array}%
$ \\ \cline{2-8}
\end{tabular}%
\medskip \newline
we have $g(u7)=y$ or $n-1$ manipulates from $u7$ to $u6$. But then $%
g(u5)\neq x$ or $n$ manipulates from $u5$ to $u7$. Combining, $g(u5)=y$%
.\medskip

\textbf{Main Thread}\medskip

\qquad We now know that $y$ is chosen at $u5$\medskip

\begin{tabular}{c|c|c|c|c|c|c|c|}
\cline{2-8}
$u5$ & $1$ & $2$ & $3$ & $\cdots $ & $n-2$ & $n-1$ & $n$ \\ \cline{2-8}
& $%
\begin{array}{c}
x \\ 
z \\ 
y%
\end{array}%
$ & $%
\begin{array}{c}
x \\ 
y \\ 
z%
\end{array}%
$ & $%
\begin{array}{c}
x \\ 
z \\ 
y%
\end{array}%
$ & $\cdots $ & $%
\begin{array}{c}
x \\ 
z \\ 
y%
\end{array}%
$ & $%
\begin{array}{c}
y \\ 
x \\ 
z%
\end{array}%
$ & $%
\begin{array}{c}
z \\ 
y \\ 
x%
\end{array}%
$ \\ \cline{2-8}
\end{tabular}%
\medskip \newline
by Result 2 and $z$ is chosen at $u4$\medskip

\begin{tabular}{c|c|c|c|c|c|c|c|}
\cline{2-8}
$u4$ & $1$ & $2$ & $3$ & $\cdots $ & $n-2$ & $n-1$ & $n$ \\ \cline{2-8}
& $%
\begin{array}{c}
z \\ 
x \\ 
y%
\end{array}%
$ & $%
\begin{array}{c}
x \\ 
y \\ 
z%
\end{array}%
$ & $%
\begin{array}{c}
x \\ 
z \\ 
y%
\end{array}%
$ & $\cdots $ & $%
\begin{array}{c}
x \\ 
z \\ 
y%
\end{array}%
$ & $%
\begin{array}{c}
y \\ 
x \\ 
z%
\end{array}%
$ & $%
\begin{array}{c}
z \\ 
y \\ 
x%
\end{array}%
$ \\ \cline{2-8}
\end{tabular}%
\medskip \newline
by Result 1. But then 1 manipulates from $u5$ to $u4$. \ \ \ $\square $%
\medskip

\bigskip

\textbf{4. N}$-$\textbf{Range Theorem: Part 2}.\medskip

In this section we prove

\textbf{Theorem 4-1}. \ (The N-Range Theorem, Part 2). \ If $m=3$ and $%
g|NP^{\ast }$ has range of just one alternative then $g$ must be of less
than full range on $NP$.\medskip

The problem here is in one sense easier and in another sense harder than in
Section 3. It is easier because we don't have to treat the many possible
comprehensive collections of decisive sets that had to be considered there.
It is harder because we have to show something a bit more complicated. In
Section 3, we showed a contrapositive:\medskip

$\qquad \qquad |Range(g|NP^{^{\ast }})|$ $=2$ implies $|Range(g)|$ $=2$%
.\medskip

But\medskip

$\qquad \qquad |Range(g|NP^{^{\ast }})|$ $=1$ does not imply $|Range(g)|$ $%
=1 $.\medskip

\textbf{Example 1.} For $n>3$, let $X=\{x,y,z\}$ and define $g$ on $NP$ as
follows: Only $x$ and $y$ are ever chosen; $x$ is chosen unless everyone in $%
\{1,2,...,n-2\}$ and one of $n-1$ and $n$ prefer $y$ to $x$, in which case, $%
y$ is chosen. Then $g$ is strategy-proof with $Range(g)=\{x,y\}$, but Range$%
(g|NP^{^{\ast }})=\{x\}$.\medskip

\textbf{Section 4-1. Lists}\medskip

We will assume $g$ is strategy-proof with $Range(g|NP^{^{\ast }})=\{x\}$ and
then construct a very short list $L$ of profiles such that if $y$ is chosen
at any profile in $NP$ it will also have to be chosen at a profile in $L$
and also a very short list $L^{^{\ast }}$ of profiles such that if $z$ is
chosen at any profile in $NP$ it will also have to be chosen at a profile in 
$L^{^{\ast }}$. Then, for each profile $u$ in $L$, we will show that if $%
y=g(u)$, then for every profile $u^{\ast }$ in $L^{^{\ast }}$, $g(u^{\ast })$
will not be $z$, so that $y=g(u)$ implies $z$ is not chosen at any profile
in $NP$ and thus that $g$ is not of full range.\medskip

To construct list $L$, we will analyze which profiles $u$ (not in $%
NP^{^{\ast }}$) could have $g(u)=y$ by paying attention to the positions of $%
x$ in $u(n-1)$ and $u(n)$. Once this is done, then $L^{^{\ast }}$ can be
constructed by interchanging $y$ and $z$ in the profiles in $L$.\medskip

\textbf{Case 1}. Suppose $g(u)=y$ and $x$ is at the bottom of both $u(n-1)$
and $u(n)$. Then if we can't switch $y$ and $z$ in $u(n-1)$ and stay in $NP$
and so $NP^{\ast }$ to get profile $u^{\prime }$ where $x$ is chosen, then
we can switch $y$ and $z$ in $u(n)$ and stay in $NP$ and so $NP^{\ast }$to
get profile $u^{\prime }$ where $x$ is chosen. Then the individual switching
has an incentive to manipulate back from $u^{\prime }$ to $u$, violating
strategy-proofness. Similarly, $g(u)=z$ would lead to a violation of
strategy-proofness.\medskip

\textbf{Case 2.} Suppose $g(u)=y$ and $x$ is at the top of both $u(n-1)$ and 
$u(n)$. Then if we can't switch $y$ and $z$ in $u(n-1)$ and stay in $NP$ and
so $NP^{\ast }$to get profile $u^{\prime }$ where $x$ is chosen, then we can
switch $y$ and $z$ in $u(n)$ and stay in $NP$ and so $NP^{\ast }$to get
profile $u^{\prime }$ where $x$ is chosen. Then the individual switching has
an incentive to manipulate from $u$ to $u^{\prime }$, violating
strategy-proofness. Similarly, $g(u)=z$ would lead to a violation of
strategy-proofness.\medskip

We now interrupt this sequence of case-by-case analyses to learn two useful
principles.\medskip

\textbf{Lemma 4-2}. If even one of $\{1,2,...,n-2\}$ has $x$ at the top at $%
u $, then $g(u)=x$.\medskip

\qquad Without loss of generality, suppose \#1 has $x$ at the top. Then
construct $u^{\prime }$ by bringing $x$ to the bottom for everyone else,
leaving everyone's ordering of $y$ and $z$ unchanged. Then $u^{\prime }$ is
in $NP$ and Case 1 implies $g(u^{\prime })=x$. Consider the standard
sequence from $u^{\prime }$ to $u$. Each profile in this sequence is in $NP$
and strategy-proofness implies $x$ is chosen at each stage. So $g(u)=x$%
.\medskip

\textbf{Lemma 4-3}. If at $u\in NP$ we have $g(u)=y$, then none of the
individuals in $\{1,2,...,n-2\}$ has $y$ at the bottom.\medskip

\qquad Suppose to the contrary $g(u)=y$; without loss of generality, suppose
\#1 has the ordering $1:zxy$ (it can't be $xzy$ by Lemma 4-2). If anyone
else has $z$ preferred to $x$, then \#1 could interchange $x$ and $z$ to get 
$x$ (by Lemma 4-2) and gain. Strategy-proofness thus implies $x\succ _{i}z$
for all $i$ in $\{2,3,...,n\}$. If any of these individuals has $y$ at top
or bottom, that individual could interchange $x$ and $z$ and still leave $y$
chosen. But then \#1 could interchange $x$ and $z$ to get $x$ (by Lemma 4-2)
and gain. Therefore everyone in $\{2,3,...,n\}$ must have the ordering $xyz$%
. But then $u(n-1)=u(n)$ and we would violate Range$(g|NP^{^{\ast }})=\{x\}$%
.\medskip

By Lemmas 4-2 and 4-3, if $g(u)=y$, then all of $1,2,...,n-2$ have $y\succ
_{i}x$ at $u$. So at least one of $n-1$ and $n$ has $x\succ _{i}y$. We now
we return to our case-by-case analysis.\medskip

\textbf{Case 3}. One of $n-1$ and $n$ has $x$ on the bottom and the other
has $x$ in the middle. Without loss of generality $u$ is\medskip

\begin{tabular}{|c|c|c|c|c|c|c|}
\hline
$1$ & $2$ & $3$ & $\cdots $ & $n-2$ & $n-1$ & $n$ \\ \hline
$%
\begin{array}{c}
y \\ 
x \\ 
\end{array}%
$ & $%
\begin{array}{c}
y \\ 
x \\ 
\end{array}%
$ & $%
\begin{array}{c}
y \\ 
x \\ 
\end{array}%
$ & $\cdots $ & $%
\begin{array}{c}
y \\ 
x \\ 
\end{array}%
$ & $%
\begin{array}{c}
y \\ 
z \\ 
x%
\end{array}%
$ & $%
\begin{array}{c}
z \\ 
x \\ 
y%
\end{array}%
$ \\ \hline
\end{tabular}%
\medskip \newline
and $g(u)=y$. Now $z$ can be brought to the bottom for 1 and 2 while for $%
3,...,n-2$, $x$ is brought to the bottom and y to the top, staying in $NP$,
and still have $y$ chosen by strategy-proofness. This is $L1$, the first
profile in list $L$:\medskip

\begin{tabular}{|c|c|c|c|c|c|c|}
\hline
$1$ & $2$ & $3$ & $\cdots $ & $n-2$ & $n-1$ & $n$ \\ \hline
$%
\begin{array}{c}
y \\ 
x \\ 
z%
\end{array}%
$ & $%
\begin{array}{c}
y \\ 
x \\ 
z%
\end{array}%
$ & $%
\begin{array}{c}
y \\ 
z \\ 
x%
\end{array}%
$ & $\cdots $ & $%
\begin{array}{c}
y \\ 
z \\ 
x%
\end{array}%
$ & $%
\begin{array}{c}
y \\ 
z \\ 
x%
\end{array}%
$ & $%
\begin{array}{c}
z \\ 
x \\ 
y%
\end{array}%
$ \\ \hline
\end{tabular}%
\medskip \newline
(along with other profiles that can be transformed from these by switching $%
x $ and $z$ everywhere; or switching preferences of $n-1$ and $n$ or, for $%
L1 $ or, later, $L2$, choosing a different individual with different
ordering from others in $\{1,2,...,n-2\}$, but all of these will be
equivalent in that if a strategy-proof rule could have $x$ chosen at one but
not on $NP^{^{\ast }}$, then a rule could be designed that had $y$ chosen at
any other one but not on $NP^{^{\ast }}$).\medskip

\textbf{Case 4}. Here $x$ is in the middle for both. Since we are not in $%
NP^{^{\ast }}$, $u(n)$ must be the inverse of $u(n-1)$. Without loss of
generality, $u$ is\medskip

\begin{tabular}{|c|c|c|c|c|c|c|}
\hline
$1$ & $2$ & $3$ & $\cdots $ & $n-2$ & $n-1$ & $n$ \\ \hline
$%
\begin{array}{c}
y \\ 
x \\ 
\end{array}%
$ & $%
\begin{array}{c}
y \\ 
x \\ 
\end{array}%
$ & $%
\begin{array}{c}
y \\ 
x \\ 
\end{array}%
$ & $\cdots $ & $%
\begin{array}{c}
y \\ 
x \\ 
\end{array}%
$ & $%
\begin{array}{c}
y \\ 
x \\ 
z%
\end{array}%
$ & $%
\begin{array}{c}
z \\ 
x \\ 
y%
\end{array}%
$ \\ \hline
\end{tabular}%
\medskip \newline
and $g(u)=y$. We will still have $y$ chosen if $z$ is brought to the bottom
for 1 and 2, while for $3,...,n-2$, $x$ is brought to the bottom and y to
the top, and then interchange $x$ and $z$ for $n-1$. But that is $L1$, so we
don't need to expand list $L$ for Case 4. That is, if $y$ is chosen at a
Case 4 profile, it must also be chosen at a Case 1 profile.\medskip

\textbf{Case 5}. One of $n-1$ and $n$ has $x$ on the top and the other has $%
x $ in the middle; without loss of generality, $n-1$ has $x$ on the top. We
split this into subcases depending on which of $y$ and $z$ is at the top of $%
n$'s ordering.\medskip

\textbf{Subcase 5}-\textbf{1}. Individual $n$ has $z$ on top; $u$ is\medskip

\begin{tabular}{|c|c|c|c|c|c|c|}
\hline
$1$ & $2$ & $3$ & $\cdots $ & $n-2$ & $n-1$ & $n$ \\ \hline
$%
\begin{array}{c}
y \\ 
x \\ 
\end{array}%
$ & $%
\begin{array}{c}
y \\ 
x \\ 
\end{array}%
$ & $%
\begin{array}{c}
y \\ 
x \\ 
\end{array}%
$ & $\cdots $ & $%
\begin{array}{c}
y \\ 
x \\ 
\end{array}%
$ & $%
\begin{array}{c}
x \\ 
\text{ } \\ 
\text{ }%
\end{array}%
$ & $%
\begin{array}{c}
z \\ 
x \\ 
y%
\end{array}%
$ \\ \hline
\end{tabular}%
\medskip \newline
with $g(u)=y$. If anyone in $\{1,2,...,n-2\}$ has $z$ above $x$, then we
could raise $x$ to the top for that individual and still have $y$ chosen,
contrary to Lemma 4-2. So everyone in $\{1,2,...,n-2\}$ has $x$ above $z$%
:\medskip 

\begin{tabular}{|c|c|c|c|c|c|c|}
\hline
$1$ & $2$ & $3$ & $\cdots $ & $n-2$ & $n-1$ & $n$ \\ \hline
$%
\begin{array}{c}
y \\ 
x \\ 
z%
\end{array}%
$ & $%
\begin{array}{c}
y \\ 
x \\ 
z%
\end{array}%
$ & $%
\begin{array}{c}
y \\ 
x \\ 
z%
\end{array}%
$ & $\cdots $ & $%
\begin{array}{c}
y \\ 
x \\ 
z%
\end{array}%
$ & $%
\begin{array}{c}
x \\ 
\text{ } \\ 
\text{ }%
\end{array}%
$ & $%
\begin{array}{c}
z \\ 
x \\ 
y%
\end{array}%
$ \\ \hline
\end{tabular}%
\medskip \newline
But then raising $y$ to the top for $n-1$ and lowering $x$ to the bottom for 
$3,...,n-2$ leaves $y$ chosen and we are back to profile $L1$, so we don't
need to expand list $L$ for Subcase 5-1.\medskip

\textbf{Subcase 5}-\textbf{2}. Individual $n$ has $y$ on top; $u$ is\medskip

\begin{tabular}{|c|c|c|c|c|c|c|}
\hline
$1$ & $2$ & $3$ & $\cdots $ & $n-2$ & $n-1$ & $n$ \\ \hline
$%
\begin{array}{c}
y \\ 
x \\ 
\end{array}%
$ & $%
\begin{array}{c}
y \\ 
x \\ 
\end{array}%
$ & $%
\begin{array}{c}
y \\ 
x \\ 
\end{array}%
$ & $\cdots $ & $%
\begin{array}{c}
y \\ 
x \\ 
\end{array}%
$ & $%
\begin{array}{c}
x \\ 
\text{ } \\ 
\text{ }%
\end{array}%
$ & $%
\begin{array}{c}
y \\ 
x \\ 
z%
\end{array}%
$ \\ \hline
\end{tabular}%
\medskip \newline
with $g(u)=y$.\medskip

\qquad If $n-1$ has $y$ preferred to $z$, then at least one of $1,2,...,n-2$
must have $z$ preferred to $y$, say individual 1. Bring $z$ to the bottom
for the others in $\{1,2,...,n-2\}$ and $y$ is chosen at say\medskip

\begin{tabular}{|c|c|c|c|c|c|c|}
\hline
$1$ & $2$ & $3$ & $\cdots $ & $n-2$ & $n-1$ & $n$ \\ \hline
$%
\begin{array}{c}
z \\ 
y \\ 
x%
\end{array}%
$ & $%
\begin{array}{c}
y \\ 
x \\ 
z%
\end{array}%
$ & $%
\begin{array}{c}
y \\ 
x \\ 
z%
\end{array}%
$ & $\cdots $ & $%
\begin{array}{c}
y \\ 
x \\ 
z%
\end{array}%
$ & $%
\begin{array}{c}
x \\ 
y \\ 
z%
\end{array}%
$ & $%
\begin{array}{c}
y \\ 
x \\ 
z%
\end{array}%
$ \\ \hline
\end{tabular}%
\medskip \newline
Then interchange $x$ and $z$ below $y$ for $3,...,n-2$ to get\medskip

\begin{tabular}{|c|c|c|c|c|c|c|}
\hline
$1$ & $2$ & $3$ & $\cdots $ & $n-2$ & $n-1$ & $n$ \\ \hline
$%
\begin{array}{c}
z \\ 
y \\ 
x%
\end{array}%
$ & $%
\begin{array}{c}
y \\ 
x \\ 
z%
\end{array}%
$ & $%
\begin{array}{c}
y \\ 
z \\ 
x%
\end{array}%
$ & $\cdots $ & $%
\begin{array}{c}
y \\ 
z \\ 
x%
\end{array}%
$ & $%
\begin{array}{c}
x \\ 
y \\ 
z%
\end{array}%
$ & $%
\begin{array}{c}
y \\ 
x \\ 
z%
\end{array}%
$ \\ \hline
\end{tabular}%
\medskip \newline
which we add as $L2$, to the list $L$ (along with related profiles as
remarked at the end of Case 3).\medskip

\qquad On the other hand, if $n-1$ has $z$ preferred to $y$, then we observe
that at least one of $1,2,...,n-2$ must have $z$ preferred to $x$, say
individual 2. Bring $z$ to the bottom for 1 and 2 while bringing $x$ to the
bottom for $3,...,n-2$ and then raise $y$ to the top for \#1 and $y$ is
still chosen at\medskip

\begin{tabular}{|c|c|c|c|c|c|c|}
\hline
$1$ & $2$ & $3$ & $\cdots $ & $n-2$ & $n-1$ & $n$ \\ \hline
$%
\begin{array}{c}
y \\ 
z \\ 
x%
\end{array}%
$ & $%
\begin{array}{c}
y \\ 
x \\ 
z%
\end{array}%
$ & $%
\begin{array}{c}
y \\ 
z \\ 
x%
\end{array}%
$ & $\cdots $ & $%
\begin{array}{c}
y \\ 
z \\ 
x%
\end{array}%
$ & $%
\begin{array}{c}
x \\ 
z \\ 
y%
\end{array}%
$ & $%
\begin{array}{c}
y \\ 
x \\ 
z%
\end{array}%
$ \\ \hline
\end{tabular}%
\medskip \newline
This is $L3$, the third profile we add to $L$ (along with related
profiles).\medskip

\textbf{Case 6}. One of $n-1$ and $n$ has $x$ on the top and the other has $%
x $ on the bottom; without loss of generality, $n-1$ has $x$ on the
top.\medskip

\begin{tabular}{|c|c|c|c|c|c|c|}
\hline
$1$ & $2$ & $3$ & $\cdots $ & $n-2$ & $n-1$ & $n$ \\ \hline
$%
\begin{array}{c}
y \\ 
x \\ 
\end{array}%
$ & $%
\begin{array}{c}
y \\ 
x \\ 
\end{array}%
$ & $%
\begin{array}{c}
y \\ 
x \\ 
\end{array}%
$ & $\cdots $ & $%
\begin{array}{c}
y \\ 
x \\ 
\end{array}%
$ & $%
\begin{array}{c}
x \\ 
\text{ } \\ 
\text{ }%
\end{array}%
$ & $%
\begin{array}{c}
\\ 
\\ 
x%
\end{array}%
$ \\ \hline
\end{tabular}%
\medskip \newline
Someone has to have $z\succ _{i}y$. There are three possibilities.\medskip

\qquad \textbf{A}. One of $\{1,2,...,n-2\}$ has $z\succ _{i}y$, say \#1.
Then bring $z$ to the bottom for everyone else:\medskip

\begin{tabular}{|c|c|c|c|c|c|c|}
\hline
$1$ & $2$ & $3$ & $\cdots $ & $n-2$ & $n-1$ & $n$ \\ \hline
$%
\begin{array}{c}
z \\ 
y \\ 
x%
\end{array}%
$ & $%
\begin{array}{c}
y \\ 
x \\ 
z%
\end{array}%
$ & $%
\begin{array}{c}
y \\ 
x \\ 
z%
\end{array}%
$ & $\cdots $ & $%
\begin{array}{c}
y \\ 
x \\ 
z%
\end{array}%
$ & $%
\begin{array}{c}
x \\ 
y \\ 
z%
\end{array}%
$ & $%
\begin{array}{c}
y \\ 
x \\ 
z%
\end{array}%
$ \\ \hline
\end{tabular}%
\medskip \newline
and $y$ is still chosen. This puts us in Case 5 and no new profiles need to
added to list $L$.\medskip

\qquad \textbf{B}. Individual $n$ has $z\succ _{n}y$. Then $z$ can be
brought to the bottom for 1, 2, and $n-1$, while bringing $x$ to the bottom
for $3,...,n-2$ and still have $y$ chosen:\medskip

\begin{tabular}{|c|c|c|c|c|c|c|}
\hline
$1$ & $2$ & $3$ & $\cdots $ & $n-2$ & $n-1$ & $n$ \\ \hline
$%
\begin{array}{c}
y \\ 
x \\ 
z%
\end{array}%
$ & $%
\begin{array}{c}
y \\ 
x \\ 
z%
\end{array}%
$ & $%
\begin{array}{c}
y \\ 
z \\ 
x%
\end{array}%
$ & $\cdots $ & $%
\begin{array}{c}
y \\ 
z \\ 
x%
\end{array}%
$ & $%
\begin{array}{c}
x \\ 
y \\ 
z%
\end{array}%
$ & $%
\begin{array}{c}
z \\ 
y \\ 
x%
\end{array}%
$ \\ \hline
\end{tabular}%
\medskip \newline
which is $L4$, the fourth profile in list $L$.\medskip

\qquad \textbf{C}. Individual $n-1$ has $z\succ _{n-1}y$ and everyone else
has $y\succ _{i}z$ (or we are back in A or B):\medskip

\begin{tabular}{|c|c|c|c|c|c|c|}
\hline
$1$ & $2$ & $3$ & $\cdots $ & $n-2$ & $n-1$ & $n$ \\ \hline
$%
\begin{array}{c}
y \\ 
\lbrack x,z] \\ 
\end{array}%
$ & $%
\begin{array}{c}
y \\ 
\lbrack x,z] \\ 
\end{array}%
$ & $%
\begin{array}{c}
y \\ 
\lbrack x,z] \\ 
\end{array}%
$ & $\cdots $ & $%
\begin{array}{c}
y \\ 
\lbrack x,z] \\ 
\end{array}%
$ & $%
\begin{array}{c}
x \\ 
z \\ 
y%
\end{array}%
$ & $%
\begin{array}{c}
y \\ 
z \\ 
x%
\end{array}%
$ \\ \hline
\end{tabular}%
\medskip \newline
where [x,z] indicates that x and z can be ordered in any manner below y for
each individual. \ Then $y$ will still be chosen if we make $z$ just above $%
x $ for \#1 and then raise $x$ just above $z$ in n's ordering. This puts us
in Case 5 and no new profiles need to added to list $L$.\medskip

Summarizing, the list $L$ consists of\medskip

\begin{tabular}{c|c|c|c|c|c|c|c|}
\cline{2-8}
$L1$ & $1$ & $2$ & $3$ & $\cdots $ & $n-2$ & $n-1$ & $n$ \\ \cline{2-8}
& $%
\begin{array}{c}
y \\ 
x \\ 
z%
\end{array}%
$ & $%
\begin{array}{c}
y \\ 
x \\ 
z%
\end{array}%
$ & $%
\begin{array}{c}
y \\ 
z \\ 
x%
\end{array}%
$ & $\cdots $ & $%
\begin{array}{c}
y \\ 
z \\ 
x%
\end{array}%
$ & $%
\begin{array}{c}
y \\ 
z \\ 
x%
\end{array}%
$ & $%
\begin{array}{c}
z \\ 
x \\ 
y%
\end{array}%
$ \\ \cline{2-8}
\end{tabular}%
\medskip \newline

\begin{tabular}{c|c|c|c|c|c|c|c|}
\cline{2-8}
$L2$ & $1$ & $2$ & $3$ & $\cdots $ & $n-2$ & $n-1$ & $n$ \\ \cline{2-8}
& $%
\begin{array}{c}
z \\ 
y \\ 
x%
\end{array}%
$ & $%
\begin{array}{c}
y \\ 
x \\ 
z%
\end{array}%
$ & $%
\begin{array}{c}
y \\ 
z \\ 
x%
\end{array}%
$ & $\cdots $ & $%
\begin{array}{c}
y \\ 
z \\ 
x%
\end{array}%
$ & $%
\begin{array}{c}
x \\ 
y \\ 
z%
\end{array}%
$ & $%
\begin{array}{c}
y \\ 
x \\ 
z%
\end{array}%
$ \\ \cline{2-8}
\end{tabular}%
\medskip \newline

\begin{tabular}{c|c|c|c|c|c|c|c|}
\cline{2-8}
$L3$ & $1$ & $2$ & $3$ & $\cdots $ & $n-2$ & $n-1$ & $n$ \\ \cline{2-8}
& $%
\begin{array}{c}
y \\ 
z \\ 
x%
\end{array}%
$ & $%
\begin{array}{c}
y \\ 
x \\ 
z%
\end{array}%
$ & $%
\begin{array}{c}
y \\ 
z \\ 
x%
\end{array}%
$ & $\cdots $ & $%
\begin{array}{c}
y \\ 
z \\ 
x%
\end{array}%
$ & $%
\begin{array}{c}
x \\ 
z \\ 
y%
\end{array}%
$ & $%
\begin{array}{c}
y \\ 
x \\ 
z%
\end{array}%
$ \\ \cline{2-8}
\end{tabular}%
\medskip \newline

\begin{tabular}{c|c|c|c|c|c|c|c|}
\cline{2-8}
$L4$ & $1$ & $2$ & $3$ & $\cdots $ & $n-2$ & $n-1$ & $n$ \\ \cline{2-8}
& $%
\begin{array}{c}
y \\ 
x \\ 
z%
\end{array}%
$ & $%
\begin{array}{c}
y \\ 
x \\ 
z%
\end{array}%
$ & $%
\begin{array}{c}
y \\ 
z \\ 
x%
\end{array}%
$ & $\cdots $ & $%
\begin{array}{c}
y \\ 
z \\ 
x%
\end{array}%
$ & $%
\begin{array}{c}
x \\ 
y \\ 
z%
\end{array}%
$ & $%
\begin{array}{c}
z \\ 
y \\ 
x%
\end{array}%
$ \\ \cline{2-8}
\end{tabular}%
\medskip \newline
along with other profiles that can be transformed from these by switching $x$
and $z$ everywhere; or switching preferences of $n-1$ and $n$ or, for $L1$
or $L2$, choosing an alternative individual with different ordering from
others in $\{1,2,...,n-2\}$, but all of these will be equivalent in that if
a strategy-proof rule could have $x$ chosen at one but not on $NP^{^{\ast }}$%
, then a rule could be designed that had $y$ chosen at any other one but not
on $NP^{^{\ast }}$.\medskip

An analogous argument allows the construction of list $L^{^{\ast }}$ such
that if $z$ is chosen at any profile in $NP$, then $z$ must be chosen at a
profile in the list $L^{^{\ast }}$. $L^{^{\ast }}$ can be obtained simply by
interchanging $y$ and $z$ in list $L$ (along with other profiles that can be
transformed from these by switching $x$ and $y$ everywhere; or switching
preferences of $n-1$ and $n$ or permuting the individuals in $%
\{1,2,...,n-2\} $, but all of these will be equivalent in that if a
strategy-proof rule could have $z$ chosen at one but not on $NP^{^{\ast }}$,
then a rule could be designed that had $z$ chosen at any other one but not
on $NP^{^{\ast }}$).\medskip

List $L^{^{\ast }}$ is:\medskip

\begin{tabular}{c|c|c|c|c|c|c|c|}
\cline{2-8}
$L1^{^{\ast }}$ & $1$ & $2$ & $3$ & $\cdots $ & $n-2$ & $n-1$ & $n$ \\ 
\cline{2-8}
& $%
\begin{array}{c}
z \\ 
x \\ 
y%
\end{array}%
$ & $%
\begin{array}{c}
z \\ 
x \\ 
y%
\end{array}%
$ & $%
\begin{array}{c}
z \\ 
y \\ 
x%
\end{array}%
$ & $\cdots $ & $%
\begin{array}{c}
z \\ 
y \\ 
x%
\end{array}%
$ & $%
\begin{array}{c}
z \\ 
y \\ 
x%
\end{array}%
$ & $%
\begin{array}{c}
y \\ 
x \\ 
z%
\end{array}%
$ \\ \cline{2-8}
\end{tabular}%
\medskip \newline

\begin{tabular}{c|c|c|c|c|c|c|c|}
\cline{2-8}
$L2^{^{\ast }}$ & $1$ & $2$ & $3$ & $\cdots $ & $n-2$ & $n-1$ & $n$ \\ 
\cline{2-8}
& $%
\begin{array}{c}
y \\ 
z \\ 
x%
\end{array}%
$ & $%
\begin{array}{c}
z \\ 
x \\ 
y%
\end{array}%
$ & $%
\begin{array}{c}
z \\ 
y \\ 
x%
\end{array}%
$ & $\cdots $ & $%
\begin{array}{c}
z \\ 
y \\ 
x%
\end{array}%
$ & $%
\begin{array}{c}
x \\ 
z \\ 
y%
\end{array}%
$ & $%
\begin{array}{c}
z \\ 
x \\ 
y%
\end{array}%
$ \\ \cline{2-8}
\end{tabular}%
\medskip \newline

\begin{tabular}{c|c|c|c|c|c|c|c|}
\cline{2-8}
$L3^{^{\ast }}$ & $1$ & $2$ & $3$ & $\cdots $ & $n-2$ & $n-1$ & $n$ \\ 
\cline{2-8}
& $%
\begin{array}{c}
z \\ 
y \\ 
x%
\end{array}%
$ & $%
\begin{array}{c}
z \\ 
x \\ 
y%
\end{array}%
$ & $%
\begin{array}{c}
z \\ 
y \\ 
x%
\end{array}%
$ & $\cdots $ & $%
\begin{array}{c}
z \\ 
y \\ 
x%
\end{array}%
$ & $%
\begin{array}{c}
x \\ 
y \\ 
z%
\end{array}%
$ & $%
\begin{array}{c}
z \\ 
x \\ 
y%
\end{array}%
$ \\ \cline{2-8}
\end{tabular}%
\medskip \newline

\begin{tabular}{c|c|c|c|c|c|c|c|}
\cline{2-8}
$L4^{^{\ast }}$ & $1$ & $2$ & $3$ & $\cdots $ & $n-2$ & $n-1$ & $n$ \\ 
\cline{2-8}
& $%
\begin{array}{c}
z \\ 
x \\ 
y%
\end{array}%
$ & $%
\begin{array}{c}
z \\ 
x \\ 
y%
\end{array}%
$ & $%
\begin{array}{c}
z \\ 
y \\ 
x%
\end{array}%
$ & $\cdots $ & $%
\begin{array}{c}
z \\ 
y \\ 
x%
\end{array}%
$ & $%
\begin{array}{c}
x \\ 
z \\ 
y%
\end{array}%
$ & $%
\begin{array}{c}
y \\ 
z \\ 
x%
\end{array}%
$ \\ \cline{2-8}
\end{tabular}%
\medskip \newline
Our goal then is, for each choice of $Li$ to assume $g(Li)=y$, to then show
that each $g(Lj^{\ast })\neq z$. We actually will establish the stronger
result that each $g(Lj^{\ast })=x$. We can simplify our analyses by
introducing yet another list, $L^{\ast \ast }$:\medskip

\begin{tabular}{c|c|c|c|c|c|c|c|}
\cline{2-8}
$L1^{^{\ast \ast }}$ & $1$ & $2$ & $3$ & $\cdots $ & $n-2$ & $n-1$ & $n$ \\ 
\cline{2-8}
& $%
\begin{array}{c}
z \\ 
x \\ 
y%
\end{array}%
$ & $%
\begin{array}{c}
z \\ 
x \\ 
y%
\end{array}%
$ & $%
\begin{array}{c}
y \\ 
z \\ 
x%
\end{array}%
$ & $\cdots $ & $%
\begin{array}{c}
y \\ 
z \\ 
x%
\end{array}%
$ & $%
\begin{array}{c}
z \\ 
y \\ 
x%
\end{array}%
$ & $%
\begin{array}{c}
y \\ 
x \\ 
z%
\end{array}%
$ \\ \cline{2-8}
\end{tabular}%
\medskip

\begin{tabular}{c|c|c|c|c|c|c|c|}
\cline{2-8}
$L2^{^{\ast \ast }}$ & $1$ & $2$ & $3$ & $\cdots $ & $n-2$ & $n-1$ & $n$ \\ 
\cline{2-8}
& $%
\begin{array}{c}
y \\ 
z \\ 
x%
\end{array}%
$ & $%
\begin{array}{c}
z \\ 
x \\ 
y%
\end{array}%
$ & $%
\begin{array}{c}
y \\ 
z \\ 
x%
\end{array}%
$ & $\cdots $ & $%
\begin{array}{c}
y \\ 
z \\ 
x%
\end{array}%
$ & $%
\begin{array}{c}
x \\ 
z \\ 
y%
\end{array}%
$ & $%
\begin{array}{c}
z \\ 
x \\ 
y%
\end{array}%
$ \\ \cline{2-8}
\end{tabular}%
\medskip

\begin{tabular}{c|c|c|c|c|c|c|c|}
\cline{2-8}
$L3^{^{\ast \ast }}$ & $1$ & $2$ & $3$ & $\cdots $ & $n-2$ & $n-1$ & $n$ \\ 
\cline{2-8}
& $%
\begin{array}{c}
z \\ 
y \\ 
x%
\end{array}%
$ & $%
\begin{array}{c}
z \\ 
x \\ 
y%
\end{array}%
$ & $%
\begin{array}{c}
y \\ 
z \\ 
x%
\end{array}%
$ & $\cdots $ & $%
\begin{array}{c}
y \\ 
z \\ 
x%
\end{array}%
$ & $%
\begin{array}{c}
x \\ 
y \\ 
z%
\end{array}%
$ & $%
\begin{array}{c}
z \\ 
x \\ 
y%
\end{array}%
$ \\ \cline{2-8}
\end{tabular}%
\medskip

\begin{tabular}{c|c|c|c|c|c|c|c|}
\cline{2-8}
$L4^{^{\ast \ast }}$ & $1$ & $2$ & $3$ & $\cdots $ & $n-2$ & $n-1$ & $n$ \\ 
\cline{2-8}
& $%
\begin{array}{c}
z \\ 
x \\ 
y%
\end{array}%
$ & $%
\begin{array}{c}
z \\ 
x \\ 
y%
\end{array}%
$ & $%
\begin{array}{c}
y \\ 
z \\ 
x%
\end{array}%
$ & $\cdots $ & $%
\begin{array}{c}
y \\ 
z \\ 
x%
\end{array}%
$ & $%
\begin{array}{c}
x \\ 
z \\ 
y%
\end{array}%
$ & $%
\begin{array}{c}
y \\ 
z \\ 
x%
\end{array}%
$ \\ \cline{2-8}
\end{tabular}%
\medskip \newline
For each $Lj^{\ast }$, the profile $Lj^{\ast \ast }$ changes zyx to yzx for
individuals $3,...,n-2$. (Alternative $x$ is also ranked at the bottom of $%
Lj(i)$ for $3\leq i\leq n-2$ for all $j$). \medskip

\textbf{Lemma 4-4}.For a strategy-proof rule g, 
\[
g(Lj^{\ast \ast })=x\text{ }implies\text{ }g(Lj^{\ast })=x.\medskip 
\]

\textbf{Proof}: Just construct a standard sequence, switching each ordering
in turn; strategy-proofness implies that x is chosen at each step. \ \ \ $%
\square $\medskip

\textbf{Lemma 4-5.} \ For a strategy-proof rule g, $g(L3^{\ast \ast })=x$
implies $g(L2^{\ast \ast })=x$ (and so, by Lemma 4-4, $g(L2^{\ast })=x$%
).\medskip

\textbf{Proof}: Suppose $g(L3^{\ast \ast })=x$ at\medskip

\begin{tabular}{c|c|c|c|c|c|c|c|}
\cline{2-8}
$L3^{^{\ast \ast }}$ & $1$ & $2$ & $3$ & $\cdots $ & $n-2$ & $n-1$ & $n$ \\ 
\cline{2-8}
& $%
\begin{array}{c}
z \\ 
y \\ 
x%
\end{array}%
$ & $%
\begin{array}{c}
z \\ 
x \\ 
y%
\end{array}%
$ & $%
\begin{array}{c}
y \\ 
z \\ 
x%
\end{array}%
$ & $\cdots $ & $%
\begin{array}{c}
y \\ 
z \\ 
x%
\end{array}%
$ & $%
\begin{array}{c}
x \\ 
y \\ 
z%
\end{array}%
$ & $%
\begin{array}{c}
z \\ 
x \\ 
y%
\end{array}%
$ \\ \cline{2-8}
\end{tabular}%
\medskip \newline
This is a 1-variant of the profile\medskip

\begin{tabular}{c|c|c|c|c|c|c|c|}
\cline{2-8}
$u1$ & $1$ & $2$ & $3$ & $\cdots $ & $n-2$ & $n-1$ & $n$ \\ \cline{2-8}
& $%
\begin{array}{c}
y \\ 
z \\ 
x%
\end{array}%
$ & $%
\begin{array}{c}
z \\ 
x \\ 
y%
\end{array}%
$ & $%
\begin{array}{c}
y \\ 
z \\ 
x%
\end{array}%
$ & $\cdots $ & $%
\begin{array}{c}
y \\ 
z \\ 
x%
\end{array}%
$ & $%
\begin{array}{c}
x \\ 
y \\ 
z%
\end{array}%
$ & $%
\begin{array}{c}
z \\ 
x \\ 
y%
\end{array}%
$ \\ \cline{2-8}
\end{tabular}%
\medskip \newline
Then $g(u1)=x$ or 1 manipulates from $L3^{\ast \ast }$ to $u1$. Next notice
that $L2^{\ast \ast }$ is an $(n-1)$-variant of $u1$:\medskip

\begin{tabular}{c|c|c|c|c|c|c|c|}
\cline{2-8}
$L2^{\ast \ast }$ & $1$ & $2$ & $3$ & $\cdots $ & $n-2$ & $n-1$ & $n$ \\ 
\cline{2-8}
& $%
\begin{array}{c}
y \\ 
z \\ 
x%
\end{array}%
$ & $%
\begin{array}{c}
z \\ 
x \\ 
y%
\end{array}%
$ & $%
\begin{array}{c}
y \\ 
z \\ 
x%
\end{array}%
$ & $\cdots $ & $%
\begin{array}{c}
y \\ 
z \\ 
x%
\end{array}%
$ & $%
\begin{array}{c}
x \\ 
z \\ 
y%
\end{array}%
$ & $%
\begin{array}{c}
z \\ 
x \\ 
y%
\end{array}%
$ \\ \cline{2-8}
\end{tabular}%
\medskip \newline
so $g(L2^{\ast \ast })=x$ or $n-1$ manipulates from $L2^{\ast \ast }$ to $u1$%
. \ \ \ \ $\square $\medskip

\textbf{Proof of Theorem 4-1}: \ Because of Lemmas 4-4 and 4-5, we will show
that if $g|NP^{^{\ast }}$ has range $\{x\}$ then for each $i=1,2,3,4$, $%
g(Li)=y$ implies $g(Lj^{\ast \ast })=x$ for all $j=1,3,4$ (but see Section
4-2-4).\medskip

\textbf{Section 4-2}. Assume $g(L1)=y$.\medskip

\qquad It continues to be necessary to check that each profile employed is
actually in $NP$. \ But one must do so with reference to the preferences of
individuals $1,2,n-1$, and $n$ only as we might have $n=4$.\medskip

\textbf{Subsection 4-2-1}. Proof that $g(L1^{\ast \ast })=x$.\medskip

\textbf{Step 1}\medskip

Let profile $u1$ be\medskip

\begin{tabular}{c|c|c|c|c|c|c|c|}
\cline{2-8}
$u1$ & $1$ & $2$ & $3$ & $\cdots $ & $n-2$ & $n-1$ & $n$ \\ \cline{2-8}
& $%
\begin{array}{c}
y \\ 
x \\ 
z%
\end{array}%
$ & $%
\begin{array}{c}
z \\ 
x \\ 
y%
\end{array}%
$ & $%
\begin{array}{c}
y \\ 
z \\ 
x%
\end{array}%
$ & $\cdots $ & $%
\begin{array}{c}
y \\ 
z \\ 
x%
\end{array}%
$ & $%
\begin{array}{c}
z \\ 
y \\ 
x%
\end{array}%
$ & $%
\begin{array}{c}
z \\ 
y \\ 
x%
\end{array}%
$ \\ \cline{2-8}
\end{tabular}%
\medskip \newline
Then $g(u1)=x$ since $u1\in NP^{^{\ast }}$. Now consider $n$-variant $u2$%
:\medskip

\begin{tabular}{c|c|c|c|c|c|c|c|}
\cline{2-8}
$u2$ & $1$ & $2$ & $3$ & $\cdots $ & $n-2$ & $n-1$ & $n$ \\ \cline{2-8}
& $%
\begin{array}{c}
y \\ 
x \\ 
z%
\end{array}%
$ & $%
\begin{array}{c}
z \\ 
x \\ 
y%
\end{array}%
$ & $%
\begin{array}{c}
y \\ 
z \\ 
x%
\end{array}%
$ & $\cdots $ & $%
\begin{array}{c}
y \\ 
z \\ 
x%
\end{array}%
$ & $%
\begin{array}{c}
z \\ 
y \\ 
x%
\end{array}%
$ & $%
\begin{array}{c}
z \\ 
x \\ 
y%
\end{array}%
$ \\ \cline{2-8}
\end{tabular}%
\medskip \newline
Still $g(u2)=x$ or $n$ would manipulate from $u1$ to $u2$. Then look at
2-variant $u3$:\medskip

\begin{tabular}{c|c|c|c|c|c|c|c|}
\cline{2-8}
$u3$ & $1$ & $2$ & $3$ & $\cdots $ & $n-2$ & $n-1$ & $n$ \\ \cline{2-8}
& $%
\begin{array}{c}
y \\ 
x \\ 
z%
\end{array}%
$ & $%
\begin{array}{c}
y \\ 
x \\ 
z%
\end{array}%
$ & $%
\begin{array}{c}
y \\ 
z \\ 
x%
\end{array}%
$ & $\cdots $ & $%
\begin{array}{c}
y \\ 
z \\ 
x%
\end{array}%
$ & $%
\begin{array}{c}
z \\ 
y \\ 
x%
\end{array}%
$ & $%
\begin{array}{c}
z \\ 
x \\ 
y%
\end{array}%
$ \\ \cline{2-8}
\end{tabular}%
\medskip \newline
Then $g(u3)\neq z$ or 2 would manipulate from $u2$ to $u3$. But now consider
profile $L1$:\medskip

\begin{tabular}{c|c|c|c|c|c|c|c|}
\cline{2-8}
$L1$ & $1$ & $2$ & $3$ & $\cdots $ & $n-2$ & $n-1$ & $n$ \\ \cline{2-8}
& $%
\begin{array}{c}
y \\ 
x \\ 
z%
\end{array}%
$ & $%
\begin{array}{c}
y \\ 
x \\ 
z%
\end{array}%
$ & $%
\begin{array}{c}
y \\ 
z \\ 
x%
\end{array}%
$ & $\cdots $ & $%
\begin{array}{c}
y \\ 
z \\ 
x%
\end{array}%
$ & $%
\begin{array}{c}
y \\ 
z \\ 
x%
\end{array}%
$ & $%
\begin{array}{c}
z \\ 
x \\ 
y%
\end{array}%
$ \\ \cline{2-8}
\end{tabular}%
\medskip \newline
We are assuming $g(L1)=y$. This profile is an $(n-1)$-variant of $u3$. If $%
g(u3)=x$, individual $n-1$ would manipulate from $u3$ to $L1$. So $g(u3)\neq
x$. Combining with $g(u3)\neq z$, we have $g(u3)=y$.\medskip

\qquad Next consider profile $u4$, a 1-variant of $u3$:\medskip

\begin{tabular}{c|c|c|c|c|c|c|c|}
\cline{2-8}
$u4$ & $1$ & $2$ & $3$ & $\cdots $ & $n-2$ & $n-1$ & $n$ \\ \cline{2-8}
& $%
\begin{array}{c}
z \\ 
y \\ 
x%
\end{array}%
$ & $%
\begin{array}{c}
y \\ 
x \\ 
z%
\end{array}%
$ & $%
\begin{array}{c}
y \\ 
z \\ 
x%
\end{array}%
$ & $\cdots $ & $%
\begin{array}{c}
y \\ 
z \\ 
x%
\end{array}%
$ & $%
\begin{array}{c}
z \\ 
y \\ 
x%
\end{array}%
$ & $%
\begin{array}{c}
z \\ 
x \\ 
y%
\end{array}%
$ \\ \cline{2-8}
\end{tabular}%
\medskip \newline
If $g(u4)=x$, then 1 would manipulate from $u4$ to $u3$; so $g(u4)\neq x$%
.\medskip

\textbf{Step 2}\medskip

Consider profile $u5$ in $NP^{^{\ast }}$:\medskip

\begin{tabular}{c|c|c|c|c|c|c|c|}
\cline{2-8}
$u5$ & $1$ & $2$ & $3$ & $\cdots $ & $n-2$ & $n-1$ & $n$ \\ \cline{2-8}
& $%
\begin{array}{c}
z \\ 
x \\ 
y%
\end{array}%
$ & $%
\begin{array}{c}
y \\ 
x \\ 
z%
\end{array}%
$ & $%
\begin{array}{c}
y \\ 
z \\ 
x%
\end{array}%
$ & $\cdots $ & $%
\begin{array}{c}
y \\ 
z \\ 
x%
\end{array}%
$ & $%
\begin{array}{c}
z \\ 
y \\ 
x%
\end{array}%
$ & $%
\begin{array}{c}
z \\ 
y \\ 
x%
\end{array}%
$ \\ \cline{2-8}
\end{tabular}%
\medskip \newline
Then $g(u5)=x$. Look at $n$-variant profile $u6$\medskip

\begin{tabular}{c|c|c|c|c|c|c|c|}
\cline{2-8}
$u6$ & $1$ & $2$ & $3$ & $\cdots $ & $n-2$ & $n-1$ & $n$ \\ \cline{2-8}
& $%
\begin{array}{c}
z \\ 
x \\ 
y%
\end{array}%
$ & $%
\begin{array}{c}
y \\ 
x \\ 
z%
\end{array}%
$ & $%
\begin{array}{c}
y \\ 
z \\ 
x%
\end{array}%
$ & $\cdots $ & $%
\begin{array}{c}
y \\ 
z \\ 
x%
\end{array}%
$ & $%
\begin{array}{c}
z \\ 
y \\ 
x%
\end{array}%
$ & $%
\begin{array}{c}
z \\ 
x \\ 
y%
\end{array}%
$ \\ \cline{2-8}
\end{tabular}%
\medskip \newline
We have $g(u6)=x$ or else $n$ would manipulate from $u5$ to $u6$. But $u6$
is also a 1-variant of $u4$. If $g(u4)=z$, 1 would manipulate from $u6$ to $%
u4$. So $g(u4)\neq z$.\medskip

\textbf{Step 3}\medskip

Combining Steps 1 and 2, $g(u4)=y$. Then look at $u7$, an $n$-variant of $u4$%
:\medskip

\begin{tabular}{c|c|c|c|c|c|c|c|}
\cline{2-8}
$u7$ & $1$ & $2$ & $3$ & $\cdots $ & $n-2$ & $n-1$ & $n$ \\ \cline{2-8}
& $%
\begin{array}{c}
z \\ 
y \\ 
x%
\end{array}%
$ & $%
\begin{array}{c}
y \\ 
x \\ 
z%
\end{array}%
$ & $%
\begin{array}{c}
y \\ 
z \\ 
x%
\end{array}%
$ & $\cdots $ & $%
\begin{array}{c}
y \\ 
z \\ 
x%
\end{array}%
$ & $%
\begin{array}{c}
z \\ 
y \\ 
x%
\end{array}%
$ & $%
\begin{array}{c}
x \\ 
y \\ 
z%
\end{array}%
$ \\ \cline{2-8}
\end{tabular}%
\medskip \newline
$g(u7)=y$ or $n$ manipulates from $u4$ to $u7$. Consider the 2-variant of $%
u7 $:\medskip

\begin{tabular}{c|c|c|c|c|c|c|c|}
\cline{2-8}
$u8$ & $1$ & $2$ & $3$ & $\cdots $ & $n-2$ & $n-1$ & $n$ \\ \cline{2-8}
& $%
\begin{array}{c}
z \\ 
y \\ 
x%
\end{array}%
$ & $%
\begin{array}{c}
y \\ 
z \\ 
x%
\end{array}%
$ & $%
\begin{array}{c}
y \\ 
z \\ 
x%
\end{array}%
$ & $\cdots $ & $%
\begin{array}{c}
y \\ 
z \\ 
x%
\end{array}%
$ & $%
\begin{array}{c}
z \\ 
y \\ 
x%
\end{array}%
$ & $%
\begin{array}{c}
x \\ 
y \\ 
z%
\end{array}%
$ \\ \cline{2-8}
\end{tabular}%
\medskip \newline
$g(u8)=y$ or 2 would manipulate from $u8$ to $u7$. Then at profile $u9$, a
1-variant of $u8$:\medskip

\begin{tabular}{c|c|c|c|c|c|c|c|}
\cline{2-8}
$u9$ & $1$ & $2$ & $3$ & $\cdots $ & $n-2$ & $n-1$ & $n$ \\ \cline{2-8}
& $%
\begin{array}{c}
z \\ 
x \\ 
y%
\end{array}%
$ & $%
\begin{array}{c}
y \\ 
z \\ 
x%
\end{array}%
$ & $%
\begin{array}{c}
y \\ 
z \\ 
x%
\end{array}%
$ & $\cdots $ & $%
\begin{array}{c}
y \\ 
z \\ 
x%
\end{array}%
$ & $%
\begin{array}{c}
z \\ 
y \\ 
x%
\end{array}%
$ & $%
\begin{array}{c}
x \\ 
y \\ 
z%
\end{array}%
$ \\ \cline{2-8}
\end{tabular}%
\medskip \newline
we have $g(u9)\neq z$ or 1 would manipulate from $u8$ to $u9$.\medskip

\textbf{Step 4}\medskip

We next want to show $g(u9)\neq y$. But look at $u10$:\medskip

\begin{tabular}{c|c|c|c|c|c|c|c|}
\cline{2-8}
$u10$ & $1$ & $2$ & $3$ & $\cdots $ & $n-2$ & $n-1$ & $n$ \\ \cline{2-8}
& $%
\begin{array}{c}
x \\ 
z \\ 
y%
\end{array}%
$ & $%
\begin{array}{c}
y \\ 
z \\ 
x%
\end{array}%
$ & $%
\begin{array}{c}
y \\ 
z \\ 
x%
\end{array}%
$ & $\cdots $ & $%
\begin{array}{c}
y \\ 
z \\ 
x%
\end{array}%
$ & $%
\begin{array}{c}
z \\ 
y \\ 
x%
\end{array}%
$ & $%
\begin{array}{c}
z \\ 
y \\ 
x%
\end{array}%
$ \\ \cline{2-8}
\end{tabular}%
\medskip \newline
$g(u10)=x$ since $u10\in NP^{^{\ast }}$. But then at $u11$, an $n$-variant
of $u10$:\medskip

\begin{tabular}{c|c|c|c|c|c|c|c|}
\cline{2-8}
$u11$ & $1$ & $2$ & $3$ & $\cdots $ & $n-2$ & $n-1$ & $n$ \\ \cline{2-8}
& $%
\begin{array}{c}
x \\ 
z \\ 
y%
\end{array}%
$ & $%
\begin{array}{c}
y \\ 
z \\ 
x%
\end{array}%
$ & $%
\begin{array}{c}
y \\ 
z \\ 
x%
\end{array}%
$ & $\cdots $ & $%
\begin{array}{c}
y \\ 
z \\ 
x%
\end{array}%
$ & $%
\begin{array}{c}
z \\ 
y \\ 
x%
\end{array}%
$ & $%
\begin{array}{c}
x \\ 
y \\ 
z%
\end{array}%
$ \\ \cline{2-8}
\end{tabular}%
\medskip \newline
we get $g(u11)=x$ or $n$ manipulates from $u10$ to $u11$ (or $u11$ to $u10)$%
. Now $u9$ is a 1-variant of $u11$. If $g(u9)=y$, individual 1 would
manipulate from $u9$ to $u11$; so $g(u9)\neq y$. Combining with Step 3, $%
g(u9)=x$.\medskip

\textbf{Final Thread}\medskip

Consider profile\medskip

\begin{tabular}{c|c|c|c|c|c|c|c|}
\cline{2-8}
$u12$ & $1$ & $2$ & $3$ & $\cdots $ & $n-2$ & $n-1$ & $n$ \\ \cline{2-8}
& $%
\begin{array}{c}
z \\ 
x \\ 
y%
\end{array}%
$ & $%
\begin{array}{c}
y \\ 
z \\ 
x%
\end{array}%
$ & $%
\begin{array}{c}
y \\ 
z \\ 
x%
\end{array}%
$ & $\cdots $ & $%
\begin{array}{c}
y \\ 
z \\ 
x%
\end{array}%
$ & $%
\begin{array}{c}
z \\ 
y \\ 
x%
\end{array}%
$ & $%
\begin{array}{c}
y \\ 
x \\ 
z%
\end{array}%
$ \\ \cline{2-8}
\end{tabular}%
\medskip \newline
$u12$ is an $n$-variant of $u9$. If $g(u12)=z$, then $n$ would manipulate
from $u12$ to $u9$. So $g(u12)\neq z$.\medskip

Next consider profile $u13$:\medskip

\begin{tabular}{c|c|c|c|c|c|c|c|}
\cline{2-8}
$u13$ & $1$ & $2$ & $3$ & $\cdots $ & $n-2$ & $n-1$ & $n$ \\ \cline{2-8}
& $%
\begin{array}{c}
z \\ 
x \\ 
y%
\end{array}%
$ & $%
\begin{array}{c}
y \\ 
z \\ 
x%
\end{array}%
$ & $%
\begin{array}{c}
y \\ 
z \\ 
x%
\end{array}%
$ & $\cdots $ & $%
\begin{array}{c}
y \\ 
z \\ 
x%
\end{array}%
$ & $%
\begin{array}{c}
y \\ 
x \\ 
z%
\end{array}%
$ & $%
\begin{array}{c}
y \\ 
x \\ 
z%
\end{array}%
$ \\ \cline{2-8}
\end{tabular}%
\medskip \newline
Here $g(u13)=x$ since $u13$ is in $NP^{^{\ast }}$. Profile $u13$ is an $%
(n-1) $-variant of $u12$. If $g(u12)=y$, individual $n-1$ would manipulate
from $u13$ to $u12$. So $g(u12)\neq y$. Therefore $g(u12)=x$. Then look at $%
L1^{\ast \ast }$, a 2-variant $u12$:\medskip

\begin{tabular}{c|c|c|c|c|c|c|c|}
\cline{2-8}
$L1^{\ast \ast }$ & $1$ & $2$ & $3$ & $\cdots $ & $n-2$ & $n-1$ & $n$ \\ 
\cline{2-8}
& $%
\begin{array}{c}
z \\ 
x \\ 
y%
\end{array}%
$ & $%
\begin{array}{c}
z \\ 
x \\ 
y%
\end{array}%
$ & $%
\begin{array}{c}
y \\ 
z \\ 
x%
\end{array}%
$ & $\cdots $ & $%
\begin{array}{c}
y \\ 
z \\ 
x%
\end{array}%
$ & $%
\begin{array}{c}
z \\ 
y \\ 
x%
\end{array}%
$ & $%
\begin{array}{c}
y \\ 
x \\ 
z%
\end{array}%
$ \\ \cline{2-8}
\end{tabular}%
\medskip \newline
Then $g(L1^{\ast \ast })=x$ or 2 would manipulate from $u12$ to $L1^{\ast
\ast }$.\medskip

\textbf{Subsection 4-2-2}. Proof that $g(L3^{\ast \ast })=x$.\medskip

\textbf{Step 1}\medskip

We start at profile\medskip

\begin{tabular}{c|c|c|c|c|c|c|c|}
\cline{2-8}
$u1$ & $1$ & $2$ & $3$ & $\cdots $ & $n-2$ & $n-1$ & $n$ \\ \cline{2-8}
& $%
\begin{array}{c}
z \\ 
x \\ 
y%
\end{array}%
$ & $%
\begin{array}{c}
y \\ 
x \\ 
z%
\end{array}%
$ & $%
\begin{array}{c}
y \\ 
z \\ 
x%
\end{array}%
$ & $\cdots $ & $%
\begin{array}{c}
y \\ 
z \\ 
x%
\end{array}%
$ & $%
\begin{array}{c}
z \\ 
x \\ 
y%
\end{array}%
$ & $%
\begin{array}{c}
z \\ 
x \\ 
y%
\end{array}%
$ \\ \cline{2-8}
\end{tabular}%
\medskip \newline
where $g(u1)=x$ since $u1\in NP^{^{\ast }}$. Then, at $u2$, an $(n-1)$%
-variant of $u1$:\medskip

\begin{tabular}{c|c|c|c|c|c|c|c|}
\cline{2-8}
$u2$ & $1$ & $2$ & $3$ & $\cdots $ & $n-2$ & $n-1$ & $n$ \\ \cline{2-8}
& $%
\begin{array}{c}
z \\ 
x \\ 
y%
\end{array}%
$ & $%
\begin{array}{c}
y \\ 
x \\ 
z%
\end{array}%
$ & $%
\begin{array}{c}
y \\ 
z \\ 
x%
\end{array}%
$ & $\cdots $ & $%
\begin{array}{c}
y \\ 
z \\ 
x%
\end{array}%
$ & $%
\begin{array}{c}
y \\ 
z \\ 
x%
\end{array}%
$ & $%
\begin{array}{c}
z \\ 
x \\ 
y%
\end{array}%
$ \\ \cline{2-8}
\end{tabular}%
\medskip \newline
we see $g(u2)\neq z$ or $n-1$ manipulates from $u1$ to $u2$. Next consider $%
u3$, a 1-variant of $u2$:\medskip

\begin{tabular}{c|c|c|c|c|c|c|c|}
\cline{2-8}
$u3$ & $1$ & $2$ & $3$ & $\cdots $ & $n-2$ & $n-1$ & $n$ \\ \cline{2-8}
& $%
\begin{array}{c}
z \\ 
y \\ 
x%
\end{array}%
$ & $%
\begin{array}{c}
y \\ 
x \\ 
z%
\end{array}%
$ & $%
\begin{array}{c}
y \\ 
z \\ 
x%
\end{array}%
$ & $\cdots $ & $%
\begin{array}{c}
y \\ 
z \\ 
x%
\end{array}%
$ & $%
\begin{array}{c}
y \\ 
z \\ 
x%
\end{array}%
$ & $%
\begin{array}{c}
z \\ 
x \\ 
y%
\end{array}%
$ \\ \cline{2-8}
\end{tabular}%
\medskip \newline
If $g(u3)=z$, then 1 manipulates from $u2$ to $u3$, so $g(u3)\neq z$. But $%
u3 $ is also a 1-variant of $L1$:\medskip

\begin{tabular}{c|c|c|c|c|c|c|c|}
\cline{2-8}
$L1$ & $1$ & $2$ & $3$ & $\cdots $ & $n-2$ & $n-1$ & $n$ \\ \cline{2-8}
& $%
\begin{array}{c}
y \\ 
x \\ 
z%
\end{array}%
$ & $%
\begin{array}{c}
y \\ 
x \\ 
z%
\end{array}%
$ & $%
\begin{array}{c}
y \\ 
z \\ 
x%
\end{array}%
$ & $\cdots $ & $%
\begin{array}{c}
y \\ 
z \\ 
x%
\end{array}%
$ & $%
\begin{array}{c}
y \\ 
z \\ 
x%
\end{array}%
$ & $%
\begin{array}{c}
z \\ 
x \\ 
y%
\end{array}%
$ \\ \cline{2-8}
\end{tabular}%
\medskip \newline
$g(L1)=y$ by assumption. If $g(u3)=x$, then 1 would manipulate from $u3$ to $%
L1$. So $g(u3)\neq x$. Combining, $g(u3)=y$.\medskip

\qquad Next consider\medskip

\begin{tabular}{c|c|c|c|c|c|c|c|}
\cline{2-8}
$u4$ & $1$ & $2$ & $3$ & $\cdots $ & $n-2$ & $n-1$ & $n$ \\ \cline{2-8}
& $%
\begin{array}{c}
z \\ 
y \\ 
x%
\end{array}%
$ & $%
\begin{array}{c}
y \\ 
x \\ 
z%
\end{array}%
$ & $%
\begin{array}{c}
y \\ 
z \\ 
x%
\end{array}%
$ & $\cdots $ & $%
\begin{array}{c}
y \\ 
z \\ 
x%
\end{array}%
$ & $%
\begin{array}{c}
y \\ 
x \\ 
z%
\end{array}%
$ & $%
\begin{array}{c}
z \\ 
x \\ 
y%
\end{array}%
$ \\ \cline{2-8}
\end{tabular}%
\medskip \newline
This is an $(n-1)$-variant of $u3$ and $g(u4)=y$ or else $n-1$ would
manipulate from $u4$ to $u3$. Then look at $u5$, a 2-variant of $u4$:\medskip

\begin{tabular}{c|c|c|c|c|c|c|c|}
\cline{2-8}
$u5$ & $1$ & $2$ & $3$ & $\cdots $ & $n-2$ & $n-1$ & $n$ \\ \cline{2-8}
& $%
\begin{array}{c}
z \\ 
y \\ 
x%
\end{array}%
$ & $%
\begin{array}{c}
y \\ 
z \\ 
x%
\end{array}%
$ & $%
\begin{array}{c}
y \\ 
z \\ 
x%
\end{array}%
$ & $\cdots $ & $%
\begin{array}{c}
y \\ 
z \\ 
x%
\end{array}%
$ & $%
\begin{array}{c}
y \\ 
x \\ 
z%
\end{array}%
$ & $%
\begin{array}{c}
z \\ 
x \\ 
y%
\end{array}%
$ \\ \cline{2-8}
\end{tabular}%
\medskip \newline
We have $g(u5)=y$ or 2 would manipulate from $u5$ to $u4$. For the last part
of this step, consider $u6$, an $(n-1)$-variant of $u5$:\medskip

\begin{tabular}{c|c|c|c|c|c|c|c|}
\cline{2-8}
$u6$ & $1$ & $2$ & $3$ & $\cdots $ & $n-2$ & $n-1$ & $n$ \\ \cline{2-8}
& $%
\begin{array}{c}
z \\ 
y \\ 
x%
\end{array}%
$ & $%
\begin{array}{c}
y \\ 
z \\ 
x%
\end{array}%
$ & $%
\begin{array}{c}
y \\ 
z \\ 
x%
\end{array}%
$ & $\cdots $ & $%
\begin{array}{c}
y \\ 
z \\ 
x%
\end{array}%
$ & $%
\begin{array}{c}
x \\ 
y \\ 
z%
\end{array}%
$ & $%
\begin{array}{c}
z \\ 
x \\ 
y%
\end{array}%
$ \\ \cline{2-8}
\end{tabular}%
\medskip \newline
Then $g(u6)\neq z$ or $n-1$ manipulates from $u6$ to $u5$.\medskip

\textbf{Step 2}

Look at profile $u7$:\medskip

\begin{tabular}{c|c|c|c|c|c|c|c|}
\cline{2-8}
$u7$ & $1$ & $2$ & $3$ & $\cdots $ & $n-2$ & $n-1$ & $n$ \\ \cline{2-8}
& $%
\begin{array}{c}
y \\ 
x \\ 
z%
\end{array}%
$ & $%
\begin{array}{c}
y \\ 
z \\ 
x%
\end{array}%
$ & $%
\begin{array}{c}
y \\ 
z \\ 
x%
\end{array}%
$ & $\cdots $ & $%
\begin{array}{c}
y \\ 
z \\ 
x%
\end{array}%
$ & $%
\begin{array}{c}
z \\ 
x \\ 
y%
\end{array}%
$ & $%
\begin{array}{c}
z \\ 
x \\ 
y%
\end{array}%
$ \\ \cline{2-8}
\end{tabular}%
\medskip \newline
$g(u7)=x$ since $u7\in NP^{^{\ast }}$. Then look at $u8$, an $(n-1)$-variant
of $u7$:\medskip

\begin{tabular}{c|c|c|c|c|c|c|c|}
\cline{2-8}
$u8$ & $1$ & $2$ & $3$ & $\cdots $ & $n-2$ & $n-1$ & $n$ \\ \cline{2-8}
& $%
\begin{array}{c}
y \\ 
x \\ 
z%
\end{array}%
$ & $%
\begin{array}{c}
y \\ 
z \\ 
x%
\end{array}%
$ & $%
\begin{array}{c}
y \\ 
z \\ 
x%
\end{array}%
$ & $\cdots $ & $%
\begin{array}{c}
y \\ 
z \\ 
x%
\end{array}%
$ & $%
\begin{array}{c}
x \\ 
y \\ 
z%
\end{array}%
$ & $%
\begin{array}{c}
z \\ 
x \\ 
y%
\end{array}%
$ \\ \cline{2-8}
\end{tabular}%
\medskip \newline
Then $g(u8)=x$ or $n-1$ manipulates from $u8$ to $u7$. But $u8$ is also a
1-variant of $u6$. If $g(u6)=y$, individual 1 would manipulate from $u8$ to $%
u6$. So $g(u6)\neq y$.\medskip

\textbf{Final Thread}\medskip

Combining Steps 1 and 2, $g(u6)=x$. But $L3^{\ast \ast }$:\medskip

\begin{tabular}{c|c|c|c|c|c|c|c|}
\cline{2-8}
$L3^{\ast \ast }$ & $1$ & $2$ & $3$ & $\cdots $ & $n-2$ & $n-1$ & $n$ \\ 
\cline{2-8}
& $%
\begin{array}{c}
z \\ 
y \\ 
x%
\end{array}%
$ & $%
\begin{array}{c}
z \\ 
x \\ 
y%
\end{array}%
$ & $%
\begin{array}{c}
y \\ 
z \\ 
x%
\end{array}%
$ & $\cdots $ & $%
\begin{array}{c}
y \\ 
z \\ 
x%
\end{array}%
$ & $%
\begin{array}{c}
x \\ 
y \\ 
z%
\end{array}%
$ & $%
\begin{array}{c}
z \\ 
x \\ 
y%
\end{array}%
$ \\ \cline{2-8}
\end{tabular}%
\medskip \newline
is a 2-variant of $u6$, so $g(L3^{\ast \ast })=x$ or 2 manipulates from $u6$
to $L3^{\ast \ast }$.\medskip

\textbf{Subsection 4-2}-\textbf{3}. Proof that $g(L4^{\ast \ast })=x$%
.\medskip

\textbf{Step 1}\medskip

We start from\medskip

\begin{tabular}{c|c|c|c|c|c|c|c|}
\cline{2-8}
$u1$ & $1$ & $2$ & $3$ & $\cdots $ & $n-2$ & $n-1$ & $n$ \\ \cline{2-8}
& $%
\begin{array}{c}
z \\ 
x \\ 
y%
\end{array}%
$ & $%
\begin{array}{c}
y \\ 
x \\ 
z%
\end{array}%
$ & $%
\begin{array}{c}
y \\ 
z \\ 
x%
\end{array}%
$ & $\cdots $ & $%
\begin{array}{c}
y \\ 
z \\ 
x%
\end{array}%
$ & $%
\begin{array}{c}
y \\ 
z \\ 
x%
\end{array}%
$ & $%
\begin{array}{c}
y \\ 
z \\ 
x%
\end{array}%
$ \\ \cline{2-8}
\end{tabular}%
\medskip \newline
where $g(u1)=x$ since $u1\in NP^{^{\ast }}$. Then at $(n-1)$-variant $u2$%
:\medskip

\begin{tabular}{c|c|c|c|c|c|c|c|}
\cline{2-8}
$u2$ & $1$ & $2$ & $3$ & $\cdots $ & $n-2$ & $n-1$ & $n$ \\ \cline{2-8}
& $%
\begin{array}{c}
z \\ 
x \\ 
y%
\end{array}%
$ & $%
\begin{array}{c}
y \\ 
x \\ 
z%
\end{array}%
$ & $%
\begin{array}{c}
y \\ 
z \\ 
x%
\end{array}%
$ & $\cdots $ & $%
\begin{array}{c}
y \\ 
z \\ 
x%
\end{array}%
$ & $%
\begin{array}{c}
x \\ 
z \\ 
y%
\end{array}%
$ & $%
\begin{array}{c}
y \\ 
z \\ 
x%
\end{array}%
$ \\ \cline{2-8}
\end{tabular}%
\medskip \newline
we also have $g(u2)=x$ or $n-1$ manipulates from $u1$ to $u2$ or from $u2$
to $u1$. Then at 2-variant profile $u3$:\medskip

\begin{tabular}{c|c|c|c|c|c|c|c|}
\cline{2-8}
$u3$ & $1$ & $2$ & $3$ & $\cdots $ & $n-2$ & $n-1$ & $n$ \\ \cline{2-8}
& $%
\begin{array}{c}
z \\ 
x \\ 
y%
\end{array}%
$ & $%
\begin{array}{c}
y \\ 
z \\ 
x%
\end{array}%
$ & $%
\begin{array}{c}
y \\ 
z \\ 
x%
\end{array}%
$ & $\cdots $ & $%
\begin{array}{c}
y \\ 
z \\ 
x%
\end{array}%
$ & $%
\begin{array}{c}
x \\ 
z \\ 
y%
\end{array}%
$ & $%
\begin{array}{c}
y \\ 
z \\ 
x%
\end{array}%
$ \\ \cline{2-8}
\end{tabular}%
\medskip \newline
we see that $g(u3)\neq y$ or 2 would manipulate from $u2$ to $u3$.\medskip

\textbf{Step 2}\medskip

We start again from a profile in $NP^{^{\ast }}$:\medskip

\begin{tabular}{c|c|c|c|c|c|c|c|}
\cline{2-8}
$u4$ & $1$ & $2$ & $3$ & $\cdots $ & $n-2$ & $n-1$ & $n$ \\ \cline{2-8}
& $%
\begin{array}{c}
y \\ 
x \\ 
z%
\end{array}%
$ & $%
\begin{array}{c}
y \\ 
z \\ 
x%
\end{array}%
$ & $%
\begin{array}{c}
y \\ 
z \\ 
x%
\end{array}%
$ & $\cdots $ & $%
\begin{array}{c}
y \\ 
z \\ 
x%
\end{array}%
$ & $%
\begin{array}{c}
z \\ 
x \\ 
y%
\end{array}%
$ & $%
\begin{array}{c}
z \\ 
x \\ 
y%
\end{array}%
$ \\ \cline{2-8}
\end{tabular}%
\medskip \newline
so $g(u4)=x$. Then at profile $u5$, an $(n-1)$-variant of $u4$:\medskip

\begin{tabular}{c|c|c|c|c|c|c|c|}
\cline{2-8}
$u5$ & $1$ & $2$ & $3$ & $\cdots $ & $n-2$ & $n-1$ & $n$ \\ \cline{2-8}
& $%
\begin{array}{c}
y \\ 
x \\ 
z%
\end{array}%
$ & $%
\begin{array}{c}
y \\ 
z \\ 
x%
\end{array}%
$ & $%
\begin{array}{c}
y \\ 
z \\ 
x%
\end{array}%
$ & $\cdots $ & $%
\begin{array}{c}
y \\ 
z \\ 
x%
\end{array}%
$ & $%
\begin{array}{c}
x \\ 
z \\ 
y%
\end{array}%
$ & $%
\begin{array}{c}
z \\ 
x \\ 
y%
\end{array}%
$ \\ \cline{2-8}
\end{tabular}%
\medskip \newline
we see $g(u5)=x$ or $n-1$ manipulates from $u5$ to $u4$. Next examine $u6$,
a 1-variant of $u5$:\medskip

\begin{tabular}{c|c|c|c|c|c|c|c|}
\cline{2-8}
$u6$ & $1$ & $2$ & $3$ & $\cdots $ & $n-2$ & $n-1$ & $n$ \\ \cline{2-8}
& $%
\begin{array}{c}
y \\ 
z \\ 
x%
\end{array}%
$ & $%
\begin{array}{c}
y \\ 
z \\ 
x%
\end{array}%
$ & $%
\begin{array}{c}
y \\ 
z \\ 
x%
\end{array}%
$ & $\cdots $ & $%
\begin{array}{c}
y \\ 
z \\ 
x%
\end{array}%
$ & $%
\begin{array}{c}
x \\ 
z \\ 
y%
\end{array}%
$ & $%
\begin{array}{c}
z \\ 
x \\ 
y%
\end{array}%
$ \\ \cline{2-8}
\end{tabular}%
\medskip \newline
We have $g(u6)\neq y$ or 1 manipulates from $u5$ to $u6$. But we saw earlier
that $g(L3^{\ast \ast })=x$ so, by Lemma 4-5, $g(L2^{\ast \ast })=x$\medskip

\begin{tabular}{c|c|c|c|c|c|c|c|}
\cline{2-8}
$L2^{^{\ast \ast }}$ & $1$ & $2$ & $3$ & $\cdots $ & $n-2$ & $n-1$ & $n$ \\ 
\cline{2-8}
& $%
\begin{array}{c}
y \\ 
z \\ 
x%
\end{array}%
$ & $%
\begin{array}{c}
z \\ 
x \\ 
y%
\end{array}%
$ & $%
\begin{array}{c}
y \\ 
z \\ 
x%
\end{array}%
$ & $\cdots $ & $%
\begin{array}{c}
y \\ 
z \\ 
x%
\end{array}%
$ & $%
\begin{array}{c}
x \\ 
z \\ 
y%
\end{array}%
$ & $%
\begin{array}{c}
z \\ 
x \\ 
y%
\end{array}%
$ \\ \cline{2-8}
\end{tabular}%
\medskip \newline

Profile $L2^{\ast \ast }$ is a 2-variant of $u6$; if $g(u6)=z$, then 2 would
manipulate from $L2^{\ast \ast }$ to $u6$. So $g(u6)\neq z$. Combining, $%
g(u6)=x$.\medskip

Now profile $u6$ is also a 1-variant of profile $u7$:\medskip

\begin{tabular}{c|c|c|c|c|c|c|c|}
\cline{2-8}
$u7$ & $1$ & $2$ & $3$ & $\cdots $ & $n-2$ & $n-1$ & $n$ \\ \cline{2-8}
& $%
\begin{array}{c}
z \\ 
x \\ 
y%
\end{array}%
$ & $%
\begin{array}{c}
y \\ 
z \\ 
x%
\end{array}%
$ & $%
\begin{array}{c}
y \\ 
z \\ 
x%
\end{array}%
$ & $\cdots $ & $%
\begin{array}{c}
y \\ 
z \\ 
x%
\end{array}%
$ & $%
\begin{array}{c}
x \\ 
z \\ 
y%
\end{array}%
$ & $%
\begin{array}{c}
z \\ 
x \\ 
y%
\end{array}%
$ \\ \cline{2-8}
\end{tabular}%
\medskip\ \newline
So $g(u7)=x$ or 1 manipulates from $u6$ to $u7$. Profile $u7$ is also an $n$%
-variant of $u3$. If $g(u3)=z$, $n$ would manipulate from $u7$ to $u3$.
Therefore, $g(u3)\neq z$.\medskip

\textbf{Final Thread}\medskip

Combining Steps 1 and 2, we get $g(u3)=x$. But $u3$ is a 2-variant of $%
L4^{\ast \ast }$:\medskip `

\begin{tabular}{c|c|c|c|c|c|c|c|}
\cline{2-8}
$L4^{\ast \ast }$ & $1$ & $2$ & $3$ & $\cdots $ & $n-2$ & $n-1$ & $n$ \\ 
\cline{2-8}
& $%
\begin{array}{c}
z \\ 
x \\ 
y%
\end{array}%
$ & $%
\begin{array}{c}
z \\ 
x \\ 
y%
\end{array}%
$ & $%
\begin{array}{c}
y \\ 
z \\ 
x%
\end{array}%
$ & $\cdots $ & $%
\begin{array}{c}
y \\ 
z \\ 
x%
\end{array}%
$ & $%
\begin{array}{c}
x \\ 
z \\ 
y%
\end{array}%
$ & $%
\begin{array}{c}
y \\ 
z \\ 
x%
\end{array}%
$ \\ \cline{2-8}
\end{tabular}%
\medskip \newline
If $g(L4^{\ast \ast })\neq x$, then 2 would manipulate from $u3$ to $%
L4^{\ast }$; therefore, $g(L4^{\ast \ast })=x$.\medskip

\textbf{Subsection 4-2-4. Other profiles.}\medskip

We have shown that if Range$(g|NP^{^{\ast }})=\{x\}$ and $g(L1)=y$ then $z$
is not chosen at any of the profiles in the list $L^{^{\ast \ast }}$. But $%
L^{^{\ast \ast }}$ was constructed by making some arbitrary choices about
profiles after having made some choices to get list $L$. Those arbitrary
choices in the construction of $L^{^{\ast \ast }}$ then may not satisfy a
\textquotedblleft without loss of generality\textquotedblright\ argument. We
need to consider what happens with other choices. Of course we would like to
treat all possibilities, but there are very many and each is dealt with
fairly straightforwardly.\medskip

\qquad For example, construct profile $T1^{\ast }$ that differs from $%
L3^{\ast \ast }$ by\medskip

\qquad (I) interchanging the preference orderings for individual 1 and 2;
and also\medskip

\qquad (II) interchanging the preference orderings for individual $n-1$ and $%
n$.\medskip

\begin{tabular}{c|c|c|c|c|c|c|c|}
\cline{2-8}
$T1^{^{\ast }}$ & $1$ & $2$ & $3$ & $\cdots $ & $n-2$ & $n-1$ & $n$ \\ 
\cline{2-8}
& $%
\begin{array}{c}
z \\ 
x \\ 
y%
\end{array}%
$ & $%
\begin{array}{c}
z \\ 
y \\ 
x%
\end{array}%
$ & $%
\begin{array}{c}
y \\ 
z \\ 
x%
\end{array}%
$ & $\cdots $ & $%
\begin{array}{c}
y \\ 
z \\ 
x%
\end{array}%
$ & $%
\begin{array}{c}
z \\ 
x \\ 
y%
\end{array}%
$ & $%
\begin{array}{c}
x \\ 
y \\ 
z%
\end{array}%
$ \\ \cline{2-8}
\end{tabular}%
\medskip \newline
We will show $g(T1^{\ast })\neq z$; in fact, we can show $g(T1^{\ast })=x$%
.\medskip

\qquad Earlier, in our analysis of $L1^{\ast \ast }$, we showed $x$ is
chosen at the following profile (which was called $u12$ back there):\medskip

\begin{tabular}{c|c|c|c|c|c|c|c|}
\cline{2-8}
$u1$ & $1$ & $2$ & $3$ & $\cdots $ & $n-2$ & $n-1$ & $n$ \\ \cline{2-8}
& $%
\begin{array}{c}
z \\ 
x \\ 
y%
\end{array}%
$ & $%
\begin{array}{c}
y \\ 
z \\ 
x%
\end{array}%
$ & $%
\begin{array}{c}
y \\ 
z \\ 
x%
\end{array}%
$ & $\cdots $ & $%
\begin{array}{c}
y \\ 
z \\ 
x%
\end{array}%
$ & $%
\begin{array}{c}
z \\ 
y \\ 
x%
\end{array}%
$ & $%
\begin{array}{c}
y \\ 
x \\ 
z%
\end{array}%
$ \\ \cline{2-8}
\end{tabular}%
\medskip \newline
Then at 2-variant profile $u2$:\medskip

\begin{tabular}{c|c|c|c|c|c|c|c|}
\cline{2-8}
$u2$ & $1$ & $2$ & $3$ & $\cdots $ & $n-2$ & $n-1$ & $n$ \\ \cline{2-8}
& $%
\begin{array}{c}
z \\ 
x \\ 
y%
\end{array}%
$ & $%
\begin{array}{c}
z \\ 
y \\ 
x%
\end{array}%
$ & $%
\begin{array}{c}
y \\ 
z \\ 
x%
\end{array}%
$ & $\cdots $ & $%
\begin{array}{c}
y \\ 
z \\ 
x%
\end{array}%
$ & $%
\begin{array}{c}
z \\ 
y \\ 
x%
\end{array}%
$ & $%
\begin{array}{c}
y \\ 
x \\ 
z%
\end{array}%
$ \\ \cline{2-8}
\end{tabular}%
\medskip \newline
we get $g(u2)=x$ or 2 manipulates from $u1$ to $u2$. Finally consider $u3$,
an $(n-1)$-variant of $u2$:\medskip

\begin{tabular}{c|c|c|c|c|c|c|c|}
\cline{2-8}
$u3$ & $1$ & $2$ & $3$ & $\cdots $ & $n-2$ & $n-1$ & $n$ \\ \cline{2-8}
& $%
\begin{array}{c}
z \\ 
x \\ 
y%
\end{array}%
$ & $%
\begin{array}{c}
z \\ 
y \\ 
x%
\end{array}%
$ & $%
\begin{array}{c}
y \\ 
z \\ 
x%
\end{array}%
$ & $\cdots $ & $%
\begin{array}{c}
y \\ 
z \\ 
x%
\end{array}%
$ & $%
\begin{array}{c}
z \\ 
x \\ 
y%
\end{array}%
$ & $%
\begin{array}{c}
y \\ 
x \\ 
z%
\end{array}%
$ \\ \cline{2-8}
\end{tabular}%
\medskip \newline
Then $g(u3)=x$ or $n-1$ manipulates from $u2$ to $u3$. But $u3$ is an $n$%
-variant of $T1^{\ast }$. If $g(T1^{\ast })\neq x$, $n$ would manipulate
from $T1^{\ast }$ to $u3$. Therefore $g(T1^{\ast })=x$.\medskip

\textbf{Section 4-3}. Assume $g(L2)=y$.\medskip

\textbf{Subsection 4-3}-\textbf{1}. Proof that $g(L1^{\ast \ast })=x$%
.\medskip

At profile $u1$:\medskip

\begin{tabular}{c|c|c|c|c|c|c|c|}
\cline{2-8}
$u1$ & $1$ & $2$ & $3$ & $\cdots $ & $n-2$ & $n-1$ & $n$ \\ \cline{2-8}
& $%
\begin{array}{c}
z \\ 
y \\ 
x%
\end{array}%
$ & $%
\begin{array}{c}
y \\ 
x \\ 
z%
\end{array}%
$ & $%
\begin{array}{c}
y \\ 
z \\ 
x%
\end{array}%
$ & $\cdots $ & $%
\begin{array}{c}
y \\ 
z \\ 
x%
\end{array}%
$ & $%
\begin{array}{c}
z \\ 
x \\ 
y%
\end{array}%
$ & $%
\begin{array}{c}
z \\ 
x \\ 
y%
\end{array}%
$ \\ \cline{2-8}
\end{tabular}%
\medskip \newline
we have $g(u1)=x$ since $u1\in NP^{^{\ast }}$. Consider $n$-variant profile $%
u2$\medskip

\begin{tabular}{c|c|c|c|c|c|c|c|}
\cline{2-8}
$u2$ & $1$ & $2$ & $3$ & $\cdots $ & $n-2$ & $n-1$ & $n$ \\ \cline{2-8}
& $%
\begin{array}{c}
z \\ 
y \\ 
x%
\end{array}%
$ & $%
\begin{array}{c}
y \\ 
x \\ 
z%
\end{array}%
$ & $%
\begin{array}{c}
y \\ 
z \\ 
x%
\end{array}%
$ & $\cdots $ & $%
\begin{array}{c}
y \\ 
z \\ 
x%
\end{array}%
$ & $%
\begin{array}{c}
z \\ 
x \\ 
y%
\end{array}%
$ & $%
\begin{array}{c}
y \\ 
x \\ 
z%
\end{array}%
$ \\ \cline{2-8}
\end{tabular}%
\medskip \newline
If $g(u2)=z$, then $n$ would manipulate from $u1$ to $u2$ (and $u_{2}$ to $%
u_{1}$). So $g(u2)\neq z$. Profile $u2$ is also an $(n-1)$-variant of $L2$%
:\medskip

\begin{tabular}{c|c|c|c|c|c|c|c|}
\cline{2-8}
$L2$ & $1$ & $2$ & $3$ & $\cdots $ & $n-2$ & $n-1$ & $n$ \\ \cline{2-8}
& $%
\begin{array}{c}
z \\ 
y \\ 
x%
\end{array}%
$ & $%
\begin{array}{c}
y \\ 
x \\ 
z%
\end{array}%
$ & $%
\begin{array}{c}
y \\ 
z \\ 
x%
\end{array}%
$ & $\cdots $ & $%
\begin{array}{c}
y \\ 
z \\ 
x%
\end{array}%
$ & $%
\begin{array}{c}
x \\ 
y \\ 
z%
\end{array}%
$ & $%
\begin{array}{c}
y \\ 
x \\ 
z%
\end{array}%
$ \\ \cline{2-8}
\end{tabular}%
\medskip \newline
and $g(L2$)$=y$ by the assumption of this section. If $g(u2)=x$, then $n-1$
manipulates from $L2$ to $u2$. Therefore $g(u2)\neq x$. Combining, $g(u2)=y$%
.\medskip

Now look at $u3$ a 2-variant of $u2$:\medskip

\begin{tabular}{c|c|c|c|c|c|c|c|}
\cline{2-8}
$u3$ & $1$ & $2$ & $3$ & $\cdots $ & $n-2$ & $n-1$ & $n$ \\ \cline{2-8}
& $%
\begin{array}{c}
z \\ 
y \\ 
x%
\end{array}%
$ & $%
\begin{array}{c}
y \\ 
z \\ 
x%
\end{array}%
$ & $%
\begin{array}{c}
y \\ 
z \\ 
x%
\end{array}%
$ & $\cdots $ & $%
\begin{array}{c}
y \\ 
z \\ 
x%
\end{array}%
$ & $%
\begin{array}{c}
z \\ 
x \\ 
y%
\end{array}%
$ & $%
\begin{array}{c}
y \\ 
x \\ 
z%
\end{array}%
$ \\ \cline{2-8}
\end{tabular}%
\medskip \newline
Then $g(u3)=y$ or 2 would manipulate from $u3$ to $u2$.\medskip

Next, 1-variant profile $u4$:\medskip

\begin{tabular}{c|c|c|c|c|c|c|c|}
\cline{2-8}
$u4$ & $1$ & $2$ & $3$ & $\cdots $ & $n-2$ & $n-1$ & $n$ \\ \cline{2-8}
& $%
\begin{array}{c}
z \\ 
x \\ 
y%
\end{array}%
$ & $%
\begin{array}{c}
y \\ 
z \\ 
x%
\end{array}%
$ & $%
\begin{array}{c}
y \\ 
z \\ 
x%
\end{array}%
$ & $\cdots $ & $%
\begin{array}{c}
y \\ 
z \\ 
x%
\end{array}%
$ & $%
\begin{array}{c}
z \\ 
x \\ 
y%
\end{array}%
$ & $%
\begin{array}{c}
y \\ 
x \\ 
z%
\end{array}%
$ \\ \cline{2-8}
\end{tabular}%
\medskip \newline
If $g(u4)=z$, then 1 manipulates from $u3$ to $u4$. So $g(u4)\neq z$. Next,
look at $u5$, an $(n-1)$-variant of $u4$:\medskip

\begin{tabular}{c|c|c|c|c|c|c|c|}
\cline{2-8}
$u5$ & $1$ & $2$ & $3$ & $\cdots $ & $n-2$ & $n-1$ & $n$ \\ \cline{2-8}
& $%
\begin{array}{c}
z \\ 
x \\ 
y%
\end{array}%
$ & $%
\begin{array}{c}
y \\ 
z \\ 
x%
\end{array}%
$ & $%
\begin{array}{c}
y \\ 
z \\ 
x%
\end{array}%
$ & $\cdots $ & $%
\begin{array}{c}
y \\ 
z \\ 
x%
\end{array}%
$ & $%
\begin{array}{c}
y \\ 
x \\ 
z%
\end{array}%
$ & $%
\begin{array}{c}
y \\ 
x \\ 
z%
\end{array}%
$ \\ \cline{2-8}
\end{tabular}%
\medskip \newline
$g(u5)=x$ since $u5\in NP^{^{\ast }}$. If $g(u4)=y$, then $n-1$ manipulates
from $u5$ to $u4$. Therefore, $g(u4)\neq y$. Combining, $g(u4)=x$.\medskip

Consider profile $u6$, also an $(n-1)$-variant of $u4$:\medskip

\begin{tabular}{c|c|c|c|c|c|c|c|}
\cline{2-8}
$u6$ & $1$ & $2$ & $3$ & $\cdots $ & \multicolumn{1}{|c|}{$n-2$} & $n-1$ & $%
n $ \\ \cline{2-8}
& $%
\begin{array}{c}
z \\ 
x \\ 
y%
\end{array}%
$ & $%
\begin{array}{c}
y \\ 
z \\ 
x%
\end{array}%
$ & $%
\begin{array}{c}
y \\ 
z \\ 
x%
\end{array}%
$ & $\cdots $ & \multicolumn{1}{|c|}{$%
\begin{array}{c}
y \\ 
z \\ 
x%
\end{array}%
$} & $%
\begin{array}{c}
z \\ 
y \\ 
x%
\end{array}%
$ & $%
\begin{array}{c}
y \\ 
x \\ 
z%
\end{array}%
$ \\ \cline{2-8}
\end{tabular}%
\medskip \newline
If $g(u6)=z$, then $n-1$ would manipulate from $u4$ to $u6$. Also $u6$ is an 
$(n-1)$-variant of $u5$. If $g(u6)=y$, then $n-1$ would manipulate from $u5$
to $u6$. Combining, $g(u6)=x$.\medskip

But $u6$ is a 2-variant of $L1^{\ast \ast }$:\medskip

\begin{tabular}{c|c|c|c|c|c|c|c|}
\cline{2-8}
$L1^{\ast \ast }$ & $1$ & $2$ & $3$ & $\cdots $ & \multicolumn{1}{|c|}{$n-2$}
& $n-1$ & $n$ \\ \cline{2-8}
& $%
\begin{array}{c}
z \\ 
x \\ 
y%
\end{array}%
$ & $%
\begin{array}{c}
z \\ 
x \\ 
y%
\end{array}%
$ & $%
\begin{array}{c}
y \\ 
z \\ 
x%
\end{array}%
$ & $\cdots $ & \multicolumn{1}{|c|}{$%
\begin{array}{c}
y \\ 
z \\ 
x%
\end{array}%
$} & $%
\begin{array}{c}
z \\ 
y \\ 
x%
\end{array}%
$ & $%
\begin{array}{c}
y \\ 
x \\ 
z%
\end{array}%
$ \\ \cline{2-8}
\end{tabular}%
\medskip \newline
$g(L1^{\ast \ast })=x$ or 2 would manipulate from $u6$ to $L1^{\ast \ast }$%
.\medskip

\textbf{Subsection 4-3}-\textbf{2}. Proof that $g(L3^{\ast \ast })=x$%
.\medskip

\textbf{Step 1}\medskip

We start with profile $u1$:\medskip

\begin{tabular}{c|c|c|c|c|c|c|c|}
\cline{2-8}
$u1$ & $1$ & $2$ & $3$ & $\cdots $ & $n-2$ & $n-1$ & $n$ \\ \cline{2-8}
& $%
\begin{array}{c}
z \\ 
y \\ 
x%
\end{array}%
$ & $%
\begin{array}{c}
x \\ 
z \\ 
y%
\end{array}%
$ & $%
\begin{array}{c}
y \\ 
z \\ 
x%
\end{array}%
$ & $\cdots $ & $%
\begin{array}{c}
y \\ 
z \\ 
x%
\end{array}%
$ & $%
\begin{array}{c}
y \\ 
z \\ 
x%
\end{array}%
$ & $%
\begin{array}{c}
y \\ 
z \\ 
x%
\end{array}%
$ \\ \cline{2-8}
\end{tabular}%
\medskip \newline
where $g(u1)=x$ since $u1\in NP^{^{\ast }}$. Then look at $(n-1)$-variant $%
u2 $:\medskip

\begin{tabular}{c|c|c|c|c|c|c|c|}
\cline{2-8}
$u2$ & $1$ & $2$ & $3$ & $\cdots $ & \multicolumn{1}{|c|}{$n-2$} & $n-1$ & $%
n $ \\ \cline{2-8}
& $%
\begin{array}{c}
z \\ 
y \\ 
x%
\end{array}%
$ & $%
\begin{array}{c}
x \\ 
z \\ 
y%
\end{array}%
$ & $%
\begin{array}{c}
y \\ 
z \\ 
x%
\end{array}%
$ & $\cdots $ & \multicolumn{1}{|c|}{$%
\begin{array}{c}
y \\ 
z \\ 
x%
\end{array}%
$} & $%
\begin{array}{c}
x \\ 
y \\ 
z%
\end{array}%
$ & $%
\begin{array}{c}
y \\ 
z \\ 
x%
\end{array}%
$ \\ \cline{2-8}
\end{tabular}%
\medskip \newline
$g(u2)=x$ or $n-1$ manipulates from $u1$ to $u2$. Next consider 2-variant
profile $u3$:\medskip

\begin{tabular}{c|c|c|c|c|c|c|c|}
\cline{2-8}
$u3$ & $1$ & $2$ & $3$ & $\cdots $ & \multicolumn{1}{|c|}{$n-2$} & $n-1$ & $%
n $ \\ \cline{2-8}
& $%
\begin{array}{c}
z \\ 
y \\ 
x%
\end{array}%
$ & $%
\begin{array}{c}
z \\ 
x \\ 
y%
\end{array}%
$ & $%
\begin{array}{c}
y \\ 
z \\ 
x%
\end{array}%
$ & $\cdots $ & \multicolumn{1}{|c|}{$%
\begin{array}{c}
y \\ 
z \\ 
x%
\end{array}%
$} & $%
\begin{array}{c}
x \\ 
y \\ 
z%
\end{array}%
$ & $%
\begin{array}{c}
y \\ 
z \\ 
x%
\end{array}%
$ \\ \cline{2-8}
\end{tabular}%
\medskip \newline
If $g(u3)=y$, then 2 manipulates from $u3$ to $u2$. So $g(u3)\neq y$.\medskip

\textbf{Step 2}\medskip

By assumption, at $L2$:\medskip

\begin{tabular}{c|c|c|c|c|c|c|c|}
\cline{2-8}
$L2$ & $1$ & $2$ & $3$ & $\cdots $ & \multicolumn{1}{|c|}{$n-2$} & $n-1$ & $%
n $ \\ \cline{2-8}
& $%
\begin{array}{c}
z \\ 
y \\ 
x%
\end{array}%
$ & $%
\begin{array}{c}
y \\ 
x \\ 
z%
\end{array}%
$ & $%
\begin{array}{c}
y \\ 
z \\ 
x%
\end{array}%
$ & $\cdots $ & \multicolumn{1}{|c|}{$%
\begin{array}{c}
y \\ 
z \\ 
x%
\end{array}%
$} & $%
\begin{array}{c}
x \\ 
y \\ 
z%
\end{array}%
$ & $%
\begin{array}{c}
y \\ 
x \\ 
z%
\end{array}%
$ \\ \cline{2-8}
\end{tabular}%
\medskip

we have $g(L2)=y$. Then at 2-variant profile $u4$:\medskip

\begin{tabular}{c|c|c|c|c|c|c|c|}
\cline{2-8}
$u4$ & $1$ & $2$ & $3$ & $\cdots $ & \multicolumn{1}{|c|}{$n-2$} & $n-1$ & $%
n $ \\ \cline{2-8}
& $%
\begin{array}{c}
z \\ 
y \\ 
x%
\end{array}%
$ & $%
\begin{array}{c}
z \\ 
y \\ 
x%
\end{array}%
$ & $%
\begin{array}{c}
y \\ 
z \\ 
x%
\end{array}%
$ & $\cdots $ & \multicolumn{1}{|c|}{$%
\begin{array}{c}
y \\ 
z \\ 
x%
\end{array}%
$} & $%
\begin{array}{c}
x \\ 
y \\ 
z%
\end{array}%
$ & $%
\begin{array}{c}
y \\ 
x \\ 
z%
\end{array}%
$ \\ \cline{2-8}
\end{tabular}%
\medskip \newline
we have $g(u4)\neq x$ or 2 manipulates from $u4$ to $L2$. Next, at profile $%
u5$:\medskip

\begin{tabular}{c|c|c|c|c|c|c|c|}
\cline{2-8}
$u5$ & $1$ & $2$ & $3$ & $\cdots $ & \multicolumn{1}{|c|}{$n-2$} & $n-1$ & $%
n $ \\ \cline{2-8}
& $%
\begin{array}{c}
z \\ 
y \\ 
x%
\end{array}%
$ & $%
\begin{array}{c}
z \\ 
x \\ 
y%
\end{array}%
$ & $%
\begin{array}{c}
y \\ 
z \\ 
x%
\end{array}%
$ & $\cdots $ & \multicolumn{1}{|c|}{$%
\begin{array}{c}
y \\ 
z \\ 
x%
\end{array}%
$} & $%
\begin{array}{c}
y \\ 
x \\ 
z%
\end{array}%
$ & $%
\begin{array}{c}
y \\ 
x \\ 
z%
\end{array}%
$ \\ \cline{2-8}
\end{tabular}%
\medskip \newline
we have $g(u5)=x$ since $u5\in NP^{^{\ast }}$. Then at $(n-1)$%
-variant\medskip

\begin{tabular}{c|c|c|c|c|c|c|c|}
\cline{2-8}
$u6$ & $1$ & $2$ & $3$ & $\cdots $ & \multicolumn{1}{|c|}{$n-2$} & $n-1$ & $%
n $ \\ \cline{2-8}
& $%
\begin{array}{c}
z \\ 
y \\ 
x%
\end{array}%
$ & $%
\begin{array}{c}
z \\ 
x \\ 
y%
\end{array}%
$ & $%
\begin{array}{c}
y \\ 
z \\ 
x%
\end{array}%
$ & $\cdots $ & \multicolumn{1}{|c|}{$%
\begin{array}{c}
y \\ 
z \\ 
x%
\end{array}%
$} & $%
\begin{array}{c}
x \\ 
y \\ 
z%
\end{array}%
$ & $%
\begin{array}{c}
y \\ 
x \\ 
z%
\end{array}%
$ \\ \cline{2-8}
\end{tabular}%
\medskip \newline
we get $g(u6)=x$ or $n-1$ manipulates from $u6$ to $u5$. But $u6$ is a
2-variant of $u4$ and if $g(u4)=z$, then 2 would manipulate from $u6$ to $u4$%
. So $g(u4)\neq z$. Combining, $g(u4)=y$.\medskip

\textbf{Final thread}\medskip

Profile $u7$ is an $n$-variant of $u4$:\medskip

\begin{tabular}{c|c|c|c|c|c|c|c|}
\cline{2-8}
$u7$ & $1$ & $2$ & $3$ & $\cdots $ & \multicolumn{1}{|c|}{$n-2$} & $n-1$ & $%
n $ \\ \cline{2-8}
& $%
\begin{array}{c}
z \\ 
y \\ 
x%
\end{array}%
$ & $%
\begin{array}{c}
z \\ 
y \\ 
x%
\end{array}%
$ & $%
\begin{array}{c}
y \\ 
z \\ 
x%
\end{array}%
$ & $\cdots $ & \multicolumn{1}{|c|}{$%
\begin{array}{c}
y \\ 
z \\ 
x%
\end{array}%
$} & $%
\begin{array}{c}
x \\ 
y \\ 
z%
\end{array}%
$ & $%
\begin{array}{c}
y \\ 
z \\ 
x%
\end{array}%
$ \\ \cline{2-8}
\end{tabular}%
\medskip \newline
Since $g(u4)=y$, we must have $g(u7)=y$. But $u7$ is also a 2-variant of $u3$%
. If $g(u3)=z$, then 2 would manipulate from $u7$ to $u3$. Combining, $%
g(u3)=x$.\medskip

Profile $L3^{\ast \ast }$ is\medskip

\begin{tabular}{c|c|c|c|c|c|c|c|}
\cline{2-8}
$L3^{\ast \ast }$ & $1$ & $2$ & $3$ & $\cdots $ & \multicolumn{1}{|c|}{$n-2$}
& $n-1$ & $n$ \\ \cline{2-8}
& $%
\begin{array}{c}
z \\ 
y \\ 
x%
\end{array}%
$ & $%
\begin{array}{c}
z \\ 
x \\ 
y%
\end{array}%
$ & $%
\begin{array}{c}
y \\ 
z \\ 
x%
\end{array}%
$ & $\cdots $ & \multicolumn{1}{|c|}{$%
\begin{array}{c}
y \\ 
z \\ 
x%
\end{array}%
$} & $%
\begin{array}{c}
x \\ 
y \\ 
z%
\end{array}%
$ & $%
\begin{array}{c}
z \\ 
x \\ 
y%
\end{array}%
$ \\ \cline{2-8}
\end{tabular}%
\medskip \newline
and so is an $n$-variant of $u3$. Then $g(L3^{\ast \ast })=x$ or $n$ would
manipulate from $u3$ to $L3^{\ast \ast }$.\medskip

\textbf{Subsection 4-3}-\textbf{3}. Proof that $g(L4^{\ast \ast }$)$=x$%
.\medskip

\textbf{Step 1}\medskip

At $u1$:\medskip

\begin{tabular}{c|c|c|c|c|c|c|c|}
\cline{2-8}
$u1$ & $1$ & $2$ & $3$ & $\cdots $ & \multicolumn{1}{|c|}{$n-2$} & $n-1$ & $%
n $ \\ \cline{2-8}
& $%
\begin{array}{c}
z \\ 
y \\ 
x%
\end{array}%
$ & $%
\begin{array}{c}
z \\ 
y \\ 
x%
\end{array}%
$ & $%
\begin{array}{c}
y \\ 
z \\ 
x%
\end{array}%
$ & $\cdots $ & \multicolumn{1}{|c|}{$%
\begin{array}{c}
y \\ 
z \\ 
x%
\end{array}%
$} & $%
\begin{array}{c}
x \\ 
y \\ 
z%
\end{array}%
$ & $%
\begin{array}{c}
x \\ 
y \\ 
z%
\end{array}%
$ \\ \cline{2-8}
\end{tabular}%
\medskip \newline
$g(u1)=x$ since $u1\in NP^{^{\ast }}$. Profile $u2$ is an $(n-1)$-variant of 
$u1$:\medskip

\begin{tabular}{c|c|c|c|c|c|c|c|}
\cline{2-8}
$u2$ & $1$ & $2$ & $3$ & $\cdots $ & \multicolumn{1}{|c|}{$n-2$} & $n-1$ & $%
n $ \\ \cline{2-8}
& $%
\begin{array}{c}
z \\ 
y \\ 
x%
\end{array}%
$ & $%
\begin{array}{c}
z \\ 
y \\ 
x%
\end{array}%
$ & $%
\begin{array}{c}
y \\ 
z \\ 
x%
\end{array}%
$ & $\cdots $ & \multicolumn{1}{|c|}{$%
\begin{array}{c}
y \\ 
z \\ 
x%
\end{array}%
$} & $%
\begin{array}{c}
x \\ 
z \\ 
y%
\end{array}%
$ & $%
\begin{array}{c}
x \\ 
y \\ 
z%
\end{array}%
$ \\ \cline{2-8}
\end{tabular}%
\medskip \newline
$g(u2)=x$ or $n-1$ manipulates from $u2$ to $u1$. Then look at $n$-variant
profile $u3$:\medskip

\begin{tabular}{c|c|c|c|c|c|c|c|}
\cline{2-8}
$u3$ & $1$ & $2$ & $3$ & $\cdots $ & \multicolumn{1}{|c|}{$n-2$} & $n-1$ & $%
n $ \\ \cline{2-8}
& $%
\begin{array}{c}
z \\ 
y \\ 
x%
\end{array}%
$ & $%
\begin{array}{c}
z \\ 
y \\ 
x%
\end{array}%
$ & $%
\begin{array}{c}
y \\ 
z \\ 
x%
\end{array}%
$ & $\cdots $ & \multicolumn{1}{|c|}{$%
\begin{array}{c}
y \\ 
z \\ 
x%
\end{array}%
$} & $%
\begin{array}{c}
x \\ 
z \\ 
y%
\end{array}%
$ & $%
\begin{array}{c}
y \\ 
x \\ 
z%
\end{array}%
$ \\ \cline{2-8}
\end{tabular}%
\medskip \newline
If $g(u3)=z$, then $n$ would manipulate from $u3$ to $u2$. So $g(u3)\neq z$%
.\medskip

\textbf{Step 2}\medskip

Now we are assuming that $y$ is chosen at $L2$:\medskip

\begin{tabular}{c|c|c|c|c|c|c|c|}
\cline{2-8}
$L2$ & $1$ & $2$ & $3$ & $\cdots $ & \multicolumn{1}{|c|}{$n-2$} & $n-1$ & $%
n $ \\ \cline{2-8}
& $%
\begin{array}{c}
z \\ 
y \\ 
x%
\end{array}%
$ & $%
\begin{array}{c}
y \\ 
x \\ 
z%
\end{array}%
$ & $%
\begin{array}{c}
y \\ 
z \\ 
x%
\end{array}%
$ & $\cdots $ & \multicolumn{1}{|c|}{$%
\begin{array}{c}
y \\ 
z \\ 
x%
\end{array}%
$} & $%
\begin{array}{c}
x \\ 
y \\ 
z%
\end{array}%
$ & $%
\begin{array}{c}
y \\ 
x \\ 
z%
\end{array}%
$ \\ \cline{2-8}
\end{tabular}%
\medskip \newline
So at 2-variant profile $u4$:\medskip

\begin{tabular}{c|c|c|c|c|c|c|c|}
\cline{2-8}
$u4$ & $1$ & $2$ & $3$ & $\cdots $ & \multicolumn{1}{|c|}{$n-2$} & $n-1$ & $%
n $ \\ \cline{2-8}
& $%
\begin{array}{c}
z \\ 
y \\ 
x%
\end{array}%
$ & $%
\begin{array}{c}
z \\ 
y \\ 
x%
\end{array}%
$ & $%
\begin{array}{c}
y \\ 
z \\ 
x%
\end{array}%
$ & $\cdots $ & \multicolumn{1}{|c|}{$%
\begin{array}{c}
y \\ 
z \\ 
x%
\end{array}%
$} & $%
\begin{array}{c}
x \\ 
y \\ 
z%
\end{array}%
$ & $%
\begin{array}{c}
y \\ 
x \\ 
z%
\end{array}%
$ \\ \cline{2-8}
\end{tabular}%
\medskip

$g(u4)\neq x$ or 2 manipulates from $u4$ to $L2$. But $u3$ is also an $(n-1)$%
-variant of $u4$ and if $g(u3)=x$, $n-1$ would manipulate from $u4$ to $u3$.
So $g(u3)\neq x$.\medskip

\textbf{Step 3}\medskip

Combining Steps 1 and 2, $g(u3)=y$. So at $u5$, an $n$-variant of $u3$%
:\medskip

\begin{tabular}{c|c|c|c|c|c|c|c|}
\cline{2-8}
$u5$ & $1$ & $2$ & $3$ & $\cdots $ & \multicolumn{1}{|c|}{$n-2$} & $n-1$ & $%
n $ \\ \cline{2-8}
& $%
\begin{array}{c}
z \\ 
y \\ 
x%
\end{array}%
$ & $%
\begin{array}{c}
z \\ 
y \\ 
x%
\end{array}%
$ & $%
\begin{array}{c}
y \\ 
z \\ 
x%
\end{array}%
$ & $\cdots $ & \multicolumn{1}{|c|}{$%
\begin{array}{c}
y \\ 
z \\ 
x%
\end{array}%
$} & $%
\begin{array}{c}
x \\ 
z \\ 
y%
\end{array}%
$ & $%
\begin{array}{c}
y \\ 
z \\ 
x%
\end{array}%
$ \\ \cline{2-8}
\end{tabular}%
\medskip \newline
$g(u5)=y$ or $n$ manipulates from $u5$ to $u3$.\medskip

\textbf{Step 4}\medskip

Profile $u6$ is in $NP^{^{\ast }}$:\medskip

\begin{tabular}{c|c|c|c|c|c|c|c|}
\cline{2-8}
$u6$ & $1$ & $2$ & $3$ & $\cdots $ & \multicolumn{1}{|c|}{$n-2$} & $n-1$ & $%
n $ \\ \cline{2-8}
& $%
\begin{array}{c}
z \\ 
y \\ 
x%
\end{array}%
$ & $%
\begin{array}{c}
x \\ 
z \\ 
x%
\end{array}%
$ & $%
\begin{array}{c}
y \\ 
z \\ 
x%
\end{array}%
$ & $\cdots $ & \multicolumn{1}{|c|}{$%
\begin{array}{c}
y \\ 
z \\ 
x%
\end{array}%
$} & $%
\begin{array}{c}
y \\ 
z \\ 
x%
\end{array}%
$ & $%
\begin{array}{c}
y \\ 
z \\ 
x%
\end{array}%
$ \\ \cline{2-8}
\end{tabular}%
\medskip \newline
so $g(u6)=x$. Then at $(n-1)$-variant profile $u7$:\medskip

\begin{tabular}{c|c|c|c|c|c|c|c|}
\cline{2-8}
$u7$ & $1$ & $2$ & $3$ & $\cdots $ & \multicolumn{1}{|c|}{$n-2$} & $n-1$ & $%
n $ \\ \cline{2-8}
& $%
\begin{array}{c}
z \\ 
y \\ 
x%
\end{array}%
$ & $%
\begin{array}{c}
x \\ 
z \\ 
y%
\end{array}%
$ & $%
\begin{array}{c}
y \\ 
z \\ 
x%
\end{array}%
$ & $\cdots $ & \multicolumn{1}{|c|}{$%
\begin{array}{c}
y \\ 
z \\ 
x%
\end{array}%
$} & $%
\begin{array}{c}
x \\ 
z \\ 
y%
\end{array}%
$ & $%
\begin{array}{c}
y \\ 
z \\ 
x%
\end{array}%
$ \\ \cline{2-8}
\end{tabular}%
\medskip \newline
$g(u7)=x$ or $n-1$ would manipulate from $u7$ to $u6$.\medskip

\textbf{Final thread}\medskip

Consider profile $u8$:\medskip

\begin{tabular}{c|c|c|c|c|c|c|c|}
\cline{2-8}
$u8$ & $1$ & $2$ & $3$ & $\cdots $ & \multicolumn{1}{|c|}{$n-2$} & $n-1$ & $%
n $ \\ \cline{2-8}
& $%
\begin{array}{c}
z \\ 
y \\ 
x%
\end{array}%
$ & $%
\begin{array}{c}
z \\ 
x \\ 
y%
\end{array}%
$ & $%
\begin{array}{c}
y \\ 
z \\ 
x%
\end{array}%
$ & $\cdots $ & \multicolumn{1}{|c|}{$%
\begin{array}{c}
y \\ 
z \\ 
x%
\end{array}%
$} & $%
\begin{array}{c}
x \\ 
z \\ 
y%
\end{array}%
$ & $%
\begin{array}{c}
y \\ 
z \\ 
x%
\end{array}%
$ \\ \cline{2-8}
\end{tabular}%
\medskip \newline
Profile $u8$ is a 2-variant of $u5$ and also a 2-variant of $u7$. By Step 3, 
$g(u5)=y$, so $g(u8)\neq z$ or 2 manipulates from $u5$ to $u8$. By Step 4, $%
g(u7)=x$, so $g(u8)\neq y$, or 2 manipulates from $u8$ to $u7$. Hence $%
g(u8)=x$. But $u8$ is a 1-variant of $L4^{\ast \ast }$:\medskip

\begin{tabular}{c|c|c|c|c|c|c|c|}
\cline{2-8}
$L4^{\ast \ast }$ & $1$ & $2$ & $3$ & $\cdots $ & \multicolumn{1}{|c|}{$n-2$}
& $n-1$ & $n$ \\ \cline{2-8}
& $%
\begin{array}{c}
z \\ 
x \\ 
y%
\end{array}%
$ & $%
\begin{array}{c}
z \\ 
x \\ 
y%
\end{array}%
$ & $%
\begin{array}{c}
y \\ 
z \\ 
x%
\end{array}%
$ & $\cdots $ & \multicolumn{1}{|c|}{$%
\begin{array}{c}
y \\ 
z \\ 
x%
\end{array}%
$} & $%
\begin{array}{c}
x \\ 
z \\ 
y%
\end{array}%
$ & $%
\begin{array}{c}
y \\ 
z \\ 
x%
\end{array}%
$ \\ \cline{2-8}
\end{tabular}%
\medskip \newline
so $g(L4^{\ast \ast })=x$ or 1 would manipulate from $u8$ to $L4^{\ast \ast
} $.\medskip

\textbf{Section 4-4}. Assume $g(L3)=y$.\medskip

\textbf{Subsection 4-4}-\textbf{1}. Proof that $g(L1^{\ast \ast })=x$%
.\medskip

At $L3$:\medskip

\begin{tabular}{c|c|c|c|c|c|c|c|}
\cline{2-8}
$L3$ & $1$ & $2$ & $3$ & $\cdots $ & \multicolumn{1}{|c|}{$n-2$} & $n-1$ & $%
n $ \\ \cline{2-8}
& $%
\begin{array}{c}
y \\ 
z \\ 
x%
\end{array}%
$ & $%
\begin{array}{c}
y \\ 
x \\ 
z%
\end{array}%
$ & $%
\begin{array}{c}
y \\ 
z \\ 
x%
\end{array}%
$ & $\cdots $ & \multicolumn{1}{|c|}{$%
\begin{array}{c}
y \\ 
z \\ 
x%
\end{array}%
$} & $%
\begin{array}{c}
x \\ 
z \\ 
y%
\end{array}%
$ & $%
\begin{array}{c}
y \\ 
x \\ 
z%
\end{array}%
$ \\ \cline{2-8}
\end{tabular}%
\medskip \newline
we have $g(L3)=y$ by assumption. Then at 1-variant profile $u1$\medskip

\begin{tabular}{c|c|c|c|c|c|c|c|}
\cline{2-8}
$u1$ & $1$ & $2$ & $3$ & $\cdots $ & \multicolumn{1}{|c|}{$n-2$} & $n-1$ & $%
n $ \\ \cline{2-8}
& $%
\begin{array}{c}
z \\ 
y \\ 
x%
\end{array}%
$ & $%
\begin{array}{c}
y \\ 
x \\ 
z%
\end{array}%
$ & $%
\begin{array}{c}
y \\ 
z \\ 
x%
\end{array}%
$ & $\cdots $ & \multicolumn{1}{|c|}{$%
\begin{array}{c}
y \\ 
z \\ 
x%
\end{array}%
$} & $%
\begin{array}{c}
x \\ 
z \\ 
y%
\end{array}%
$ & $%
\begin{array}{c}
y \\ 
x \\ 
z%
\end{array}%
$ \\ \cline{2-8}
\end{tabular}%
\medskip \newline
we get $g(u1)\neq x$ or 1 would manipulate from $u1$ to $L3$. Next, at $u2$,
an $n$-variant of $u1$:\medskip

\begin{tabular}{c|c|c|c|c|c|c|c|}
\cline{2-8}
$u2$ & $1$ & $2$ & $3$ & $\cdots $ & \multicolumn{1}{|c|}{$n-2$} & $n-1$ & $%
n $ \\ \cline{2-8}
& $%
\begin{array}{c}
z \\ 
y \\ 
x%
\end{array}%
$ & $%
\begin{array}{c}
y \\ 
x \\ 
z%
\end{array}%
$ & $%
\begin{array}{c}
y \\ 
z \\ 
x%
\end{array}%
$ & $\cdots $ & \multicolumn{1}{|c|}{$%
\begin{array}{c}
y \\ 
z \\ 
x%
\end{array}%
$} & $%
\begin{array}{c}
x \\ 
z \\ 
y%
\end{array}%
$ & $%
\begin{array}{c}
x \\ 
z \\ 
y%
\end{array}%
$ \\ \cline{2-8}
\end{tabular}%
\medskip \newline
$g(u2)=x$ since $u2\in NP^{^{\ast }}$. If $g(u1)=z$, then $n$ would
manipulate from $u1$ to $u2$. So $g(u1)\neq z$. Combining, $g(u1)=y$.\medskip

The next profile, $u3$, is an $(n-1)$-variant of $u1$:\medskip

\begin{tabular}{c|c|c|c|c|c|c|c|}
\cline{2-8}
$u3$ & $1$ & $2$ & $3$ & $\cdots $ & $n-2$ & $n-1$ & $n$ \\ \cline{2-8}
& $%
\begin{array}{c}
z \\ 
y \\ 
x%
\end{array}%
$ & $%
\begin{array}{c}
y \\ 
x \\ 
z%
\end{array}%
$ & $%
\begin{array}{c}
y \\ 
z \\ 
x%
\end{array}%
$ & $\cdots $ & $%
\begin{array}{c}
y \\ 
z \\ 
x%
\end{array}%
$ & $%
\begin{array}{c}
z \\ 
x \\ 
y%
\end{array}%
$ & $%
\begin{array}{c}
y \\ 
x \\ 
z%
\end{array}%
$ \\ \cline{2-8}
\end{tabular}%
\medskip \newline
$g(u3)=y$ or $n-1$ would manipulate from $u1$ to $u3$.\medskip

\textbf{Step 2}\medskip

Now look at profile $u4$, a 2-variant of $u3$:\medskip

\begin{tabular}{c|c|c|c|c|c|c|c|}
\cline{2-8}
$u4$ & $1$ & $2$ & $3$ & $\cdots $ & \multicolumn{1}{|c|}{$n-2$} & $n-1$ & $%
n $ \\ \cline{2-8}
& $%
\begin{array}{c}
z \\ 
y \\ 
x%
\end{array}%
$ & $%
\begin{array}{c}
y \\ 
z \\ 
x%
\end{array}%
$ & $%
\begin{array}{c}
y \\ 
z \\ 
x%
\end{array}%
$ & $\cdots $ & \multicolumn{1}{|c|}{$%
\begin{array}{c}
y \\ 
z \\ 
x%
\end{array}%
$} & $%
\begin{array}{c}
z \\ 
x \\ 
y%
\end{array}%
$ & $%
\begin{array}{c}
y \\ 
x \\ 
z%
\end{array}%
$ \\ \cline{2-8}
\end{tabular}%
\medskip \newline
Then $g(u4)=y$ or 2 manipulates from $u4$ to $u3$. Next, look at 1-variant $%
u5$:\medskip

\begin{tabular}{c|c|c|c|c|c|c|c|}
\cline{2-8}
$u5$ & $1$ & $2$ & $3$ & $\cdots $ & \multicolumn{1}{|c|}{$n-2$} & $n-1$ & $%
n $ \\ \cline{2-8}
& $%
\begin{array}{c}
z \\ 
x \\ 
y%
\end{array}%
$ & $%
\begin{array}{c}
y \\ 
z \\ 
x%
\end{array}%
$ & $%
\begin{array}{c}
y \\ 
z \\ 
x%
\end{array}%
$ & $\cdots $ & \multicolumn{1}{|c|}{$%
\begin{array}{c}
y \\ 
z \\ 
x%
\end{array}%
$} & $%
\begin{array}{c}
z \\ 
x \\ 
y%
\end{array}%
$ & $%
\begin{array}{c}
y \\ 
x \\ 
z%
\end{array}%
$ \\ \cline{2-8}
\end{tabular}%
\medskip \newline
If $g(u5)=z$, 1 would manipulate from $u4$ to $u5$. So $g(u5)\neq z$. Then
look at $u6$, an $(n-1)$-variant of $u5$:\medskip

\begin{tabular}{c|c|c|c|c|c|c|c|}
\cline{2-8}
$u6$ & $1$ & $2$ & $3$ & $\cdots $ & \multicolumn{1}{|c|}{$n-2$} & $n-1$ & $%
n $ \\ \cline{2-8}
& $%
\begin{array}{c}
z \\ 
x \\ 
y%
\end{array}%
$ & $%
\begin{array}{c}
y \\ 
z \\ 
x%
\end{array}%
$ & $%
\begin{array}{c}
y \\ 
z \\ 
x%
\end{array}%
$ & $\cdots $ & \multicolumn{1}{|c|}{$%
\begin{array}{c}
y \\ 
z \\ 
x%
\end{array}%
$} & $%
\begin{array}{c}
y \\ 
x \\ 
z%
\end{array}%
$ & $%
\begin{array}{c}
y \\ 
x \\ 
z%
\end{array}%
$ \\ \cline{2-8}
\end{tabular}%
\medskip \newline
$g(u6)=x$ since $u6\in NP^{^{\ast }}$. If $g(u5)=y$, then $n-1$ manipulates
from $u6$ to $u5$. So $g(u5)\neq y$. Combining, $g(u5)=x$.\medskip

\textbf{Final thread}\medskip

Consider profile $u7$:\medskip

\begin{tabular}{c|c|c|c|c|c|c|c|}
\cline{2-8}
$u7$ & $1$ & $2$ & $3$ & $\cdots $ & \multicolumn{1}{|c|}{$n-2$} & $n-1$ & $%
n $ \\ \cline{2-8}
& $%
\begin{array}{c}
z \\ 
x \\ 
y%
\end{array}%
$ & $%
\begin{array}{c}
y \\ 
z \\ 
x%
\end{array}%
$ & $%
\begin{array}{c}
y \\ 
z \\ 
x%
\end{array}%
$ & $\cdots $ & \multicolumn{1}{|c|}{$%
\begin{array}{c}
y \\ 
z \\ 
x%
\end{array}%
$} & $%
\begin{array}{c}
z \\ 
y \\ 
x%
\end{array}%
$ & $%
\begin{array}{c}
y \\ 
x \\ 
z%
\end{array}%
$ \\ \cline{2-8}
\end{tabular}%
\medskip \newline
$u7$ is an $(n-1)$-variant of $u6$ and of $u5$. If $g(u7)=y$, then $n-1$
manipulates from $u6$ to $u7$; if $g(u7)=z$, then $n-1$ manipulates from $u5$
to $u7$. So $g(u7)=x$.\medskip

Now $L1^{\ast \ast }$ is given by:\medskip

\begin{tabular}{c|c|c|c|c|c|c|c|}
\cline{2-8}
$L1^{\ast \ast }$ & $1$ & $2$ & $3$ & $\cdots $ & \multicolumn{1}{|c|}{$n-2$}
& $n-1$ & $n$ \\ \cline{2-8}
& $%
\begin{array}{c}
z \\ 
x \\ 
y%
\end{array}%
$ & $%
\begin{array}{c}
z \\ 
x \\ 
y%
\end{array}%
$ & $%
\begin{array}{c}
y \\ 
z \\ 
x%
\end{array}%
$ & $\cdots $ & \multicolumn{1}{|c|}{$%
\begin{array}{c}
y \\ 
z \\ 
x%
\end{array}%
$} & $%
\begin{array}{c}
z \\ 
y \\ 
x%
\end{array}%
$ & $%
\begin{array}{c}
y \\ 
x \\ 
z%
\end{array}%
$ \\ \cline{2-8}
\end{tabular}%
\medskip \newline
so $g(L1^{\ast \ast })=x$ or else 2 would manipulate from $u7$ to $L1^{\ast
\ast }$.\medskip

\textbf{Subsection 4-4}-\textbf{2}. Proof that $g(L3^{\ast \ast })=x$%
.\medskip

\textbf{Step 1}\medskip

At profile $u1$ in $NP^{^{\ast }}$:\medskip

\begin{tabular}{c|c|c|c|c|c|c|c|}
\cline{2-8}
$u1$ & $1$ & $2$ & $3$ & $\cdots $ & \multicolumn{1}{|c|}{$n-2$} & $n-1$ & $%
n $ \\ \cline{2-8}
& $%
\begin{array}{c}
z \\ 
y \\ 
x%
\end{array}%
$ & $%
\begin{array}{c}
x \\ 
y \\ 
z%
\end{array}%
$ & $%
\begin{array}{c}
y \\ 
z \\ 
x%
\end{array}%
$ & $\cdots $ & \multicolumn{1}{|c|}{$%
\begin{array}{c}
y \\ 
z \\ 
x%
\end{array}%
$} & $%
\begin{array}{c}
y \\ 
z \\ 
x%
\end{array}%
$ & $%
\begin{array}{c}
y \\ 
z \\ 
x%
\end{array}%
$ \\ \cline{2-8}
\end{tabular}%
\medskip \newline
$g(u1)=x$. Then at $u2$, an $(n-1)$-variant of $u1$:\medskip

\begin{tabular}{c|c|c|c|c|c|c|c|}
\cline{2-8}
$u2$ & $1$ & $2$ & $3$ & $\cdots $ & \multicolumn{1}{|c|}{$n-2$} & $n-1$ & $%
n $ \\ \cline{2-8}
& $%
\begin{array}{c}
z \\ 
y \\ 
x%
\end{array}%
$ & $%
\begin{array}{c}
x \\ 
y \\ 
z%
\end{array}%
$ & $%
\begin{array}{c}
y \\ 
z \\ 
x%
\end{array}%
$ & $\cdots $ & \multicolumn{1}{|c|}{$%
\begin{array}{c}
y \\ 
z \\ 
x%
\end{array}%
$} & $%
\begin{array}{c}
x \\ 
y \\ 
z%
\end{array}%
$ & $%
\begin{array}{c}
y \\ 
z \\ 
x%
\end{array}%
$ \\ \cline{2-8}
\end{tabular}%
\medskip \newline
$g(u2)=x$ or $n-1$ manipulates from $u2$ to $u1$. Then consider $u3$, a
2-variant of $u2$:\medskip

\begin{tabular}{c|c|c|c|c|c|c|c|}
\cline{2-8}
$u3$ & $1$ & $2$ & $3$ & $\cdots $ & \multicolumn{1}{|c|}{$n-2$} & $n-1$ & $%
n $ \\ \cline{2-8}
& $%
\begin{array}{c}
z \\ 
y \\ 
x%
\end{array}%
$ & $%
\begin{array}{c}
z \\ 
x \\ 
y%
\end{array}%
$ & $%
\begin{array}{c}
y \\ 
z \\ 
x%
\end{array}%
$ & $\cdots $ & \multicolumn{1}{|c|}{$%
\begin{array}{c}
y \\ 
z \\ 
x%
\end{array}%
$} & $%
\begin{array}{c}
x \\ 
y \\ 
z%
\end{array}%
$ & $%
\begin{array}{c}
y \\ 
z \\ 
x%
\end{array}%
$ \\ \cline{2-8}
\end{tabular}%
\medskip \newline
If $g(u3)=y$, then 2 would manipulate from $u3$ to $u2$. So $g(u3)\neq y$%
.\medskip

\textbf{Step 2}\medskip

Now look at profile $u4$:\medskip

\begin{tabular}{c|c|c|c|c|c|c|c|}
\cline{2-8}
$u4$ & $1$ & $2$ & $3$ & $\cdots $ & \multicolumn{1}{|c|}{$n-2$} & $n-1$ & $%
n $ \\ \cline{2-8}
& $%
\begin{array}{c}
z \\ 
y \\ 
x%
\end{array}%
$ & $%
\begin{array}{c}
z \\ 
x \\ 
y%
\end{array}%
$ & $%
\begin{array}{c}
y \\ 
z \\ 
x%
\end{array}%
$ & $\cdots $ & \multicolumn{1}{|c|}{$%
\begin{array}{c}
y \\ 
z \\ 
x%
\end{array}%
$} & $%
\begin{array}{c}
y \\ 
x \\ 
z%
\end{array}%
$ & $%
\begin{array}{c}
y \\ 
x \\ 
z%
\end{array}%
$ \\ \cline{2-8}
\end{tabular}%
\medskip \newline
$g(u4)=x$ since $u4\in NP^{^{\ast }}$. Then consider $(n-1)$-variant profile 
$u5$:\medskip

\begin{tabular}{c|c|c|c|c|c|c|c|}
\cline{2-8}
$u5$ & $1$ & $2$ & $3$ & $\cdots $ & \multicolumn{1}{|c|}{$n-2$} & $n-1$ & $%
n $ \\ \cline{2-8}
& $%
\begin{array}{c}
z \\ 
y \\ 
x%
\end{array}%
$ & $%
\begin{array}{c}
z \\ 
x \\ 
y%
\end{array}%
$ & $%
\begin{array}{c}
y \\ 
z \\ 
x%
\end{array}%
$ & $\cdots $ & \multicolumn{1}{|c|}{$%
\begin{array}{c}
y \\ 
z \\ 
x%
\end{array}%
$} & $%
\begin{array}{c}
x \\ 
z \\ 
y%
\end{array}%
$ & $%
\begin{array}{c}
y \\ 
x \\ 
z%
\end{array}%
$ \\ \cline{2-8}
\end{tabular}%
\medskip \newline
$g(u5)=x$ or $n-1$ manipulates from $u5$ to $u4$. Next look at $u6$, a
2-variant of $u5$:\medskip

\begin{tabular}{c|c|c|c|c|c|c|c|}
\cline{2-8}
$u6$ & $1$ & $2$ & $3$ & $\cdots $ & \multicolumn{1}{|c|}{$n-2$} & $n-1$ & $%
n $ \\ \cline{2-8}
& $%
\begin{array}{c}
z \\ 
y \\ 
x%
\end{array}%
$ & $%
\begin{array}{c}
z \\ 
y \\ 
x%
\end{array}%
$ & $%
\begin{array}{c}
y \\ 
z \\ 
x%
\end{array}%
$ & $\cdots $ & \multicolumn{1}{|c|}{$%
\begin{array}{c}
y \\ 
z \\ 
x%
\end{array}%
$} & $%
\begin{array}{c}
x \\ 
z \\ 
y%
\end{array}%
$ & $%
\begin{array}{c}
y \\ 
x \\ 
z%
\end{array}%
$ \\ \cline{2-8}
\end{tabular}%
\medskip \newline
If $g(u6)=z$, then 2 manipulates from $u5$ to $u6$. So $g(u6)\neq z$. This
profile $u6$ is also a 2-variant of $u7$:\medskip

\begin{tabular}{c|c|c|c|c|c|c|c|}
\cline{2-8}
$u7$ & $1$ & $2$ & $3$ & $\cdots $ & \multicolumn{1}{|c|}{$n-2$} & $n-1$ & $%
n $ \\ \cline{2-8}
& $%
\begin{array}{c}
z \\ 
y \\ 
x%
\end{array}%
$ & $%
\begin{array}{c}
y \\ 
x \\ 
z%
\end{array}%
$ & $%
\begin{array}{c}
y \\ 
z \\ 
x%
\end{array}%
$ & $\cdots $ & \multicolumn{1}{|c|}{$%
\begin{array}{c}
y \\ 
z \\ 
x%
\end{array}%
$} & $%
\begin{array}{c}
x \\ 
z \\ 
y%
\end{array}%
$ & $%
\begin{array}{c}
y \\ 
x \\ 
z%
\end{array}%
$ \\ \cline{2-8}
\end{tabular}%
\medskip \newline
It was shown in Subsection 4-4-1 that $g(L3)=y$ implies $g(u7)=y$ (in that
subsection, this profile was called $u1$). If $g(u6)=x$, then 2 would
manipulate from $u6$ to $u7$. So $g(u6)\neq x$. Combining, $g(u6)=y$.\medskip

\textbf{Step 3}\medskip

Then at profile $u8$, an $(n-1)$-variant of $u6$:\medskip

\begin{tabular}{c|c|c|c|c|c|c|c|}
\cline{2-8}
$u8$ & $1$ & $2$ & $3$ & $\cdots $ & \multicolumn{1}{|c|}{$n-2$} & $n-1$ & $%
n $ \\ \cline{2-8}
& $%
\begin{array}{c}
z \\ 
y \\ 
x%
\end{array}%
$ & $%
\begin{array}{c}
z \\ 
y \\ 
x%
\end{array}%
$ & $%
\begin{array}{c}
y \\ 
z \\ 
x%
\end{array}%
$ & $\cdots $ & \multicolumn{1}{|c|}{$%
\begin{array}{c}
y \\ 
z \\ 
x%
\end{array}%
$} & $%
\begin{array}{c}
x \\ 
y \\ 
z%
\end{array}%
$ & $%
\begin{array}{c}
y \\ 
x \\ 
z%
\end{array}%
$ \\ \cline{2-8}
\end{tabular}%
\medskip \newline
we have $g(u8)=y$ or $n-1$ manipulates from $u6$ to $u8$. Now look at
profile $u9$, an $n$-variant of $u8$:\medskip

\begin{tabular}{c|c|c|c|c|c|c|c|}
\cline{2-8}
$u9$ & $1$ & $2$ & $3$ & $\cdots $ & \multicolumn{1}{|c|}{$n-2$} & $n-1$ & $%
n $ \\ \cline{2-8}
& $%
\begin{array}{c}
z \\ 
y \\ 
x%
\end{array}%
$ & $%
\begin{array}{c}
z \\ 
y \\ 
x%
\end{array}%
$ & $%
\begin{array}{c}
y \\ 
z \\ 
x%
\end{array}%
$ & $\cdots $ & \multicolumn{1}{|c|}{$%
\begin{array}{c}
y \\ 
z \\ 
x%
\end{array}%
$} & $%
\begin{array}{c}
x \\ 
y \\ 
z%
\end{array}%
$ & $%
\begin{array}{c}
y \\ 
z \\ 
x%
\end{array}%
$ \\ \cline{2-8}
\end{tabular}%
\medskip \newline
$g(u9)=y$ or $n$ manipulates from $u9$ to $u8$. But $u9$ is a 2-variant of $%
u3$. If $g(u3)=z$, 2 manipulates from $u9$ to $u3$. So $g(u3)\neq z$.\medskip

\textbf{Final thread}\medskip

From Step 1, we got $g(u3)\neq y$; from step 3, $g(u3)\neq z$. So $g(u3)=x$.
But $u3$ is an $n$-variant of $L3^{\ast \ast }$:\medskip

\begin{tabular}{c|c|c|c|c|c|c|c|}
\cline{2-8}
$L3^{\ast \ast }$ & $1$ & $2$ & $3$ & $\cdots $ & \multicolumn{1}{|c|}{$n-2$}
& $n-1$ & $n$ \\ \cline{2-8}
& $%
\begin{array}{c}
z \\ 
y \\ 
x%
\end{array}%
$ & $%
\begin{array}{c}
z \\ 
x \\ 
y%
\end{array}%
$ & $%
\begin{array}{c}
y \\ 
z \\ 
x%
\end{array}%
$ & $\cdots $ & \multicolumn{1}{|c|}{$%
\begin{array}{c}
y \\ 
z \\ 
x%
\end{array}%
$} & $%
\begin{array}{c}
x \\ 
y \\ 
z%
\end{array}%
$ & $%
\begin{array}{c}
z \\ 
x \\ 
y%
\end{array}%
$ \\ \cline{2-8}
\end{tabular}%
\medskip \newline
$g(L3^{\ast \ast })=x$ or $n$ would manipulate from $L3^{\ast \ast }$ to $u3$
(and $u3$ to $L3^{\ast \ast }$).\medskip

\textbf{Subsection 4-4}-\textbf{3}. Proof that $g(L4^{\ast \ast })=x$%
.\medskip

In the previous subsection, we found $x$ is chosen at profile $u1$ (called $%
u3$ there):\medskip

\begin{tabular}{c|c|c|c|c|c|c|c|}
\cline{2-8}
$u1$ & $1$ & $2$ & $3$ & $\cdots $ & \multicolumn{1}{|c|}{$n-2$} & $n-1$ & $%
n $ \\ \cline{2-8}
& $%
\begin{array}{c}
z \\ 
y \\ 
x%
\end{array}%
$ & $%
\begin{array}{c}
z \\ 
x \\ 
y%
\end{array}%
$ & $%
\begin{array}{c}
y \\ 
z \\ 
x%
\end{array}%
$ & $\cdots $ & \multicolumn{1}{|c|}{$%
\begin{array}{c}
y \\ 
z \\ 
x%
\end{array}%
$} & $%
\begin{array}{c}
x \\ 
y \\ 
z%
\end{array}%
$ & $%
\begin{array}{c}
y \\ 
z \\ 
x%
\end{array}%
$ \\ \cline{2-8}
\end{tabular}%
\medskip \newline
Consider 1-variant profile $u2$:\medskip

\begin{tabular}{c|c|c|c|c|c|c|c|}
\cline{2-8}
$u2$ & $1$ & $2$ & $3$ & $\cdots $ & \multicolumn{1}{|c|}{$n-2$} & $n-1$ & $%
n $ \\ \cline{2-8}
& $%
\begin{array}{c}
z \\ 
x \\ 
y%
\end{array}%
$ & $%
\begin{array}{c}
z \\ 
x \\ 
y%
\end{array}%
$ & $%
\begin{array}{c}
y \\ 
z \\ 
x%
\end{array}%
$ & $\cdots $ & \multicolumn{1}{|c|}{$%
\begin{array}{c}
y \\ 
z \\ 
x%
\end{array}%
$} & $%
\begin{array}{c}
x \\ 
y \\ 
z%
\end{array}%
$ & $%
\begin{array}{c}
y \\ 
z \\ 
x%
\end{array}%
$ \\ \cline{2-8}
\end{tabular}%
\medskip \newline
$g(u2)=x$ or 1 would manipulate from $u1$ to $u2$. But $u2$ is an $(n-1)$%
-variant of $L4^{\ast \ast }$:\medskip

\begin{tabular}{c|c|c|c|c|c|c|c|}
\cline{2-8}
$L4^{\ast \ast }$ & $1$ & $2$ & $3$ & $\cdots $ & \multicolumn{1}{|c|}{$n-2$}
& $n-1$ & $n$ \\ \cline{2-8}
& $%
\begin{array}{c}
z \\ 
x \\ 
y%
\end{array}%
$ & $%
\begin{array}{c}
z \\ 
x \\ 
y%
\end{array}%
$ & $%
\begin{array}{c}
y \\ 
z \\ 
x%
\end{array}%
$ & $\cdots $ & \multicolumn{1}{|c|}{$%
\begin{array}{c}
y \\ 
z \\ 
x%
\end{array}%
$} & $%
\begin{array}{c}
x \\ 
z \\ 
y%
\end{array}%
$ & $%
\begin{array}{c}
y \\ 
z \\ 
x%
\end{array}%
$ \\ \cline{2-8}
\end{tabular}%
\medskip \newline
So $g(L4^{\ast \ast })=x$ or $n-1$ would manipulate from $L4^{\ast \ast }$
to $u2$.\medskip

\textbf{Section 4-5.} Assume $g(L4)=y$.\medskip

\textbf{Subsection 4-5}-\textbf{1}. Proof that $g(L1^{\ast \ast })=x$%
.\medskip

\textbf{Step 1}\medskip

We start from\medskip

\begin{tabular}{c|c|c|c|c|c|c|c|}
\cline{2-8}
$L4$ & $1$ & $2$ & $3$ & $\cdots $ & $n-2$ & $n-1$ & $n$ \\ \cline{2-8}
& $%
\begin{array}{c}
y \\ 
x \\ 
z%
\end{array}%
$ & $%
\begin{array}{c}
y \\ 
x \\ 
z%
\end{array}%
$ & $%
\begin{array}{c}
y \\ 
z \\ 
x%
\end{array}%
$ & $\cdots $ & $%
\begin{array}{c}
y \\ 
z \\ 
x%
\end{array}%
$ & $%
\begin{array}{c}
x \\ 
y \\ 
z%
\end{array}%
$ & $%
\begin{array}{c}
z \\ 
y \\ 
x%
\end{array}%
$ \\ \cline{2-8}
\end{tabular}%
\medskip \newline
$g(L4)=y$ by assumption. Then consider $(n-1)$-variant profile $u1$:\medskip

\begin{tabular}{c|c|c|c|c|c|c|c|}
\cline{2-8}
$u1$ & $1$ & $2$ & $3$ & $\cdots $ & \multicolumn{1}{|c|}{$n-2$} & $n-1$ & $%
n $ \\ \cline{2-8}
& $%
\begin{array}{c}
y \\ 
x \\ 
z%
\end{array}%
$ & $%
\begin{array}{c}
y \\ 
x \\ 
z%
\end{array}%
$ & $%
\begin{array}{c}
y \\ 
z \\ 
x%
\end{array}%
$ & $\cdots $ & \multicolumn{1}{|c|}{$%
\begin{array}{c}
y \\ 
z \\ 
x%
\end{array}%
$} & $%
\begin{array}{c}
z \\ 
x \\ 
y%
\end{array}%
$ & $%
\begin{array}{c}
z \\ 
y \\ 
x%
\end{array}%
$ \\ \cline{2-8}
\end{tabular}%
\medskip \newline
If $g(u1)=x$, then $n-1$ manipulates from $L4$ to $u1$. So $g(u1)\neq x$%
.\medskip

\textbf{Step 2}\medskip

At profile $u2$:\medskip

\begin{tabular}{c|c|c|c|c|c|c|c|}
\cline{2-8}
$u2$ & $1$ & $2$ & $3$ & $\cdots $ & $n-2$ & $n-1$ & $n$ \\ \cline{2-8}
& $%
\begin{array}{c}
y \\ 
x \\ 
z%
\end{array}%
$ & $%
\begin{array}{c}
x \\ 
z \\ 
y%
\end{array}%
$ & $%
\begin{array}{c}
y \\ 
z \\ 
x%
\end{array}%
$ & $\cdots $ & $%
\begin{array}{c}
y \\ 
z \\ 
x%
\end{array}%
$ & $%
\begin{array}{c}
z \\ 
y \\ 
x%
\end{array}%
$ & $%
\begin{array}{c}
z \\ 
y \\ 
x%
\end{array}%
$ \\ \cline{2-8}
\end{tabular}%
\medskip \newline
$g(u2)=x$ since $u2\in NP^{^{\ast }}$. Then look at $(n-1)$-variant $u3$%
:\medskip

\begin{tabular}{c|c|c|c|c|c|c|c|}
\cline{2-8}
$u3$ & $1$ & $2$ & $3$ & $\cdots $ & \multicolumn{1}{|c|}{$n-2$} & $n-1$ & $%
n $ \\ \cline{2-8}
& $%
\begin{array}{c}
y \\ 
x \\ 
z%
\end{array}%
$ & $%
\begin{array}{c}
x \\ 
z \\ 
y%
\end{array}%
$ & $%
\begin{array}{c}
y \\ 
z \\ 
x%
\end{array}%
$ & $\cdots $ & \multicolumn{1}{|c|}{$%
\begin{array}{c}
y \\ 
z \\ 
x%
\end{array}%
$} & $%
\begin{array}{c}
z \\ 
x \\ 
y%
\end{array}%
$ & $%
\begin{array}{c}
z \\ 
y \\ 
x%
\end{array}%
$ \\ \cline{2-8}
\end{tabular}%
\medskip \newline
$g(u3)=x$ or $n-1$ would manipulate from $u2$ to $u3$. But $u3$ is a
2-variant of $u1$. If $g(u1)=z$, then $n-1$ would manipulate from $u3$ to $%
u1 $. So $g(u1)\neq z$. Combining with Step 1, $g(u1)=y$.\medskip

But then consider profile $u4$, an $n$-variant of $u1$,\medskip

\begin{tabular}{c|c|c|c|c|c|c|c|}
\cline{2-8}
$u4$ & $1$ & $2$ & $3$ & $\cdots $ & \multicolumn{1}{|c|}{$n-2$} & $n-1$ & $%
n $ \\ \cline{2-8}
& $%
\begin{array}{c}
y \\ 
x \\ 
z%
\end{array}%
$ & $%
\begin{array}{c}
y \\ 
x \\ 
z%
\end{array}%
$ & $%
\begin{array}{c}
y \\ 
z \\ 
x%
\end{array}%
$ & $\cdots $ & \multicolumn{1}{|c|}{$%
\begin{array}{c}
y \\ 
z \\ 
x%
\end{array}%
$} & $%
\begin{array}{c}
z \\ 
x \\ 
y%
\end{array}%
$ & $%
\begin{array}{c}
y \\ 
x \\ 
z%
\end{array}%
$ \\ \cline{2-8}
\end{tabular}%
\medskip \newline
$g(u4)=y$ or $n$ would manipulate from $u4$ to $u1$. But then at 1-variant $%
u5$:\medskip

\begin{tabular}{c|c|c|c|c|c|c|c|}
\cline{2-8}
$u5$ & $1$ & $2$ & $3$ & $\cdots $ & \multicolumn{1}{|c|}{$n-2$} & $n-1$ & $%
n $ \\ \cline{2-8}
& $%
\begin{array}{c}
z \\ 
y \\ 
x%
\end{array}%
$ & $%
\begin{array}{c}
y \\ 
x \\ 
z%
\end{array}%
$ & $%
\begin{array}{c}
y \\ 
z \\ 
x%
\end{array}%
$ & $\cdots $ & \multicolumn{1}{|c|}{$%
\begin{array}{c}
y \\ 
z \\ 
x%
\end{array}%
$} & $%
\begin{array}{c}
z \\ 
x \\ 
y%
\end{array}%
$ & $%
\begin{array}{c}
y \\ 
x \\ 
z%
\end{array}%
$ \\ \cline{2-8}
\end{tabular}%
\medskip \newline
we see that if $g(u5)=x$, then 1 would manipulate from $u5$ to $u4$. So $%
g(u5)\neq x$. But also $g(u5)\neq z$ as can be seen by considering $u6$, an $%
n$-variant of $u5$:\medskip

%
%
\medskip \newline
$g(u6)=x$ since $u6\in NP^{^{\ast }}$. If $g(u5)=z$, then $n$ would
manipulate from $u6$ to $u5$. Therefore $g(u5)\neq z$. Combining, $g(u5)=y$%
.\medskip

\textbf{Step 3}\medskip

Knowing $g(u5)=y$, we look at a 2-variant of $u5$:\medskip

%
%
\medskip \newline
$g(u7)=y$ or 2 would manipulate from $u7$ to $u5$. Then at $u8$, a 1-variant
of $u7$:\medskip

%
%
\medskip \newline
we see that if $g(u8)=z$, then 1 would manipulate from $u7$ to $u8$. So $%
g(u8)\neq z$. Next look at $(n-1)$-variant $u9$:\medskip

%
%
\medskip \newline
If $g(u9)=z$, then $n-1$ would manipulate from $u8$ to $u9$. So $g(u9)\neq z$%
.\medskip

\textbf{Final thread}\medskip

Consider profile $u10$:\medskip

%
%
\medskip \newline
$g(u10)=x$ since $u10\in NP^{^{\ast }}$. But $u10$ is an $(n-1)$-variant of $%
u9$. If $g(u9)=y$, then $n-1$ would manipulate from $u10$ to $u9$. So $%
g(u9)\neq y$. Combining with Step 3, $g(u9)=x$.\medskip

Finally, look at $L1^{\ast \ast }$, a 2-variant of $u9$:\medskip

%
%
\medskip \newline
$g(L1^{\ast \ast })=x$ or 2 would manipulate from $u9$ to $L1^{\ast \ast }$%
.\medskip

\textbf{Subsection 4-5}-\textbf{2}. Proof that $g(L3^{\ast \ast })=x$%
.\medskip

\textbf{Step 1}\medskip

We start from $g(L4)=y$ at\medskip

%
%
\medskip \newline
if $g(u1)=x$, 1 would manipulate from $u1$ to $L4$. So $g(u1)\neq x$.\medskip

\textbf{Step 2}\medskip

At profile $u2$ in $NP^{^{\ast }}$:\medskip

%
%
\medskip \newline
$g(u3)=x$ or $n-1$ would manipulate from $u2$ to $u3$. But $u3$ is a
2-variant of $u1$. If $g(u1)=z$, then 2 manipulates from $u1$ to $u3$. So $%
g(u1)\neq z$.\medskip

\textbf{Step 3}\medskip

Combining Steps 1 and 2, we have $g(u1)=y$. Look at 2-variant profile $u4$%
:\medskip

%
%
\medskip \newline
$g(u4)=y$ or 2 would manipulate from $u4$ to $u1$. Now consider $u5$, an $n$%
-variant of $u4$:\medskip

%
%
\medskip \newline
has $g(u7)=x$ or $n-1$ would manipulate from $u7$ to $u6$. So $g(u7)\neq x$.
But $u7$ is also a 1-variant of $u5$. If $g(u5)=y$, then 1 would manipulate
from $u7$ to $u5$. So $g(u5)\neq y$.\medskip

\textbf{Final thread}\medskip

Combining Steps 3 and 4, we get $g(u5)=x$. But $u5$ is a 2-variant of $%
L3^{\ast \ast }$:\medskip

%
%
\medskip \newline
$g(L3^{\ast \ast })=x$ or 2 would manipulate from $u5$ to $L3^{\ast \ast }$%
.\medskip

\textbf{Subsection 4-5}-\textbf{3}. Proof that $g(L4^{\ast \ast })=x$%
.\medskip

\textbf{Step 1}\medskip

We have already shown that $g(L4)=y$ implies $g(L3^{\ast \ast })=x$ and thus
(by Lemma 4-5), $g(L2^{\ast \ast })=x$ at\medskip

%
%
\medskip \newline
we see that if $g(u1)=z$, $n$ would manipulate from $L2^{\ast \ast }$ to $u1$%
. So $g(u1)\neq x$.\medskip

\textbf{Step 2}\medskip

Next consider profile $u2$:\medskip

%
%
\medskip \newline
$g(u2)=x$ since $u2\in NP^{^{\ast }}$. But then look at $(n-1)$-variant
profile $u3$:\medskip

%
%
\medskip \newline
$g(u3)=x$ or $n-1$ will manipulate from $u2$ to $u3$. But $u3$ is a
1-variant of $u1$. If $g(u1)=y$ then 1 would manipulate from $u3$ to $u1$.
So $g(u1)\neq y$.\medskip

\textbf{Final thread}\medskip

Combining Steps 1 and 2, $g(u1)=x$. But $u1$ is a 1-variant of $L4^{\ast
\ast }$:\medskip

%
%
\medskip \newline
$g(L4^{\ast \ast })=x$ or 1 would manipulate from $u1$ to $L4^{\ast \ast }$.
\ \ \ \ $\square $\medskip

Summarizing, we have shown:\medskip

\qquad $g$ strategy-proof on $NP\Rightarrow $ [ Range($g$) $%
=\{x,y,z\}\Rightarrow $ Range($g|NP^{\ast }$) $=\{x,y,z\}$ ]\medskip

We have actually addressed\medskip

\qquad $g$ strategy-proof on $NP\Rightarrow $ [Range($g|NP^{\ast }$) $\neq
\{x,y,z\}\Rightarrow $ Range($g$) $\neq \{x,y,z\}$ ]\medskip

Now Range($g|NP^{\ast }$) $\neq \{x,y,z\}$ can happen in two ways: the range
can contain two alternatives or one.\medskip

\qquad \textbf{Case 1}. $|$Range$(g|NP^{\ast })|=2$, say Range$(g|NP^{\ast
})=\{y,z\}$. We first find a list $\{Li\}$ of profiles such that if $x$ is
in the Range($g$), then $g(Li)=x$ for some $i$. Then we show that $g(Li)=x$
implies a violation of strategy-proofness.

\qquad \lbrack This takes up Part 1. What makes this complicated is that all
of those violations of strategy-proofness are carried out separately for
different decisiveness structures.]\medskip

\qquad \textbf{Case 2}. $|$Range$(g|NP^{\ast })|=1$, say Range$(g|NP^{\ast
})=\{x\}$. [It is not possible to do the same analysis as Case 1, because
there do exist Range-two rules on NP that have a range of one alternative on 
$NP^{\ast }$.] We first find \medskip

\qquad (1) a list $\{Li\}$ of profiles such that if $y$ is in Range$(g)$,
then $g(Li)=y$ for some $i$, and

\qquad (2) another list $\{Lj^{\ast }\}$ of profiles such that if $z$ is in
Range$(g)$, then $g(Lj^{\ast })=z$ for some $Lj^{\ast }$. \medskip

Then, we show that for each pair, $Li$, $Lj^{\ast }$, if both $g(Li)=y$ and $%
g(Lj^{\ast })=z$, then there must be a violation of strategy-proofness.

\qquad \lbrack In the details, we actually work with a related list $%
\{Lj^{\ast }\}$, and show that strategy-proofness implies $g(Lj^{\ast \ast
})=x$ and so NOT $z$.]

\bigskip

\textbf{5. M}-\textbf{Range Theorem}.\medskip

We adopt the following construction from companion paper Campbell and Kelly
(2014b). \ Let $g$ be a given strategy-proof social choice function on $%
NP(n,m+1)$ that has full range. Now we define a rule $g^{\ast }$ based on $g$%
. Select arbitrary, but distinct, $w$ and $z$ in $X$. Let $NP^{wz}(n,m+1)$
be the set of profiles in $NP(n,m+1)$ such that alternatives $w$ and $z$ are
contiguous in each individual ordering. Choose some alternative $x^{\ast }$
that does not belong to $X$ and set $X^{\ast }=\{x^{\ast }\}\cup X\backslash
\{w,z\}$. Then $g^{\ast }$ will have domain $D^{\ast }$ by which we mean the
domain $NP(n,m)$ when the feasible set is $X^{\ast }$. To define $g^{\ast }$
we begin by selecting arbitrary profile $p\in D^{\ast }$, and then we choose
some profile $r\in NP^{wz}(n,m+1)$ such that\textbf{\medskip }

\qquad 1.\qquad $r|X\backslash \{w,z\}=p|X\backslash \{w,z\}$, and

\qquad 2.\qquad for any $i\in \{1,2,...,n\}$, we have 
\[
\{x\in X\backslash \{w,z\}:x\succ _{r(i)}w\}=\{x\in X\backslash
\{w,z\}:x\succ _{p(i)}x^{\ast }\}.
\]%
In words, we create $r$ from $p$ by replacing $x^{\ast }$ with $w$ and $z$
so that $w$ and $z$ are contiguous in each $r(i)$, and $r$ does not exhibit
any Pareto domination, and in each $r(i)$ either $w$ or $z$ occupies the
same rank as $x^{\ast }$ in $p(i)$. In Campbell and Kelly (2014b), we show
that the selected alternative, which we can denote $f(p)$, is independent of
the choice of profile r and so $g^{\ast }$ is well defined.\medskip 

\textbf{Lemma 5-1}. If $g$ is strategy-proof with range $X$, then Range($%
g^{\ast }$) $=X^{\ast }$.\medskip

Before embarking on a proof of this theorem, we have to establish a lemma
about moving alternatives $w$ and $z$ closer together in a profile. \ First,
some more terminology and notation.\medskip

The \textit{position} of alternative $y\in X$ in the linear ordering $\succ $
on $X$ is the number $1+|\{y^{\prime }\in X:y^{\prime }\succ y\}|$. We say
that $a$ ranks above $b$ in $\succ $ if $a$ has a lower position number than 
$b$.\medskip

\qquad For any two distinct alternatives $a$ and $b$ in $X$ and any subset $%
Y $ of $X\backslash \{a,b\}$ we let $a\succ Y\succ b$ denote the fact that $%
a\succ y\succ b$ holds for every $y\in Y$. \medskip

\qquad Let $\sigma _{i}(r,a,b)$ denote the number of alternatives strictly
between $a$ and $b$ in $r(i)$. That is, $\sigma _{i}(r,a,b)$ is the
cardinality of the set\medskip 
\[
\{y\in X:a\succ y\succ b\text{ or }b\succ y\succ a\}. 
\]%
\medskip

\qquad Given arbitrary alternatives $a$ and $b$ and profile $r$, Lemma 5-2
exhibits a technique for moving alternatives around to reduce $\sigma
(r,a,b) $ without changing the selected alternative. Given $j\in N$ and $%
r\in NP(n,m+1)$ with $a\succ _{r(j)}b$ we say that profile $s$ is obtained
from $r $ by moving $b$ up in $r(j)$ but not above $a$ (resp., moving a down
in $r(j) $ but not below $b$) if $b$ ranks higher in $s(j)$ than in $r(j)$
(resp., a ranks lower in $s(j)$ than in $r(j)$), and $s(i)=r(i)$ for all $%
i\in N\backslash \{j\}$, and $a$ ranks above $b$ in $s(j)$. The proof of the
lemma only considers modifications to $r(j)$ that do not change the position
of any alternative above $a$ or below $b$ in $r(j)$, to ensure that the
selected alternative does not change. Such modifications will yield enough
information to allow us to prove, in the subsequent lemma, that the range of 
$g^{\ast }$ is $X^{\ast }$.\medskip

\textbf{Lemma 5-2}:\qquad Given $j\in N$, $r\in NP(n,m+1)$, and $a,b\in X$
such that $a\succ _{r(j)}b$, let $Y$ denote the set $\{y\in X:a\succ
_{r(j)}y\succ _{r(j)}b\}$.

\qquad 1. If $g(r)\notin Y$ and there exists no $u\in NP(n,m+1)$ such that $%
g(u)=g(r)$ and $u$ is obtained from $r$ by moving $b$ up in $r(j)$, and $%
\sigma _{j}(u,a,b)<\sigma _{j}(r,a,b)$ then $b\succ _{r(i)}Y$ holds for all $%
i\in N\backslash \{j\}$.\medskip

\qquad 2. If $g(r)\notin Y\cup \{a\}$ and there exists no $u\in NP(n,m+1)$
such that $g(u)=g(r)$ and $u$ is obtained from $r$ by moving $a$ down in $%
r(j)$, and $\sigma _{j}(u,a,b)<\sigma _{j}(r,a,b)$ then $Y\succ _{r(i)}a$
holds for all $i\in N\backslash \{j\}$.\medskip

\textbf{Proof:} We number the members of $Y$ so that $Y=%
\{y^{1},y^{2},...,y^{T}\}$ and $y^{t}\succ _{r(j)}y^{t+1}$ for $t\in
\{1,2,...,k-1\}$. Let $x$ denote $g(r)$. We will create profile $s$ from $r$
by switching the order of some alternatives in $r(j)$, keeping $s(i)=r(i)$
for all $i\in N\backslash \{j\}$. We will create $s(j)$ in a way that
guarantees $g(s)=g(r)=x$, provided that $s$ belongs to $NP(n,m+1)$. (In
proving 1 we will only change the ordering of the members of $Y\cup \{b\}$
relative to each other, and $b$ can only move up relative to any member of $%
Y $. Therefore, the set of alternatives preferred to $x$ by person $j$ will
not expand, it can only shrink. The fact that $x\notin Y$ makes that easy to
check. In proving 2 we will only change the ordering of the members of $%
Y\cup \{a\}$ relative to each other. We will move $a$ down relative to some
or all members of $Y$ but not in a way that changes the set of alternatives
preferred to $x$ by person $j$.)\medskip

\qquad Create profile $s$ from $r$ by switching $y^{T}$ and $b$ in $r(j)$.
If $s\in NP(n,m+1)$ then we have $g(s)=x$ and $\sigma _{j}(s,a,b)<\sigma
_{j}(r,a,b)$. If $s\notin NP(n,m+1)$ then $b\succ _{r(i)}y^{T}$ for all $%
i\in N\backslash \{j\}$. Suppose that we have $b\succ _{r(i)}y^{t}$ for $%
t=\ell ,\ell +1,...,T$ and all $i\in N\backslash \{j\}$. If $y^{k}\succ
_{r(i)}y^{\ell -1}$ holds for some $k\geq \ell $ and all $i\in N\backslash
\{j\}$ then transitivity implies that $b\succ _{r(i)}y^{t}$ holds for $%
t=\ell -1,\ell ,\ell +1,...,T$ and all $i\in N\backslash \{j\}$. Suppose
that for each $t\in \{\ell ,\ell +1,...,T\}$ we have $y^{\ell -1}\succ
_{r(i)}y^{t}$ for some $i\in N\backslash \{j\}$. Then we can create $q(j)$
from $r(j)$ by moving $y^{\ell -1}$ down until it ranks just above $b$,
keeping $r(j)$ otherwise unchanged. If $q(i)=r(i)$ for all $i\in N\backslash
\{j\}$ then we have $q\in NP(n,m+1)$ and $g(q)=x$. Now create $u$ from $q$
by switching $y^{\ell -1}$ and $b$ in $q(j)$, leaving $q$ otherwise
unchanged. If $u\in NP(n,m+1)$ then we have $g(u)=x$ and $\sigma
_{j}(u,a,b)<\sigma _{j}(r,a,b)$. If $u\notin NP(n,m+1)$ then $b\succ
_{r(i)}y^{\ell -1}$ for all $i\in N\backslash \{j\}$. If this process of
moving $y^{\ell }$ down in person $j$'s ordering for successively smaller
values of $\ell $ does not yield a profile $u\in NP(n,m+1)$ such that $%
g(u)=x $ and $\sigma _{1}(u,a,b)<\sigma _{1}(r,a,b)$ then we will have
established that\medskip

\qquad \qquad $b\succ _{r(i)}y^{t}$ for all $t\in \{1,2,...,T\}$ and all $%
i\in N\backslash \{j\}$.\medskip

\qquad If we do find the desired profile $u$ then $\{x^{\prime }\in
X:x^{\prime }\succ _{u(j)}x\}=\{x^{\prime }\in X:x^{\prime }\succ _{r(j)}x\}$
if $x\neq b$ and $\{x^{\prime }\in X:x^{\prime }\succ _{u(j)}x\}\subset
\{x^{\prime }\in X:x^{\prime }\succ _{r(j)}x\}$ if $x=b$, and hence $g(u)=x$%
.\medskip

2. Again we create profile $s$ from $r$ by changing the order of one or more
alternatives in $r(j)$, keeping $s(i)=r(i)$ for all $i\in N\backslash \{j\}$%
. Because we create $s$ by lifting a member of $Y$ just above alternative $a$
we cannot guarantee that $x=g(r)$ will still be selected unless $x\notin
Y\cup \{a\}$.\medskip

\qquad Create profile $s$ from $r$ by switching $y^{1}$ and $a$ in $r(j)$.
If $s\in NP(n,m+1)$ then we have $g(s)=x$ and $\sigma _{1}(s,a,b)<\sigma
_{1}(r,a,b)$. If $s\notin NP(n,m+1)$ then $y^{1}\succ _{r(i)}a$ for all $%
i\in N\backslash \{j\}$. We can proceed as we did in Part 1 but we must
apply the proof of Part 1 to the inverse of $r(i)$ for all $i\in N$ ---
i.e., turn the orderings of $r$ upside down --- and with alternative $a$ in
the role of alternative $b$. That is, of we let $r^{\prime }(i)$ denote the
inverse of $r(i)$ for each $i\in N$, and set $a^{\prime }=b$ and $b^{\prime
}=a$ then we can apply the proof of part 1 to $r^{\prime }$, $a^{\prime }$,
and $b^{\prime }$ by moving $b^{\prime }$ up in $r^{\prime }(j)$. But we
also have to change the names of the members of $Y$ so that\medskip

\qquad \qquad $a\succ _{r(j)}y^{T}\succ _{r(j)}y^{T-1}\succ _{r(j)}...\succ
_{r(j)}y^{2}\succ _{r(j)}y^{1}\succ _{r(j)}b$.\medskip

\qquad We will either find a profile $u\in NP(n,m+1)$ with $g(u)=x$ , $%
u(i)=r(i)$ for $i\in N\backslash \{j\}$, and $\sigma _{1}(u,a,b)<\sigma
_{1}(r,a,b)$ or else we will establish that\medskip

\qquad \qquad $y^{t}\succ _{r(i)}a$ holds for all $t\in \{1,2,...,T\}$ and
all $i\in N\backslash \{j\}$.\medskip

\qquad If we do find the desired profile $u$ then $\{x^{\prime }\in
X:x^{\prime }\succ _{u(j)}x\}=\{x^{\prime }\in X:x^{\prime }\succ _{r(j)}x\}$
and hence $g(u)=x$. \ \ \ \ $\square $\medskip

\qquad The social choice function $g^{\ast }$ derived from $g$ is defined
for two fixed alternatives $w$ and $z$. Therefore, the remaining lemma will
refer to $\sigma _{i}(p)$ instead of $\sigma _{i}(p,w,z)$. We let $\sigma
(p) $ denote the sum of the $\sigma _{i}(p)$:\medskip

\qquad \qquad $\sigma (p)=\sigma _{1}(p)+\sigma _{2}(p)+...+\sigma
_{n-1}(p)+\sigma _{n}(p)$.\medskip

\textbf{Proof of Lemma 5-1}: To establish that the range of $g^{\ast }$ is $%
X^{\ast }$ let $r$ be an arbitrary profile in $NP(n,m+1)$. It suffices to
prove that if $g(r)\in X\backslash \{w,z\}$ and $\sigma (r)>0$ there is a
profile $u\in NP(n,m+1)$ such that $g(u)=g(r)$ and $\sigma (u)<\sigma (r)$,
and if $g(r)\in \{w,z\}$ there is a profile $u\in NP(n,m+1)$ such that g(u) $%
\in $ \{w,z\} and $\sigma (u)<\sigma (r)$.\medskip

\qquad Choose $x\in X$ and some profile $r\in NP(n,m+1)$ such that $g(r)=x$.
Suppose that $\sigma (r)>0$.\medskip

\textbf{Case 1}: \qquad There exists $j\in N$ such that $\sigma _{j}(r)\geq
\sigma _{i}(r)$ for all $i\in N\backslash \{j\}$ and, $X\notin \{w,z\}$ and
for any $y\in X$, $w\succ _{r(j)}y\succ _{r(j)}z$ or $z\succ _{r(j)}y\succ
_{r(j)}w$ implies $y\neq x$.\medskip

Without loss of generality $w\succ _{r(j)}z$. Let $Y=\{y\in X:w\succ
_{r(j)}y\succ _{r(j)}z\}$. If we create $s$ from $r$ by moving $w$ down in $%
r(j)$, but not below $z$, or moving $z$ up in $r(j)$, but not above $z$ and $%
s$ belongs to $NP(n,m+1)$ then $g(s)=x$ because $x\notin Y\cup \{w,z\}$. It
follows that if we cannot find a profile $u\in NP(n,m+1)$ such that $g(u)=x$%
, $u(i)=r(i)$ for all $i\in N\backslash \{j\}$, and $\sigma (u)<\sigma (r)$
then we have\medskip

\qquad \qquad $z\succ _{r(i)}Y\succ _{r(i)}w$ for all $i\in N\backslash
\{j\} $\medskip

by Lemma 5-2. Because $\sigma _{j}(r)\geq \sigma _{i}(r)$ for all $i\in
N\backslash \{j\}$, all of the alternatives ranking between $z$ and $w$ in $%
r(i)$ must belong to $Y$ for each $i\in N\backslash \{j\}$. Choose any $h\in
N\backslash \{j\}$. Let $y^{\ast }$ be the alternative just above $w$ in $%
r(h)$. Create $u$ from $r$ by switching $y^{\ast }$ and $w$ in $r(h)$
leaving everything else unchanged. We have $u\in NP(n,m+1)$ because $n>2$.
Then $g(u)=x$ because $x$ does not belong to $Y\cup \{w,z\}$., and $\sigma
(u)<\sigma (r)$ because we have moved $w$ closer to $z$ in person $h$'s
ordering.\medskip

\textbf{Case 2}: $x\in Y$ for $Y$ defined at the beginning of Case 1. (Note
that $x\notin \{w,z\}$.)\medskip

The remainder of the proof of Lemma 5 does not actually require $\sigma
_{j}(r)\geq \sigma _{i}(r)$ for all $i\in N\backslash \{j\}$. This is
important because we will have different individuals playing the role of
person $j$. We have $w\succ _{r(j)}x\succ _{r(j)}z$. Let $A=\{a\in X:w\succ
_{r(j)}a\succ _{r(j)}x\}$ and $B=\{b\in X:x\succ _{r(j)}b\succ _{r(j)}z\}$.
Then $w\succ _{r(j)}A\succ _{r(j)}x\succ _{r(j)}B\succ _{r(j)}z$. We can
assume that we have moved $x$ up as far as possible in person $j$'s ordering
without creating Pareto dominance, and without moving it above $w$.\medskip

\textbf{Part 1}: $A\neq \varnothing $.\medskip

Then $x\succ _{r(i)}A$ for all $i\in N\backslash \{j\}$ by Lemma 5-2. If we
can move $w$ down in $r(j)$ then we can reduce $\sigma $. Otherwise $A\succ
_{r(i)}w$ for all $i\in N\backslash \{j\}$. We can reduce $\sigma $ if we
can move $z$ up in $r(j)$. Otherwise $z\succ _{r(i)}B$ for all $i\in
N\backslash \{j\}$. Therefore, we assume that\medskip

\qquad \qquad $x\succ _{r(i)}A\succ _{r(i)}w$ and $z\succ _{r(i)}B$ for all $%
i\in N\backslash \{j\}$.\medskip

Let $a^{i}$ denote the member of $A$ that ranks lowest in $r(i)$.\medskip

(I)\qquad Suppose there exists $h\in N$ such that $z\succ _{r(h)}w$ and $%
a^{h}\succ _{r(h)}w$, and no member of $X$ ranks between $a^{h}$ and $w$ in $%
r(h)$.\medskip

Clearly, $z\notin A$ and $x\succ _{r(h)}A\succ _{r(h)}w$ and hence $z\succ
_{r(h)}a^{h}\succ _{r(h)}w$. Now move $a^{h}$ just below $w$ in $r(h)$.
Alternative $x$ will be selected at the new profile, which belongs to $%
NP(n.m+1)$ because $n>2$ and $A\succ _{r(i)}w$ for all $i\in N\backslash
\{j\}$. We have thereby reduced $\sigma $.\medskip

(II)\qquad Suppose there exists $h\in N$ such that $z\succ _{r(h)}w$ and $%
a^{h}\succ _{r(h)}w$ and $z\succ _{r(h)}a^{h}\succ _{r(h)}C\succ _{r(h)}w$,
with $C=\{c\in X:a^{h}\succ _{r(h)}c\succ _{r(h)}w\}\neq \varnothing $%
.\medskip

Note that $A\cap C=\varnothing $ by definition of $a^{h}$. If we can move $w$
up in $r(h)$ above some members of $C$ then we can reduce $\sigma $ without
changing the selected alternative. ($x\succ _{r(h)}a^{h}\succ _{r(h)}C\succ
_{r(h)}w$ and thus $x\notin C$.) If we cannot reduce $\sigma $ in this
manner then $w\succ _{r(i)}C$ for all $i\in N\backslash \{h\}$ by Lemma 5-2.
Then $A\succ _{r(i)}C$ for all $i\in N\backslash \{h,j\}$. We can thus
create $r^{\prime }$ from $r$ by moving all of the members of $C$ above $%
a^{h}$ in $r(h)$ without creating Pareto domination, and without changing $%
\sigma $, provided that $r^{\prime }(h)|C=r(h)|C$. Then (I) holds if we
replace $r$ in that statement with $r^{\prime }$, and hence there exists $%
u\in NP(n,m+1)$ such that $\sigma (u)<\sigma (r^{\prime })=\sigma (r)$%
.\medskip

(III)\qquad But suppose that $z\succ _{r(h)}w$ and $x\succ _{r(h)}a^{h}\succ
_{r(h)}C^{1}\succ _{r(h)}z\succ _{r(h)}C^{2}\succ _{r(h)}w$, with $C^{1}$
(resp., $C^{2}$) containing all of the alternatives ranking between $a^{h}$
and $z$ (resp., $z$ and $w$).\medskip

Note that $C$ (from statement II) equals $C^{1}\cup C^{2}\cup \{z\}$.
Suppose that $C^{2}\neq \varnothing $. If we can move $w$ up or $z$ down in $%
r(h)$, without causing Pareto dominance and without letting $z$ rank above $%
w $, then we can reduce $\sigma $. If we cannot reduce $\sigma $ in this
matter we have $w\succ _{r(i)}C^{2}\succ _{r(i)}z$ for all i $\in $ N%
\TEXTsymbol{\backslash}\{h\} by Lemma 5-2. Clearly, $C^{2}\cap A=\varnothing 
$ , so $x\succ _{r(j)}C^{2}$ because $w\succ _{r(j)}C^{2}$ and $A$ contains
all of the alternatives ranking between $w$ and $x$ in $r(j)$. We have $x$ $%
\succ _{r(h)}$ $C^{2}$, and for all $i\in N\backslash \{h,j\}$ we
have\medskip

\qquad \qquad $x\succ _{r(i)}A\succ _{r(i)}w\succ _{r(i)}C^{2}$.\medskip

Thus, $x$ Pareto dominates the members of $C^{2}$, contradicting $r\in
NP(n,m+1)$. Therefore, $C^{2}=\varnothing $ if we cannot reduce $\sigma $ by
moving $w$ up or $z$ down in $r(h)$.\medskip

\qquad Statement III and $C^{2}=\varnothing $ imply that there exists $%
a^{h}\in A$ such that\medskip

\qquad \qquad $x\succ _{r(h)}a^{h}\succ _{r(h)}C\succ _{r(h)}z\succ _{r(h)}w$
and $z$ and $w$ are contiguous in $r(h)$, and\medskip

\qquad \qquad $C$ contains all of the alternatives ranking between $a^{h}$
and $z$.\medskip

\qquad Of course, $C\cap A=\varnothing $. Let $b^{h}$ denote the highest
ranking member of $B$ in $r(h)$. Recall that $z\succ _{r(i)}B$ for all $i\in
N\backslash \{j\}$ and hence $w\succ _{r(h)}B$ because $w$ and $z$ are
contiguous in $r(h)$ and $w\notin B$ by definition.\medskip

(IV)\qquad II and $w\succ _{r(h)}b^{h}$ both hold, and no member of $X$
ranks between $w$ and $b^{h}$ in $r(h)$,\medskip

Create profile $q$ from $r$ by setting $q(i)=r(i)$ for all $i\in N\backslash
\{j\}$ and setting $q(j)|B=r(h)^{-1}|B$, with each member of $X\backslash B$
occupying the same position in $q(j)$ as in $r(j)$. Then $\sigma (q)=\sigma
(r)$ and $q\in NP(n,m+1)$ and $b^{h}$ is the lowest ranking member of $B$ in 
$q(j)$. Obviously, $g(q)=x$ and $\sigma (q)=\sigma (r)$. Now, create $s$
from $q$ by moving $b^{h}$ just above $z$ in $q(h)$ while preserving the
position of every member of $X\backslash \{b^{h},z,w\}$ in $q(h)=r(h)$. Then 
$s\in NP(n,m+1)$ because $w\succ _{r(j)}B$, and $n>2$ and $z\succ _{r(i)}B$
for all $i\in N\backslash \{j\}$. And $g(s)=x$ because $x\succ
_{r(h)}a^{h}\succ _{r(h)}z$. Now we switch $z$ and $b^{h}$ in $s(j)=q(j)$.
This new profile belongs to $NP(n,m+1)$ and it will have a lower value of $%
\sigma $ than profile $q$.\medskip

\qquad But suppose that at least one member of $X\backslash B$ ranks between 
$w$ and $b^{h}$ in $r(h)$. Then we have\medskip

(V)\qquad $x\succ _{r(h)}a^{h}\succ _{r(h)}z\succ _{r(h)}w\succ
_{r(h)}D\succ _{r(h)}b^{h}$, and $z$ and $w$ are contiguous in $r(h)$%
.\medskip

Here $D=\{y\in X:w\succ _{r(h)}y\succ _{r(h)}b^{h}\}\neq \varnothing $ . If
we can't move $b^{h}$ above a member of $D$ without creating Pareto
domination then we have $b^{h}\succ _{r(i)}D$ for all $i\in N\backslash
\{h\} $. Alternative $b^{h}$ is the highest ranking member of $B$ in $r(h)$.
Therefore, $D\succ _{r(h)}b^{h}$ implies $D\cap B=\varnothing $. Hence $%
b^{h}\succ _{r(j)}D$ and $B\succ _{r(j)}z$ imply $z\succ _{r(j)}D$. We have $%
z\succ _{r(h)}D$ and, for all $i\in N\backslash \{h,j\}$, $z\succ _{r(i)}B$
and $b^{h}\succ _{r(i)}D$ and thus $z\succ _{r(i)}D$. Therefore, $z$ Pareto
dominates $D$, contradicting $D\neq \varnothing $. Hence, IV holds if we
can't move $b^{h}$ above a member of $D$ without creating Pareto domination.
As we have seen, this implies the existence of a profile with a lower value
of $\sigma $ than $r$ but with alternative $x$ still being selected.\medskip

\qquad If we can move $b^{h}$ above a member of $D$ then we will have
statement $V$ with a new profile in the role of $r$ and a proper subset $%
D^{\prime }$ of $D$ substituting for $D$. We then apply the argument of the
previous paragraph, eventually arriving at statement IV with a new profile $%
r^{\prime \prime }$ in place of $r$, and with $\sigma (r^{\prime \prime
})=\sigma (r)$. This implies the existence of a profile $u$ such that $%
\sigma (u)<\sigma (r)$ and $g(u)=x$.\medskip

\textbf{Part 2}: $A=\varnothing $\medskip

If $\sigma (r)>0$ and for some $i\in N$ there exists an $a\in X$ such that $%
w\succ _{r(i)}a\succ _{r(i)}x\succ _{r(i)}z$ or $z\succ _{r(i)}a\succ
_{r(i)}x\succ _{r(i)}w$ the argument of Part 1 implies that there exists a
profile in $NP(n,m+1)$ with a lower value of $\sigma $ than $\sigma (r)$.
(Part 1 assumed that $w\succ _{r(j)}a\succ _{r(j)}x\succ _{r(j)}z$ holds for
some $j\in N$ but by switching the roles of $w$ and $z$ we can also
establish the existence of a profile $u$ such that $\sigma (u)<\sigma (r)$
and $g(u)=x$ if we know that $z\succ _{r(i)}a\succ _{r(i)}x\succ _{r(i)}w$
holds for some $i\in N$.)\medskip

\qquad Let $J=\{i\in N:w\succ _{r(i)}z\}$ and $H=\{i\in N:z\succ _{r(i)}w\}$%
. Of course, $J\neq \varnothing \neq H$.\medskip

\qquad We may assume that for any $j\in J$ there is no $a\in X$ such that $%
w\succ _{r(j)}a\succ _{r(j)}x$ and for any $h\in H$ there is no $a\in X$
such that $z\succ _{r(h)}a\succ _{r(h)}x$.\medskip

\qquad If $j\in J$, let $A^{j}=\{a\in X:w\succ _{r(j)}a\succ _{r(j)}z\}$ and
if $h\in H$ and $A^{h}=\{a\in X:z\succ _{r(h)}a\succ _{r(h)}w\}$.\medskip

We have the following:\medskip

\qquad If $x\in A^{j}$ and $j\in J$ we have $w\succ _{r(j)}x\succ
_{r(j)}A^{j}\backslash \{x\}\succ _{r(i)}z$. ($w$ and $x$ are contiguous in $%
r(j)$.)\medskip

\qquad If $x\in A^{h}$ and $h\in H$ we have $z\succ _{r(h)}x\succ
_{r(h)}A^{h}\backslash \{x\}\succ _{r(h)}w$. ($z$ and $x$ are contiguous in $%
r(h)$.)\medskip

Choose any $j\in J$. Suppose that $A^{j}\neq \varnothing $ and $x\notin
A^{j} $.\medskip

\qquad If we cannot reduce $\sigma $ by moving $z$ up in $r(j)$, but not
above $w$, or $w$ down (note that $x$ would still be selected as a result)
then $z\succ _{r(i)}A^{j}\succ _{r(i)}w$ for all $i\in N\backslash \{j\}$ by
Lemma 5-2. Hence $z\succ _{r(i)}w$ for all $i\in N\backslash \{j\}$, and
thus $H=N\backslash \{j\}$ and $A^{j}\subset A^{i}$ for all $i\in H$.
Therefore, $A^{i}\neq \varnothing $ for any $i\in H$. Suppose $x\notin A^{h}$
and $h\in H$. If we cannot move $w$ up or $z$ down without changing the
selected alternative or creating Pareto domination we have $w\succ
_{r(i)}A^{h}\succ _{r(i)}z$, and thus $w\succ _{r(i)}z$, for all $i\in
N\backslash \{h\}$. This contradicts $n>2$ and $H=N\backslash \{j\}$.
Therefore,\medskip

\qquad $j\in J$ and $x\notin A^{j}\neq \varnothing $ implies $H=N\backslash
\{j\}$ and $x\in A^{h}$ for all $h\in H$.\medskip

\qquad Continuing to assume that $x\notin A^{j}\neq \varnothing $, we
have\medskip

\qquad $z\succ _{r(i)}x\succ _{r(i)}A^{i}\backslash \{x\}\succ _{r(i)}w$,
and $z$ and $x$ are contiguous in $r(i)$, for all $i\in H=N\backslash \{j\}$%
.\medskip

Choose any $h\in H$ and any $a\in A^{j}$. Then $a\neq x$. Because $%
A^{j}\subset A^{i}$ for all $i\in H$ we have $z\succ _{r(i)}a\succ _{r(i)}w$
for all $i\in H$.\medskip

Choose any two distinct $h$ and $k\in H$. If we cannot reduce $\sigma $ by
moving $w$ up in $r(k)$ without changing the selected alternative then $%
w\succ _{r(i)}A^{k}$ holds for all $i\in N\backslash \{k\}$. But $a\in A^{k}$
and thus we have $z\succ _{r(h)}a\succ _{r(h)}w\succ _{r(h)}a$,
contradicting transitivity of $r(h)$.\medskip

\qquad We are forced to conclude that for all $i\in N$, if $A^{i}\neq
\varnothing $ then $x\in A^{i}$. (If $z\succ _{r(i)}A^{i}\succ _{r(i)}w$ and 
$x\notin A^{j}\neq \varnothing $ then we also arrive at a contradiction if
we assume that we cannot reduce $\sigma $ without changing the selected
alternative.) We are assuming that $A=\varnothing $ which means that $w\succ
_{r(i)}x$ and $w$ and $x$ are contiguous in $r(i)$ for all $i\in J$, and $%
z\succ _{r(i)}x$ and $z$ and $x$ are contiguous in $r(i)$ for all $i\in H$.
For $i\in J$, let $B^{i}=\{b\in X:x\succ _{r(i)}b\succ _{r(i)}z\}$, and for $%
i\in H$ let $B^{i}=\{b\in X:x\succ _{r(i)}b\succ _{r(i)}w\}$. If $j\in J$
and we cannot reduce $\sigma $ by moving $z$ up in $r(j)$ without changing
the selected alternative then $z$ $\succ _{r(h)}$ $B^{j}$ for all $h\in
N\backslash \{j\}$. If $h\in H$ and we cannot reduce $\sigma $ by moving $w$
up in $r(h)$ without changing the selected alternative then $w\succ
_{r(j)}B^{h}$ for all $j\in N\backslash \{h\}$. Therefore\medskip

\qquad \qquad $z\succ _{r(h)}B^{j}$ for all $j\in J$ and $h\in H$,
and\medskip

\qquad \qquad $w\succ _{r(j)}B^{h}$ for all $j\in J$ and $h\in H$.\medskip

\qquad Suppose that $B^{j}\neq \varnothing $ for some $j\in J$. Let $b^{\ast
}$ denote the member of $B^{j}$ ranked lowest in $r(j)$. If $i\in J$ and $%
b^{\ast }\succ _{r(i)}w$ then $b^{\ast }\succ _{r(i)}z$ and we can switch $z$
and $b^{\ast }$ in r(j) thus reducing $\sigma $ without changing the
selected alternative or creating Pareto domination. Therefore, we may assume 
$w$ $\succ _{r(i)}$ $b^{\ast }$ for all $i\in J$. Because $w$ and $x$ are
contiguous in $r(i)$ for all $i\in J$ we have $x\succ _{r(i)}b^{\ast }$ for
all $i\in J$. Then $x\succ _{r(i)}b^{\ast }$ for all $i\in N$ because $%
z\succ _{r(h)}B^{j}$for all $h\in H$, and $z$ and $x$ are contiguous in $%
r(h) $ for all $h\in H$. Similarly, if $B^{h}\neq \varnothing $ for some $%
h\in H$ then there is an instance of Pareto domination at profile $r$%
.\medskip

\qquad Assume, then, that there exist $j\in J$ and $h\in H$ such that $%
B^{j}\cap B^{h}\neq \varnothing $. For $b\in B^{j}\cap B^{h}$and $j\in J$ we
have $w\succ _{r(j)}x\succ _{r(j)}b\succ _{r(j)}z\succ _{r(j)}b$,
contradicting transitivity. (Note that $x\succ _{r(j)}b$ holds because $j\in
J$ and $b\in B^{j}$. And $z\succ _{r(j)}b$ holds because $h\in H$ and $b\in
B^{h}$.)\medskip

\qquad Therefore, if we can't reduce $\sigma $ without changing the selected
alternative we have\medskip

\qquad \qquad $w\succ _{r(j)}x\succ _{r(j)}z$ for all $j\in J$ and $z\succ
_{r(h)}x\succ _{r(h)}w$ for all $h\in H$, and\medskip

\qquad \qquad for all $j\in J$ and all $y\in X\backslash \{w,x,z\}$, if $%
w\succ _{r(j)}y\succ _{r(j)}z$ then

$\qquad \qquad y=x$, and\medskip

\qquad \qquad for all $h\in H$ and all $y\in X\backslash \{w,x,z\}$, if $%
z\succ _{r(h)}y\succ _{r(h)}w$ then

$\qquad \qquad y=x$.\medskip

(No alternatives rank between $w$ and $x$ or between $x$ and $z$ for any $%
i\in N$.)\medskip

\qquad We now use $g$ and $r$ to define a social choice function $\mu $ with
domain $NP(n,3)$ and $X=\{w,x,z\}$. Given profile $\rho \in NP(3,3)$ let $%
p\in NP(n,m+1)$ be the profile for which, for all $i\in N$,\medskip 
\[
p(i)|\{\alpha ,\beta ,\gamma \}=\rho (i), 
\]%
\newline
and for all $y\in X\backslash \{w,x,z\}$ alternative $y$ has the same
position in $p(i)$ as $r(i)$.\medskip

Refer to $p$ as the extension of $\rho $. For any $\rho ^{\prime }\in
NP(n,3) $ set $\mu (\rho ^{\prime })=g(p^{\prime })$ for the extension $%
p^{\prime }$ of $\rho ^{\prime }$.\medskip

\qquad Suppose that $|H|$ $\geq 2$. Create $u$ from $r$ by switching $x$ and 
$z$ in $r(h)$ for some $h\in H$. We will have $u\in NP(n,m+1)$ and $g(u)=x$.
Now create s from $u$ by switching $x$ and $z$ for some $j\in J$. Create
profile $t$ from $u$ by switching $x$ and $w$ for some $h\in H$. If $g(s)=x$
or $g(t)=x$ then we have reduced $\sigma $ and hence are finished the proof
of Case 2. If $g(s)\neq x$ and $g(t)\neq x$ then strategy-proofness of g
implies that $g(s)=z$ and $g(t)=w$, in which case the range of $\mu $ is $%
\{x,z,w\}$. Then the rule $\mu $ is dictatorial because $m=3$. If at $\rho
^{\prime }$ the dictator has $x$ ranked above both $w$ and $z$ and the other
members of $N$ have the opposite ranking of the three alternatives we will
have $\mu (\rho ^{\prime })=x$. Because $g$ is strategy proof it will select 
$x$ at the extension of $\rho ^{\prime }$. Because $x$ ranks between $w$ and 
$z$ in $r(i)$ for each $i\in N$ we have reduced $\sigma $. Similarly, if $%
|J| $ $\geq 2$ we can reduce $\sigma $, without changing the selected
alternative $x$.\medskip

\textbf{Case 3}: $x\in \{w,z\}$.\medskip

To prove that $x^{\ast }$ is in the range of $g^{\ast }$ we only need to
show that a member of $\{w,z\}$ is selected by $g$ at some profile in $%
NP\ast (n,m+1)$. Let $r$ be any profile in $NP(n,m+1)$ such that $g(r)=w$%
.\medskip

\qquad Suppose that for any $i\in N$ such that $w\succ _{r(i)}z$ the
alternatives $w$ and $z$ are contiguous in $r(i)$. Then $\sigma (r)>0$
implies that there exists $h\in N$ such that $z\succ _{r(h)}C\succ _{r(h)}w$
for some nonempty subset $C$ of $X\backslash \{w,z\}$. If there is no
profile $u\in NP(n,m+1)$ such that $g(u)=w$ and $\sigma (u)<\sigma (r)$ then
we cannot create a profile $u$ from $r$ by moving $w$ up in $r(j)$ (ensuring
that $w$ will still be selected) and hence, by Lemma 5-2, we have $w\succ
_{r(i)}C$ for all $i\in N\backslash \{h\}$. Then $w\succ _{r(i)}z$ implies $%
w\succ _{r(i)}z\succ _{r(i)}C$ because $w\succ _{r(i)}C$ and $w$ and $z$ are
contiguous in $r(i)$. And $z\succ _{r(i)}w$ for $i\neq h$ implies $z\succ
_{r(i)}w\succ _{r(i)}C$ and hence $z\succ _{r(i)}C$. We also have $z\succ
_{r(h)}C$ and thus $z\succ _{r(i)}C$ for all $i\in N$, contradicting $C\neq
\varnothing $ and $r$ $\in $ $NP(n,m+1)$. We have proved that if $w\succ
_{r(i)}z$ implies that $w$ and $z$ are contiguous in $r(i)$ then $\sigma
(r)=0$.\medskip

\qquad Suppose then that there exists $j\in N$ such that $w\succ _{r(j)}z$
and $w$ and $z$ are not contiguous in $r(j)$. Let $Y$ denote the nonempty
set $\{y\in X:w\succ _{r(j)}y\succ _{r(j)}z\}$. If we cannot create a
profile $u\in NP(n,m+1)$ from $r$ by moving $z$ up in $r(j)$ but not above $%
w $ --- guaranteeing that the selected alternative does not change and $%
\sigma $ decreases --- then, by Lemma 5-2, we have $z\succ _{r(i)}Y$ for all 
$i\neq j$. If for any $i\in N$ such that $z\succ _{r(i)}w$ the alternatives $%
w$ and $z$ are contiguous in $r(i)$ then $z\succ _{r(i)}w$ implies $z\succ
_{r(i)}w\succ _{r(i)}Y$ because $z\succ _{r(i)}Y$. If $w\succ _{r(i)}z$ and $%
i\neq j$ then $w\succ _{r(i)}z\succ _{r(i)}Y$ and hence $w\succ _{r(i)}Y$.
Because $Y\neq \varnothing $ and we also have $w\succ _{r(j)}Y$ we have
contradicted the fact that $r$ exhibits no Pareto domination. Therefore,
there exists $k\in N\backslash \{j\}$ such that $z\succ _{r(k)}D\succ
_{r(k)}w$ for some nonempty subset $D$ of $X$. If there is no profile $u\in
NP(n,m+1)$ such that $g(u)=w$ and $\sigma (u)<\sigma (r)$ then we cannot
create a profile $u$ from $r$ by moving $w$ up (ensuring that $w$ will still
be selected) or $z$ down in r$(k)$ then, by Lemma 5-2, we have $w\succ
_{r(i)}D\succ _{r(i)}z$ for all $i\in N\backslash \{k\}$. In particular, $%
w\succ _{r(j)}D\succ _{r(j)}z$ and hence $D\subset Y$. But $n>2$, and for $%
i\in N\backslash \{j,k\}$ we have $w\succ _{r(i)}D\succ _{r(i)}z\succ
_{r(i)}Y$, contradicting transitivity of $r(i)$ and the fact that $D$ is a
nonempty subset of $Y$. Therefore, there must exist a profile $u\in
NP(n,m+1) $ such that $g(u)=w$ and $\sigma (u)<\sigma (r)$. $\ \ \ \ \square 
$\medskip \medskip \medskip

\textbf{References}\medskip

Barber\`{a} S, Berga D, Moreno B (2010) \textquotedblleft Individual versus
group strategy-

\qquad proofness: When do they coincide?"\ \textit{J Econ Theory} 145:
1648-1674.\medskip

Campbell DE, Kelly JS (2010a) "Losses due to the manipulation of social

\qquad choice rules," \textit{Econ Theory}, 45: 453-467.\medskip

Campbell DE, Kelly JS (2014a) Universally beneficial manipulation: a 

\qquad characterization," \textit{Soc. Choice and Welf.}, 43: 329-355.

Campbell DE, Kelly JS (2014b) "Strategy-proofness on the Non-Paretian

\qquad Subdomain".\medskip 

\bigskip 

Campbell: \ Department of Economics and The Program in Public Policy,

The College of William and Mary, Williamsburg, VA 23187-8795, USA

E-mail: decamp@wm.edu\bigskip

Kelly: Department of Economics, Syracuse University, 

Syracuse, NY 13244-1020, USA

E-mail: jskelly@maxwell.syr.edu

\end{document}